\documentclass[12pt]{article}
\usepackage[T1]{fontenc}
\usepackage{lmodern}

\usepackage{amsmath,amsthm,amscd,amssymb}
\usepackage[margin=.5 in]{geometry}
\usepackage{subfig}
\usepackage{graphicx}
\usepackage[bookmarks=false]{hyperref}
\usepackage[labelfont=bf]{caption}
\allowdisplaybreaks
\usepackage{color}
\usepackage{cancel}
\usepackage{bbm}
\usepackage{dsfont}
\usepackage{hyperref}

\usepackage{amsfonts}
\usepackage[bbgreekl]{mathbbol}
\DeclareSymbolFontAlphabet{\mathbbm}{bbold}
\DeclareSymbolFontAlphabet{\mathbb}{AMSb}%

\allowdisplaybreaks
\usepackage{mark_style}

\theoremstyle{plain}

\theoremstyle{definition}

\newtheorem{defn}{Definition}[section]
\theoremstyle{remark}
\newtheorem{rem}{Remark}[section]

\newcommand{\mc}{\mathcal}

\newcommand{\R}{\mathbb{R}}

\newcommand{\T}{{0 \leq t \leq T}}

\newcommand{\nn}{\nonumber}

\theoremstyle{plain}
\newtheorem{thm}{Theorem}[section]
\newtheorem*{thm*}{Theorem}

\title{The lifted functional approach to \\ mean field games with common noise}

\date{}

\begin{document}
\author{Mark Cerenzia, Aaron Zeff Palmer\footnote{Electronic contact: \texttt{cerenziam@math.uchicago.edu, azp6@cornell.edu}}}
\maketitle

\abstract{
We introduce a new path-by-path approach to mean field games with common noise that recovers duality at the pathwise level.  We verify this perspective by explicitly solving some difficult examples with linear-quadratic data, including control in the volatility coefficient of the common noise as well as the constraint of partial information. As an application, we establish the celebrated separation principle in the latter context.  In pursuing this program, we believe we have made a crucial contribution to clarifying the notion of regular solution in the path dependent PDE literature.
}

\tableofcontents

\section{Introduction}

This paper offers a new perspective on certain
classes of forward-backward systems
of stochastic partial differential equations 
that arise naturally in mean field game theory and the theory of optimal control with partial information.
The systems arising from either of these fields share the following major difficulty: 
although
the noise is exogenously given in the forward equation describing the state dynamics, 
the noise is endogenously determined in the backward HJB equation characterizing optimality.

We propose a novel path-by-path interpretation that exhibits duality 
between the equations of such systems at the pathwise level. 
This paper introduces and verifies this approach through significant examples, some of which
we have not yet found solved explicitly elsewhere in the literature.

Mean field games with common noise have attracted much attention
due to their practical and theoretical interest.
Indeed, it is a natural modeling assumption that all agents in a game are subject to common random shocks in addition to possible individual shocks. 
On the other hand, the problem is notoriously
difficult because the corresponding mean field game consistency condition
now features a stochastic equilibrium measure flow that must coincide with the flow of conditional laws of an optimally controlled process given the common noise.

For the PDE approach, the breakthrough work \cite{cdll1} of Cardaliaguet-Delarue-Lasry-Lions
interprets the mean field game system with common noise (see the system \eqref{sMFG} below)
as the characteristics
for the so-called \emph{master equation}, a certain PDE on Wasserstein space.
For the probabilistic approach, 
Carmona-Delarue \cite{bible2} interpret a similar class of such PDEs on Wasserstein space
as determining decoupling fields for forward-backward systems of stochastic differential equations that characterize mean field equilibria, whether for a probabilistic representation of the value function or of its gradient (the latter being the content of the Pontryagin maximum principle). 
Either of these perspectives offers ways of achieving
wellposedness for the mean field game problem in the presence of common noise, 
and further can yield explicit solutions for certain data; see
Sections 3.5 and 4.5 of Carmona-Delarue \cite{bible2} for some linear-quadratic examples featuring a common noise. 

By contrast, 
the topic of control in the volatility coefficient of the
common noise has not been explored much in the mean field game theory literature.
The only paper we have found on the topic is 
the recent work of Barrasso-Touzi \cite{barrasso2020controlled};
otherwise, some general expressions and equations in Carmona-Delarue \cite{bible1, bible2} account for the possibility
of controlled volatility coefficients, so 
the abstract theory still applies insofar as one
can characterize equilibria based on dynamic programming (leading to a system of stochastic PDEs) or based on the Pontryagin maximum principle (leading to an FBSDE). 
However, wellposedness results and explicit solutions do not seem to be available yet in the literature.

The topic of optimal control with partial information has a long history and an accordingly large literature.
We refer the reader to the book \cite{bensoussan1992stochastic} of Bensoussan
and references therein.
Mean field games with common noise and with partial information seem to be largely unexplored, 
even though the recent paper of Bensoussan-Yam \cite{benyam} that motivated our calculations clearly takes 
inspiration from these authors' own work on mean field games. 
See also the earlier paper of Bandini-Cosso-Fuhrman-Pham \cite{bandini2019randomized} that approaches the partial information problem (without mean field interactions) using viscosity solutions on Wasserstein space.
The unpublished work of Huang-Wang \cite{huang2014class} attempts to pursue this problem
via the Pontryagin maximum principle, and although we believe this probabilistic approach can work, 
the authors' calculations here do not appear to satisfy the \emph{separation principle}, a standard litmus test for such a solution. Roughly speaking, this principle says that
to go from the optimal feedback control in the case of full information to the case of partial information,
one just needs to replace the state with the best guess of the state given 
the common noise and partial observation.
A main result of this paper is that the lifted functional approach
can be used to establish this principle for mean field games with common noise and partial information; see 
the end of the 
final Section \ref{MFGpartialinfosection}
for the theorem statement and discussion.

One apparent difficulty with the dynamic programming approach to mean field games with partial information is that one must account for both the common and observational noises, so each of these must be endogenously determined
in the stochastic backward HJB equation to ensure non-anticipativity of the value function and optimal feedback control.
Another, more subtle, difficulty that arises here is that the probability measure 
with respect to which one formulates a typical player's control problem with a partial information constraint differs 
from the probability measure with respect to which one derives and articulates the forward-backward 
system of stochastic PDEs; see the system \eqref{general} below for how one may handle this issue.

Finally, if one drops the mean field coupling and partial information constraint, 
the resulting 
backward stochastic HJB equations of the various systems \eqref{sMFG}, \eqref{generalsHJB}, and \eqref{general} that we consider
are well-known to be related to so-called \emph{path dependent PDEs} (see Section 11.3.5 of Zhang \cite{zhang2017backward}).
We refer the reader to the early work 
of Ekren-Keller-Touzi-Zhang \cite{ppdefirst} and Ekren-Touzi-Zhang \cite{part1,part2}
for the first accepted notion of viscosity solution for path dependent 
PDEs, but otherwise point to
the bibliographical notes of Chapter 11 of Zhang \cite{zhang2017backward}.

\subsection*{Conceptual Overview}

On the one hand, the main concepts
in this paper were inspired by careful manipulations involving the \emph{functional It\^o formula}
for path dependent functionals (see Dupire 
\cite{dupire2019functional} and Cont-Fourni\'e 
\cite{cont2013functional}).
On the other hand, we do not know of references from the path dependent PDE literature that systematically
explore explicit solutions.
We believe this gap speaks to one of the main benefits of the lifted functional approach 
as a complementary perspective on path dependent PDEs, namely,
that it more concretely and quickly emphasizes the connection to classical PDE theory.

To explain the concept it is perhaps best to start with the \emph{path dependent heat equation}. The random solutions for which satisfy, for $t\in [0,T]$ and $\omega\in \Omega=C([0,T])$, $u(t,\omega) = \EE [ G(\bW^{t,\omega})]$, where $\bW^{t,\omega} 
	= (W^{t,\omega}_s)_\T$ is defined by
	\[
	W^{t,\omega}_s := \omega_s 1_{[0,t)}(s) + [\omega_t + W_s - W_t] 1_{[t,T]}(s).
	\]
 The path dependent PDE can be derived from the functional It\^{o} formula and expressed in terms of the vertical and horizontal derivatives (defined in Section \ref{notations} below) as
	\begin{equation}
		\label{ppdeheat}
	\begin{cases}
		-\partial_t^H u(t,\omega) - \frac{1}{2} \partial_{\omega \omega}^V u(t,\omega) = 0  &\ \omega\in \Omega,\ t\in [0,T)\\
		u(\omega) = G(\omega) &\ \omega\in \Omega_T.
	\end{cases}
	\end{equation}

Suppose next that 
the solution $u$ to \eqref{ppdeheat}
admits a lifted functional form  $\hat{u}:[0,T]\times \Omega \times \RR^d\rightarrow \RR$ such that
$$
u(t,\omega) = \hat{u}(t,\omega, \omega_t),
$$
where $(t,\omega, y) \mapsto \hat{u}(t, \omega, y)$ is strictly non-anticipative in $\omega$ (see Section \ref{notations} below for the definition).
The initial problem precipitating 
our paper is the following question:
for each fixed path $\omega \in \Omega$,
what classical parabolic PDE, if any, does $(t,y) \mapsto \hat{u}(t,\omega,y)$ satisfy?

Our understanding is 
that even experts in the path dependent PDE literature
thought this problem was already resolved by 
the work of Bion-Nadal \cite{nadal} (see the definition of ``regular solution'' in Section 2.2 therein).
Roughly speaking, the 
conclusion of \cite{nadal} was that $\hat{u}(t,\omega,y)$
should satisfy a strict non-anticipativity condition along with
the \emph{naive} classical backward heat equation 
\begin{equation}
	\label{naiveheateqn}
\begin{cases}
	-\partial_t \hat{u}(t,\omega,y) - \frac{1}{2} \partial_{yy} \hat{u}(t,\omega,y) = 0 \\
	\hat{u}(T,\omega,y) = \hat{G}(\omega, y),
\end{cases}
\end{equation}
where we assume $G(\omega)$ admits the form $G(\omega) = \hat{G}(\omega,\omega_T)$
for some lifted functional $\hat{G}(\omega,y)$.

However, \emph{this approach does not work}: the strict non-anticipativity must be enforced by an additional term in the PDE.  This omission is unfortunately 
a significant error; 
indeed, consider a simple example, e.g., the path dependent heat equation 
with terminal condition $G(\omega) = \int_0^T \omega_s \, ds$ at time $T \geq 0$. 
The correct lifted functional solution here
is well-known to be given by $\hat{u}(t,\omega, y) = \int_0^t \omega_s \, ds + (T - t) y$ (see Example 11.1.2 of Zhang \cite{zhang2017backward}), which does not satisfy equation \eqref{naiveheateqn};
instead, with the path $\omega$ fixed, it solves the completely classical PDE
\begin{equation}
\begin{cases}
	-\partial_t \hat{u}(t,\omega,y) - \frac{1}{2} \partial_{yy} \hat{u}(t,\omega,y) = y - \omega_t \\
	\hat{u}(T,\omega,y) = \hat{G}(\omega,y).
\end{cases}
\end{equation}
The additional term ``$y - \omega_t$'' admits the nice interpretation of serving to naturally enforce the strict non-anticipativity condition.

More generally, we expect 
$\hat{u}(t,\omega,y) := \EE [ G(\bW^{t,\omega,y}) ]$,
where
$\bW^{t,\omega, y} 
	= (W^{t,\omega, y}_s)_\T$ is defined by
	\begin{equation}
		\label{rightcont}
	W^{t,\omega,y}_s := \omega_s 1_{[0,t)}(s) + [y + W_s - W_t] 1_{[t,T]}(s).
	\end{equation}
	The discovery that was the impetus for this paper 
	is that $\hat{u}(t,\omega,y)$ so-defined generally
	satisfies the naive classical PDE \eqref{naiveheateqn} but
	compensated by a term enforcing strict non-anticipativity:
	\begin{equation}
		\label{compensatedheateqnexplicit}
	\begin{cases}
		-\partial_t \hat{u}(t,\omega,y) - \frac{1}{2} \partial_{yy} \hat{u}(t,\omega,y) = \mathcal{D}_\omega^y\hat{u}(t,\omega,y) \\
		\hat{u}(T,\omega,y) = \hat{G}(\omega,y),
	\end{cases}
	\end{equation}
	where the new \emph{compensator term} is defined by
		\begin{equation}
	\label{compensator}
	\calD_{\omega}^y \hat{\varphi}(t,\omega,y)  := \lim_{\epsilon \downarrow 0}  \epsilon^{-1} \left [\hat{\varphi}\big(t+\epsilon, \omega\, 1_{[0,t)} + y\, 1_{[t,T]}, y\big) - \hat{\varphi}(t+\epsilon,\omega, y) \right ], \ \ \ (t,\omega,y) \in [0,T)\times \Omega\times \RR^d.
\end{equation}
Evidently, the compensator requires evaluating at paths with a jump and we will address this in the body of our paper where we take a more rigorous approach. We will also find it is more convenient to combine the time derivative and compensator into a \emph{compensated time derivative} $D_{t, \omega}^{y} = \partial_t + \mathcal{D}_{\omega}^y$.

Our main findings are that for key examples when a path dependent PDE 
admits a lifted functional form, the underlying 
lifted functional satisfies the 
associated naive parabolic PDE 
compensated by a term enforcing strict non-anticipativity. We apply this approach to the more involved stochastic PDE resulting from mean field games with common noise and related path dependent features. Our main desire is for the lifted functional perspective to help bring important insights from the well-developed 
deterministic mean field game theory
to bear on \emph{strong} solutions for 
mean field games with common noise of various types.

\subsection*{Reader's Guide}
The main discovery of this paper extends the insight just reviewed 
for the prototypical case of the heat equation 
to the forward-backward system 
characterizing a Nash equilibrium 
for a mean field game with common noise. 
Namely, the forward-backward system of stochastic PDEs 
\begin{align} \label{sMFG}
    d u_t(x) =&\ \bigg ( - \frac{1}{2} \tr \left [(\sigma \sigma^\top +\sigma^0 (\sigma^0)^\top) \partial_{xx}^2 u_t(x) \right ] -\tr\big[ \sigma^0 \partial_x v_t(x) \big]dt
    + H\big(t,x,\partial_x u_t(x),\mu_t\big)  \bigg ) dt +  v_t(x) \cdot dW^0_t,   \nonumber \\
    d m_t(x) =&\ \bigg ( \frac{1}{2} \tr \left [ (\sigma \sigma^\top +\sigma^0 (\sigma^0)^\top) \partial_{xx}^2  m_t(x) \right ] + \text{div}_x \Big  [ m_t(x) \partial_p H\big(t,x,\partial_x u_t(x),\mu_t \big)  \Big ]  \bigg )dt - \text{div}_x \big [ m_t(x) \sigma^0   dW^0_t  \big ], \nonumber\\
    u_T(x) =&\ g(x,\mu_T), \ \  m_0(x) = \ell(x), \ \ \mu_t(dx) := m_t(x)dx.
\end{align}
can be recast\footnote{
	A general argument 
based on the seminal work of Pardoux-Peng \cite{pardoux1992backward} 
is provided in Appendix \ref{apx:derivation}.
} 
as a completely classical mean field game system of PDEs that is
appropriately compensated (analogous to \eqref{compensatedheateqnexplicit}): 
\begin{align} \label{liftedsystemintro}
-\partial_t \hat{u}(t,x,\omega,y) & - \frac{1}{2} \tr\big[(\sigma \sigma^\top +\sigma^0 (\sigma^0)^\top) \partial^2_{xx} \hat{u}(t,x,\omega,y) + 2 \sigma^0 \partial_{x}\partial_{y}\hat{u}(t,x,\omega,y) + \partial^2_{yy}\hat{u}(t,x,\omega,y)\big] \nonumber \\
  &\ + H\big(t,x, \partial_x \hat{u}(t,x,\omega,y),\hat{\mu}(t,\omega,y)\big)\nonumber   =  \calD_{\omega}^y\hat{u}(t,x,\omega,y), \nonumber  \\
		  \partial_t \hat{m}(t,x,\omega, y) & -  \frac{1}{2} \tr  \big [(\sigma \sigma^\top +\sigma^0 (\sigma^0)^\top) \partial^2_{xx} \hat{m}(t,x,\omega,y) + 2 \sigma^0 \partial_{x}\partial_{y}\hat{m}(t,x,\omega,y) + \partial^2_{yy}\hat{m}(t,x,\omega,y)\big] \nonumber \\
		    & - \text{div}_x \left [ \hat{m}(t,x,\omega,y) \partial_p H\big(t,x,\partial_x \hat{u}(t,x,\omega, y),\hat{\mu}(t,\omega,y)\big)   \right ] = -\calD_{\omega}^y\hat{m}(t,x,\omega,y), \nonumber\\
 \hat{u}(T,x,\omega,y) =&\ g(x,\hat{\mu}(T,\omega,y)), \ \ \ 
\hat{m}(0,x,\omega,y) = \ell(x-\sigma^0 y), \ \ \ \hat{\mu}(t,\omega,y)(dx) := \hat{m}(t,x,\omega,y)dx.
\end{align}
Equation (\ref{sMFG}) is a prototype forward-backward system of stochastic PDEs that characterizes a mean field game equilibrium in the presence of common noise. By lifting to a system with an additional independent variable, Equation (\ref{liftedsystemintro}) provides a novel perspective on the system that separates out the classical derivatives from the ``compensator''.

To review our program in a nutshell, we first develop technical preliminaries (Sections \ref{notations} and \ref{liftedfuntionalsection}) to make sense of the lifted functional system (\ref{liftedsystemintro}). We then show how this lifted functional approach recovers known
results in optimal control and mean field games with common noise (Sections \ref{pathdepcostprob} and \ref{MFGtypical}).
Emboldened by this consistency,
we next pursue more substantial and uncharted examples 
of mean field games with controlled common noise and partial information 
(Sections \ref{controlledvolatilitysection} and \ref{MFGpartialinfosection}, respectively).
As a sanity check, after some admittedly grueling calculations 
in Section \ref{MFGpartialinfosection}, 
we are rewarded by confirmation of the separation principle, extending its reach into new territory.

A more detailed outline of the paper is as follows.
Before we can articulate the lifted functional approach, we develop technical preliminaries and introduce notation to handle the compensator in Section \ref{notations}. We also demonstrate how the compensator relates to other pathwise derivatives studied in the literature. In particular, we restate a pathwise It\^{o} formula in Lemma \ref{lem:ito}. In Section \ref{liftedfuntionalsection}, we then cover the details of the lifted functional system \eqref{liftedsystemintro}. Here, we provide a notion of strong solution, and we prove the strong solutions of \eqref{liftedsystemintro} that satisfy a flow transformation condition are strong solutions of \eqref{sMFG} in Theorem \ref{thm:strong_solutions}.
In Section \ref{problemformulations}, 
we present in straightforward settings the problem formulations associated with the various stochastic PDE systems \eqref{sMFG}, \eqref{generalsHJB}, and \eqref{general} studied in the paper; a reader experienced in the interpretations of such systems may wish to skip this section.
In Section \ref{pathdepcostprob} we work out a first example of a mean field control problem with path dependence. This can be considered a warm-up in a simpler setting for the more involved calculations of later sections; 
nevertheless, this example also confirms the consistency of the lifted functional approach with more classical approaches of the optimal control theory literature.

Section \ref{MFGtypical} finally employs the lifted functional
approach to explicitly solve a linear quadratic mean field game with common noise; 
a reader who is pressed for time may wish to focus on Sections \ref{liftedfuntionalsection} and \ref{MFGtypical}, once acquainted with the notation of the 
\emph{compensator} \eqref{compensator} and \emph{compensated time derivative} \eqref{compensatedtimederivative} below.
However, turning to applications that constitute new results,
Sections \ref{controlledvolatilitysection} and \ref{MFGpartialinfosection} adapt the lifted functional approach
to solve mean field games with common noise featuring, respectively, control in the volatility coefficient and the constraint of a partially observed state.

\section{Notation and technical preliminaries} \label{notations}

Throughout the paper, we work on a filtered probability space 
$(\Omega', \calF, \FF, \PP)$ supporting independent standard 
$d$-dimensional Brownian motions
$\bW = (W_t)_\T$ and $\bW^0 = (W_t^0)_\T$. 
We write $\FF^{\bY} := (\calF^{\bY}_t)_\T $
with $\calF^{\bY}_t := \sigma(\cup_{0 \leq s \leq t} \ \sigma(Y_s) )$
for the filtration generated by a given stochastic process $\bY = (Y_t)_\T$. 
Finally, we write
\[
\Omega := D([0,T]; \R^d)
\]
to be the space of c\`{a}dl\`{a}g paths, whose elements serve as fixed realizations of the common noise $\bW^0 = (W^0_t)_\T$. We use the notation $1_{[a,b]}$ to denote the characteristic function of $[a,b]\subset [0,T]$, which is useful in operations between paths of $\Omega$. When differentiating vector valued functions, we use the notational convention such that, for $v:\mathbb{R}^d\rightarrow \mathbb{R}^d$ and $a\in \mathbb{R}^d$, the directional derivative is given by
$$
\frac{d}{d\epsilon}\Big|_{\epsilon=0} v(x+\epsilon\, a) = \partial_x v(x)\, a.
$$

\subsection*{Strict non-anticipativity, lifted functionals, path derivatives, compensator, and functional It\^{o} formula}

We begin by isolating the dependence of functionals $\varphi:[0,T]\times\Omega\rightarrow \mathbb{R}$ on the present value of the path through the notion of being  \emph{strictly non-anticipative}.
\begin{defn}
We 
say that $\varphi : [0,T] \times \Omega \to \RR$ is \emph{strictly non-anticipative}
if for all $t \in [0,T]$ and for all paths $\omega, \eta \in \Omega$,
$\varphi(t,\omega) = \varphi(t,\eta)$ 
whenever $\omega_s = \eta_s$ for all $0 \leq s < t$.

We say that $\varphi : [0,T] \times \Omega \to \RR$ is merely \emph{non-anticipative} 
if for all $t \in [0,T]$ and for all paths $\omega, \eta \in\Omega $,
$\varphi(t,\omega) = \varphi(t,\eta)$ 
whenever $\omega_s = \eta_s$ for all $0 \leq s \leq t$.
\end{defn}

We now consider lifted functionals that depend explicitly on the present value of the path and are strictly non-anticipative with respect to the history.
\begin{defn}\label{def:lifted}
We say that $\hat{\varphi}:[0,T]\times \Omega\times \mathbb{R}^d$ is a \emph{lifted functional} if $(t,\omega)\mapsto \hat{\varphi}(t,\omega,y)$ is strictly non-anticipative for each $y\in \RR^d$.
\end{defn}

We will also relate to the function It\^{o} calculus developed in \cite{cont2013functional} and \cite{dupire2019functional}. In these works the space considered is $\Lambda = \cup_{t\in [0,T]}\Omega_t$ where $\Omega_t = D([0,t];\RR^d)$.

\begin{lemma}\label{lem:lifting}
The functionals $\phi:\Lambda\rightarrow \RR$ are in one-to-one correspondence with merely non-anticipative functionals $\varphi:[0,T]\times \Omega \rightarrow \RR$.  Furthermore, the merely non-anticipative functionals $\varphi:[0,T]\times \Omega \rightarrow \RR$ are in one-to-one correspondence with the lifted functionals $\hat{\varphi}:[0,T]\times \Omega\times \RR^d\rightarrow \RR$.
\end{lemma}
\begin{proof}
Suppose that $\phi:\Lambda\rightarrow \RR$, then $\varphi:[0,T]\times \Omega \rightarrow \RR$ defined by
$$
\varphi(t,\omega) = \phi\big((\omega_s)_{0\leq s\leq t}\big)
$$
for $(t,\omega)\in [0,T]\times \Omega$ is non-anticipative. For the other direction, suppose that $\varphi:[0,T]\times \Omega \rightarrow \RR$ is non-anticipative. We may define $\omega^h\in \Omega_{t+h}$ by extending $\omega$ by a constant, i.e., $\omega_s^h=\omega_t$ for $s\in [t,t+h]$.  Given $\omega\in \Omega_t$ we can define $\phi:\Lambda\rightarrow \RR$ by
$$
\phi(\omega) = \varphi(t,\omega^{T-t}).
$$

Suppose again that  $\varphi : [0,T] \times \Omega \to \RR$ is non-anticipative. Then we define $\hat{\varphi}:[0,T]\times \Omega\times \RR^d\rightarrow \RR$ by
$$
\hat{\varphi}(t,\omega,y) = \varphi(t, \omega\, 1_{[0,t)} + y\, 1_{[t,T]})
$$
for $(t,\omega,y)\in  [0,T] \times \Omega \to \RR$, and easily verify that it is a lifted functional as it depends only on the strict history of $\omega$.

Suppose that  $\hat{\varphi}:[0,T]\times \Omega\times \RR^d\rightarrow \RR$ is a lifted functional then we define $\varphi:[0,T]\times \Omega\rightarrow \RR$ by
\begin{equation}\label{eqn:lifted}
\varphi(t,\omega) = \hat{\varphi}(t,\omega,\omega_t)
\end{equation}
for $(t,\omega)\in [0,T]\times \Omega$, which is easily seen to be non-anticipative.
\end{proof}

We now state the path derivatives as defined in \cite{dupire2019functional}.
\begin{defn}
Let $\phi : \Lambda \to \RR$. The \emph{horizontal derivative} of $\phi$ at time $t \in [0,T)$ and $\omega\in \Omega_t$ is written as 
	\[
	\partial_t^H \phi(\omega) := \lim_{h \downarrow 0} \frac{ \phi\big( \omega^{t+h}\big) - \phi(\omega) }{h} 
	\]
	when the limit exists.

	For each $1 \leq i \leq d$, let $e_i$ be the $i$th standard basis vector in $\RR^d$.
	The \emph{vertical derivatives} of $\phi$ are given by
	\[
	\begin{aligned}
	\partial^V_{\omega^i} \phi(\omega) & := 
	 \lim_{h \to 0} \frac{ \phi(\omega + h\, e_i \, 1_{\{t\}}) -\phi(\omega)}{h} \\
	\partial^V_{\omega^i \omega^j} \phi(\omega) & :=   \partial^V_{\omega_j}  (\partial^V_{\omega_i} \phi)(\omega)
	\end{aligned}
	\]
	when the limits exist. 
	We write $\partial^V_{\omega} \phi(\omega) := (\partial_{\omega_i} \phi(\omega))_{i=1}^d$ as a column vector
	and $\partial^V_{\omega \omega} \phi(\omega) : = (\partial_{\omega_i \omega_j} \phi(\omega))_{i,j=1}^d$ as a $d \times d$ matrix. 
\end{defn}

Using the correspondence of Lemma \ref{lem:lifting}, we may relate the path derivatives to partial derivatives of the lifted functionals along with a compensator term.
\begin{defn}
 For a lifted functional $\hat{\varphi}$,
the \emph{compensated time derivative} is defined by
\begin{equation}
	\label{compensatedtimederivative}
	 D_{t,\omega}^y  \hat{\varphi}(t,\omega,y)  := \lim_{\epsilon \downarrow 0}  \epsilon^{-1} \left [\hat{\varphi}\Big(t+\epsilon,\omega\, 1_{[0,t)} + y\, 1_{[t,T]}, y\Big) - \hat{\varphi}(t,\omega, y) \right ], \ \ \ (t,\omega, y) \in [0,T)\times \Omega \times\RR^d.
\end{equation}
It is easily seen that $D_{t,\omega}^y = \partial_t + \mathcal{D}^y_{\omega}$ for the \emph{compensator} $\mathcal{D}$ defined in \eqref{compensator} provided all limits exist.
\end{defn}

Using the notion of lifted functional form, the derivatives are given a more tractable form. These derivatives, however, turn out to be equivalent to the vertical and horizontal derivatives. The key difference in the lifted functional form is to have isolated the dependence on the present value of $\omega_t$. The \emph{strong} form of the lifted equations involves an additional independent variable for any present value $y$.
\begin{lemma}\label{lem:derivative_equivalence}
If $\phi:\Lambda\rightarrow \RR$  
has lifted functional 
form $\hat{\varphi}$ from Lemma \ref{lem:lifting} and $\partial_t^H\phi$, $\partial_\omega^V$, $\partial_{\omega\omega}^V$ exist,
then for $(t,\omega)\in [0,T]\times \Omega$ we have 
\begin{align*}
	\partial_t^H \phi\big((\omega_s)_{0\leq s\leq t}\big) =&\  
	D_{t,\omega}^y\hat{\varphi}(t,\omega, \omega_t) \quad \quad \Big[= \partial_t \hat{\varphi}(t,\omega, \omega_t)\Big], \\
	\partial_\omega^V \phi\big((\omega_s)_{0\leq s\leq t}\big)
	=&\ (\partial_y \hat{\varphi})(t,\omega, \omega_t),\\
	\partial_{\omega\omega}^V \phi\big((\omega_s)_{0\leq s\leq t}\big)
	=&\ (\partial_{yy} \hat{\varphi})(t,\omega, \omega_t),
\end{align*}
where $\partial_y$ is the 
classical spatial partial derivative in $y \in \RR^d$ (in the sense of Leibniz-Newton). [See Remark \ref{rem:first} below for how, modulo technicalities, the compensated time derivative reduces to the classical time partial derivative when evaluated at $y=\omega_t$.]
\end{lemma}
\begin{proof}
We recall that $\hat{\varphi}(t,\omega,\omega_t) = \phi(\omega)$ when $\omega \in \Omega_t$, so
\begin{align*}
\partial_t^H \phi\big((\omega_s)_{0\leq s\leq t}\big) =&\ \lim_{h \downarrow 0} \frac{ \phi\big( \omega^{t+h}\big) - \phi(\omega) }{h}\\
=&\  \lim_{h \downarrow 0} \frac{ \hat{\varphi}(t, \omega\, 1_{[0,t)} + \omega_t\, 1_{[t,T]}, \omega_t) - \hat{\varphi}(t, \omega , \omega_t) }{h}\\
=&\ D_{t,\omega}^y\hat{\varphi}(t,\omega,\omega_t).
\end{align*}
Similarly,
\begin{align*}
\partial_{\omega_i}^V \phi\big((\omega_s)_{0\leq s\leq t}\big) =&\ \lim_{h \to 0} \frac{ \phi(\omega + h\, e_i \, 1_{\{t\}}) -\phi(\omega)}{h} \\
=&\  \lim_{h \downarrow 0} \frac{ \hat{\varphi}(t, \omega, \omega_t + h\, e_i) - \hat{\varphi}(t, \omega , \omega_t) }{h}\\
=&\ \partial_{y_i}\hat{\varphi}(t,\omega,\omega_t).
\end{align*}
The second-derivatives follow the same calculation.
\end{proof}

\begin{rem}\label{rem:first}
If we  consider $\varphi(t,\omega,y) = \int_0^t \omega_s\, ds$, it is straightforward to differentiate using (\ref{compensatedtimederivative}) to find that
$$
D_{t,\omega}^y \varphi(t,\omega,y) = y.
$$
When the time derivative exists, we find that we will always have the compensator vanish when $y=\omega_t$, i.e.,
$$
\mathcal{D}_{\omega}^y \varphi(t,\omega,\omega_t)=0.
$$
In particular, in the correspondence of Lemma \ref{lem:lifting}, the horizontal time derivative corresponds to the classical partial time derivative of the lifted functional
$$
\partial_t^H \phi\big((\omega_s)_{0\leq s\leq t}\big) = D_{t,\omega}^y\hat{\varphi}(t,\omega, \omega_t) = \partial_t \hat{\varphi}(t,\omega,\omega_t).
$$

However, $\partial_t \varphi(t,\omega,y)$ does not exist when $\omega$ has a discontinuity at $t$. Other approaches are also possible, and an approach to work only with continuous paths $\omega$ and functionals that extend uniquely to a path with a single jump is reviewed in Appendix \ref{apx:unique_extensions_to_paths_with_a_single_jump}. 
\end{rem}

We now turn to the functional It\^{o} formula for the lifted functionals.
\begin{defn}
We let $\hat{\mathbb{C}}^{1,2}$ denote the lifted functionals $\hat{\varphi}$ such that $D_{t,\omega}^y \hat{\varphi}$, $\partial_y \hat{\varphi}$, and $\partial^2_{yy}\hat{\varphi}$ all exist and are along with $\hat{\varphi}$ are continuous on $[0,T]\times \Omega \times \mathbb{R}^d$ with the standard topology using the sup norm on $\Omega$. We further assume that when $\omega$ is continuous, we have
\begin{align}\label{eqn:continuous_omega}
D_{t,\omega}^y\hat{\varphi}(t,\omega, \omega_t) = \partial_t \hat{\varphi}(t,\omega,\omega_t).
\end{align}
\end{defn}
\begin{lemma}\label{lem:ito}
Suppose that $\hat{\varphi}\in \hat{\mathbb{C}}^{1,2}$. Then, given any continuous semimartingale $\bY = (Y_t)_\T$ supported on a filtered probability space $(\Omega', \calF, \FF, \PP)$, we have, $\PP$-almost surely, 
\begin{align*}
\hat{\varphi}(t,\bY,Y_t) =&\ \phi(0,\bY,Y_0) + \int_0^t \partial_t\hat{\varphi}(s,\bY,Y_s)ds\\
&\ + \int_0^t \partial_y\hat{\varphi}(s,\bY,Y_s) \cdot dY_s + \frac{1}{2}\int_0^t \partial^2_{yy}\hat{\varphi}(s,\bY,Y_s) \cdot d[Y,Y]_s.
\end{align*}
\end{lemma}
\begin{proof}
We relate this to the functional It\^{o} formula of \cite{dupire2019functional}, which is proven for functionals $\phi:\Lambda\rightarrow \mathbb{R}$ that are continuous, along with their derivatives with respect to the metric, for $\omega\in \Omega_t$ and $\eta\in \Omega_s$ with $t\leq s$,
$$
d(\omega,\eta) =|t-s| +  \sup_{r\in [0,t]} |\omega^{s-t}(r) - \eta(r)|.
$$

The It\^{o} formula then follows from the correspondence of derivatives from Lemma \ref{lem:derivative_equivalence} given that we show that $\phi:\Lambda\rightarrow \RR$ that arises from $\varphi\in \hat{\mathbb{C}}^{1,2}$ by Lemma \ref{lem:lifting} has sufficient regularity. Indeed, suppose that $\omega^i\in \Omega_{s^i}\rightarrow \omega\in \Omega_s$. This implies that $s^i\rightarrow s$ and $\omega^{i, T-s^i}$ converges to $\omega^{T-s}$ in the sup norm topology.  It follows that $\phi$ and its derivatives are $\Lambda$-continuous.  Thus, the It\^{o} lemma follows from Theorem 3.1 of \cite{dupire2019functional} in the form
\begin{align*}
\hat{\varphi}(t,\bY,Y_t) =&\ \phi(0,\bY,Y_0) + \int_0^t D_{t,\omega}^y\hat{\varphi}(s,\bY,Y_s)ds\\
&\ + \int_0^t \partial_y\hat{\varphi}(s,\bY,Y_s) \cdot dY_s + \frac{1}{2}\int_0^t \partial^2_{yy}\hat{\varphi}(s,\bY,Y_s) \cdot d[Y,Y]_s.
\end{align*}
Finally, condition \eqref{eqn:continuous_omega} and that $\bY$ is continuous implies $D_{t,\omega}^y\hat{\varphi}(s,\bY,Y_s)=\partial_t\hat{\varphi}(s,\bY,Y_s)$.
\end{proof}

As we will see below, the name derives from the interpretation that it is exactly the term to ``compensate'' the naive classical backward HJB equation to enforce strict non-anticipativity.
We remark that although $(t,\omega)\mapsto \hat{\varphi}(t,\omega,y)$ is strictly non-anticipative,
the compensated time derivative $D_{t,\omega}^y\hat{\varphi}(t,\omega,\omega_t)$ will in general extract 
the dependence on the present value $\omega_t$; indeed, one expects $D_{t,\omega}^y \hat{\varphi}(t,\omega,\omega_t) = \partial_t \hat{\varphi}(t,\omega,\omega_t)$ when $\omega$ is continuous.

\begin{rem}
To relate to more classical notions, consider an adapted c\`{a}dl\`{a}g process $\mathbf{Y}= (Y_t)_{0\leq t\leq T}$ supported on a filtered probability space $(\Omega',\mathcal{F},\mathbb{F},\mathbb{P})$. Then for $\hat{\varphi}\in \hat{\mathbb{C}}^{1,2}$,  the process $\hat{\varphi}(t,\mathbf{Y},Y_t)$ is c\~{a}dl\`{a}g and adapted. If $\mathbf{Y}$ is a continuous semi-martingale, then $\hat{\varphi}(t,\mathbf{Y},Y_t)$ is also a continous semi-martingale as follows from Lemma \ref{lem:ito}.
\end{rem}

	\begin{rem}
		The viscosity solution theory of path dependent Hamilton-Jacobi equations has been developed in \cite{zhou2023viscosity} and \cite{tang2012path}. It should be noted that these papers use the same horizontal and vertical path derivatives of \cite{dupire2019functional} and \cite{cosso2019crandall}. Our lifted functional approach departs from these works by adopting the classical partial derivatives along with the compensator, giving a new perspective on similar path dependent control and mean field game problems. In particular, we show how path dependent PDE theory can be pushed to reveal a hidden duality at the pathwise level and duality is crucial for uniqueness, stability, numerical schemes, etc. 
	\end{rem}

	\begin{rem}
		The compensated HJB method is not valid for functional data that is too sensitive to a jump \emph{nearby} a fixed time.	
		Indeed, the choice $G(\omega) = \hat{G}(\omega,\omega_T) = \sup_{0\leq  s < T} |\omega_s| \vee |\omega_T|$ provides a counterexample. 
		Although one can show $\hat{u}(t,\omega,y) := \EE[ G(\bW^{t,\omega,y}) ]$ satisfies a certain classical heat equation (see Section 3.2 of Cosso-Russo \cite{cosso2016functional} or Example 11.1.2(iii) of Zhang \cite{zhang2017backward}),
		 the uniform metric is very sensitive to a jump nearby a fixed time, so one cannot 
		compute the compensator $\calD_\omega^y \hat{u}(t,\omega,y)$. Thus, $\hat{u}(t,\omega,y)$ cannot be realized as a solution of a \emph{compensated} heat equation.

		However, Laplace's principle allows us to approximate the uniform metric as 
		\[
		\sup_{0\leq  t \leq T} |\omega_s| = \lim_{N \to \infty} N^{-1} \log \int_0^T e^{N\omega_s}ds,
		\] 
		where each approximating terminal data $G_N(\omega):= N^{-1} \log \int_0^T e^{N\omega_s}ds$ is not too sensitive to jumps. 
		One can show that the candidates $\hat{u}_N(t,\omega,y) := \EE[ G_N(\bW^{t,\omega,y})] $ are $C^2$ solutions to compensated heat equations that converge to the viscosity solution $\hat{u}(t,\omega,y) := \EE \sup_{0\leq  s \leq T} | W^{t,\omega,y}_s| $ of the path dependent heat equation with terminal data $G(\omega) = \sup_{0\leq  s \leq T} |\omega_s|$.
	\end{rem}

\section{The lifted functional approach} \label{liftedfuntionalsection}

We let $b:[0,T]\times \RR^d\times \mathcal{P}_2(\RR^d)\times \mathbb{A}\rightarrow \RR^d$ be a drift function on the control set $\mathbb{A}$, $f:[0,T]\times \RR^d\times \mathcal{P}_2(\RR^d)\times \mathbb{A}\rightarrow \RR$ be a mean field running cost, $g: \RR^d\times \mathcal{P}_2(\RR^d)\rightarrow \RR$ be a mean field terminal cost and $\sigma,\sigma^0$ be $d\times d$ matrix-valued diffusion coefficients. We define the Hamiltonian $H:[0,T]\times \RR^d\times\RR^d\times \mathcal{P}_2(\RR^d)\rightarrow \RR$ by
$$
H(t,x,p,\mu) = \sup_{\alpha \in \mathbb{A}}\big\{-p\cdot b(t,x,\mu,\alpha) - f(t,x,\mu,\alpha)\big\}.
$$
A typical mean field game system with common noise
can be stated as follows: given any probability measure $\lambda \in \calP_2(\RR^d)$
with density $\ell:\RR^d\rightarrow \RR$, find an $\FF^{\bW^0}$-adapted triple
$(u,v,m)$ of random fields on $[0,T] \times \RR^d$ satisfying the forward-backward system of
stochastic PDEs \eqref{sMFG}; see \cite{cardaliaguet2025meanfieldgamescommon} for additional analysis of the notions of solutions.

Intuitively, the forward conservation law in \eqref{sMFG}
describes how the mass density $m_t(x)$ of some agents, such as a flock of birds,
evolves in time when subject to a random environment $W^0_t$,
while the backward HJB equation
determines the value function $``u_t(x)$'' of a typical agent responding optimally
to the random evolution of the mass.
The somewhat mysterious random field $v_t(x)$ 
is part of the unknowns and
plays the role of ensuring that $u_t(x)$ is $\FF^{\bW^0}$-adapted; e.g., for a flock of birds buffeted by wind,
a typical bird at time $t$ has only observed
the behavior of the wind $(W_s^0)_{0 \leq s \leq t}$,
but is not allowed to anticipate the future 
behavior of the wind $(W_s^0)_{t < s \leq T}$ when optimizing.
A solution to the system \eqref{sMFG} can naturally be cast as a fixed
point and admits
the interpretation as characterizing a continuum version of Nash optimality.

Motivated by the literature on the functional It\^o formula, Lemma \ref{lem:ito} and
path dependent PDE theory (see Chapter 11 of Zhang \cite{zhang2017backward} and references therein), we have discovered that 
if $\hat{u}(t,x,\omega,y)$, $\hat{m}(t,x,\omega,y)$ are lifted functionals
of the $\FF^{\bW^0}$-adapted random fields $(u_t(x), m_t(x))_{0\leq t \leq T}$
from \eqref{sMFG}, then 
for each fixed path $\omega \in \Omega$,
the functions $(t,x,y) \mapsto \hat{u}(t,x,\omega,y), \hat{m}(t,x,\omega,y)$
are determined by a rough forward conservation law
coupled with a classical HJB equation that is ``compensated''
by the operator \eqref{compensator} applied to $\hat{u}(t,x,\omega,y)$.

More precisely,
the \emph{lifted functional approach
to mean field games with common noise} asserts that
the solution
of \eqref{sMFG} can be reduced
to a pair of lifted functionals 
$(\hat{u}(t,x,\omega,y), \hat{m}(t,x,\omega,y))$
satisfying the system \eqref{liftedsystemintro}.

In the backward equation of \eqref{liftedsystemintro}, the compensator term ``$\calD_\omega^y \hat{u}(t,x,\omega,y)$'' serves to enforce the strict non-anticipativity condition in 
the path variable.
However, the rather unexpected appearance of ``$\calD_\omega^y$'' in the forward equation exactly serves to exhibit 
duality
between the equations.
To our knowledge,
this duality at the pathwise level 
appears to be new and is nonobvious to illustrate otherwise.
To be more precise, 
exploiting the duality of the 
original system \eqref{sMFG}
requires taking an expectation, i.e.,
averaging over the path.

We state a theorem to make the connection between \eqref{liftedsystemintro} and \eqref{sMFG} precise. In particular, we require the forward equation to have a specific form that $\hat{m}(t,x,\omega,y)= \widetilde{r}(t,x-\sigma^0\, y,\omega)$ for some function $\tilde{r}$, or, equivalently, the flow transformation condition
\begin{align}\label{eqn:flow_transformation_condition}
-(\sigma^0)^{\top}\, \partial_x \hat{m}(t,x,\omega,y) = \partial_y \hat{m}(t,x,\omega,y).
\end{align}
This is consistent with the flow transformation method of Lions-Souganidis \cite{lions1998fully}, which reduces solving for $\hat{m}$ in the form of $\hat{m}(t,x,\omega,y)= \widetilde{r}(t,x-\sigma^0\, y,\omega)$. To illustrate we consider the case without mean field dependence, where $\widetilde{r}$ solves the random but classical PDE
\begin{equation}
	\label{forwardPDEintro}
\begin{cases}
\partial_t  \widetilde{r}(t,x,\omega)
 = \frac{1}{2} \tr  \big [ \sigma \sigma^\top  \partial_{xx}^2  \widetilde{r}(t,x,\omega) \big ]  + \text{div}_x \left [\widetilde{r}(t,x,\omega) \, \partial_p H(t,x+\sigma^0\omega_t,\partial_x \hat{u}(t,x+\sigma^0\omega_t,\omega, \omega_t))   \right ]  \\
 \widetilde{r}(0,x,\omega) = \ell(x).
\end{cases}
\end{equation}

We extend the definition of smooth lifted functionals to $\hat{\mathbb{C}}^{1,2,2}$, which now also includes twice differentiability in the $x$ variable.
\begin{thm}\label{thm:strong_solutions}
Suppose that $\hat{u}, \hat{m}\in  \hat{\mathbb{C}}^{1,2,2}$ solve (\ref{liftedsystemintro}) pointwise. Furthermore, suppose that $\hat{m}$ satisfies the flow transformation condition \eqref{eqn:flow_transformation_condition}.  
Then the triple of random fields given by
\begin{equation} \label{solutiontriple}
\big(u^*_t(x),v^*_t(x),m^*_t(x)\big):= 
	\big(\hat{u}^*(t,x,\bW^0,W^0_t), \ \partial_y \hat{u}^*(t,x,\bW^0,W^0_t), \  \hat{m}^*(t,x,\bW^0,W^0_t)\big),
\end{equation}
is a \emph{strong solution} of the original mean field game system with common noise \eqref{sMFG}.
\end{thm}
\begin{proof}
We observe from the functional It\^{o} formula of Lemma \ref{lem:ito} that, composing with the common noise process,
\begin{align} \label{keycalc2}
d u^*_t(x)  =&  \big ( \partial_t \hat{u}^* + \frac{1}{2} \partial^2_{yy} \hat{u}^*     \big )(t,x,\bW^0, W^0_t)  dt + \partial_y \hat{u}^*(t,x,\bW^0, W^0_t) \cdot dW^0_t \\
=&\ \Big( - \frac{1}{2} \tr\big[ (\sigma\sigma^\top + \sigma^0(\sigma^0)^\top) \partial_{xx}^2 u^*_t(x)  - \sigma^0 \partial_x v^*_t(x) \big] + H\big(t,x,\partial_x u^*_t(x),\mu^*_t\big) \Big) dt + v^*_t(x) \cdot dW_t^0,\nonumber
\end{align}
where we used \eqref{liftedsystemintro} in the second equality and define $\mu^*_t$ by $\mu^*_t(dx) = m^*_t(x)dx$. Similarly, 
\begin{align}\label{keycalc3}
	d   m^*_t(x) =&  \big ( \partial_t\hat{m}^*+ \frac{1}{2} \partial^2_{yy} \hat{m}^*     \big )(t,x,\bW^0, W^0_t) dt + \partial_y \hat{m}^*(t,x,\bW^0, W^0_t) \cdot dW^0_t \\
=&\ \bigg ( \frac{1}{2} \tr \left [ (\sigma\sigma^\top + \sigma^0(\sigma^0)^\top) \ \partial_{xx}^2  m_t(x) \right ] + \text{div}_x \Big  [ m_t(x) \partial_p H\big(t,x,\partial_x u_t(x),\mu^*_t \big)  \Big ]  \bigg )dt\nonumber \\
&\ + \bigg(\partial^2_{yy} \hat{m}^* (t,x,\bW^0, W^0_t) + \sigma^0\partial_x\partial_y \hat{m}^*(t,x,\bW^0, W^0_t)\bigg) dt  + \partial_y \hat{m}^*(t,x,\bW^0, W^0_t) \cdot dW^0_t. \nonumber
\end{align}

Plugging in the flow transformation condition \eqref{eqn:flow_transformation_condition}, we see that 
$$\partial_{yy}^2 \hat{m}^* (t,x,\bW^0, W^0_t) = -\sigma^0\partial_x\partial_y \hat{m}^*(t,x,\bW^0, W^0_t)$$
 and 
$$
\partial_y \hat{m}^*(t,x,\bW^0, W^0_t) \cdot dW^0_t = - (\sigma^0)^{\top}\, \partial_x m^*_t(x)\cdot dW^0_t.
$$
Since the paths $\bW^0$ are continuous we have $\partial_t \hat{u}^*(t,x,\bW^0, W^0_t) = D_{t,\omega}^y \hat{u}^*(t,x,\bW^0, W^0_t)$ and \\
$\partial_t \hat{m}^*(t,x,\bW^0, W^0_t) = D_{t,\omega}^y \hat{m}^*(t,x,\bW^0, W^0_t)$, i.e., their compensators are $0$ along the common noise. Thus by plugging in equation \eqref{liftedsystemintro}, the triple of random fields
\eqref{solutiontriple} is a strong solution of \eqref{sMFG}.
\end{proof}

We conclude this section with a couple of technical remarks relating to other approaches to solving the mean field game system of equations \eqref{sMFG}.
\begin{rem}
	For readers familiar with the notion of the master equation from Cardaliaguet-Delarue-Lasry-Lions \cite{cdll1},
	some important quantities from \eqref{liftedsystemintro}
	take on a particularly simple form.
	To see this, 
	suppose $u_t(x) = \calU(t,x,\mu_t)$ for some nice $\calU : [0,T] \times \RR^d \times \calP_2(\RR^d) \to \RR$ and that $\mu_t(dx)=m_t(x)dx$ where $m_t(x) = \hat{m}(t,x,\bW^0, W^0_t)$ for a lifted functional satisfying $\hat{m}(t,x,\omega,y) = \widetilde{r}(t,x-\sigma^0 y,\omega)$. 

	We write $$\partial_\mu \calU(t,x,\mu)(v) := \partial_v (\delta_\mu \calU)(t,x,\mu)(v),$$ where  ``$\delta_\mu$'' denotes the linear functional derivative
	and ``$\partial_\mu$'' denotes the Wasserstein gradient.
	Then, recalling $v_t(x) = \partial_y\hat{u}(t,x,\bW^0,W^0_t)$ and $\hat{\mu}(t,\omega,y)(dx) := \hat{m}(t,x,\omega,y)dx$, 
	we can compute
\begin{equation}
	\label{ibpcalc}
	\begin{aligned}
		v_t(x) = \partial_y \hat{u}(t,x,\bW^0, W^0_t) &=  
		\partial_y [\calU(t,x,\hat{\mu}(t,\bW^0,W^0_t))] \\
		& = \int_{\RR^d} (\delta_\mu \calU)(t,x,\mu_t)(v)
		\cdot (\partial_y \hat{m})(t,v,\bW^0, W^0_t) \, dv \\
		& = -\int_{\RR^d} (\delta_\mu \calU)(t,x,\mu_t)(v)
		\cdot (\sigma^0)^\top(\partial_x\widetilde{r})(t,v-\sigma^0 W^0_t,\bW^0)\, dv \\
		& = \int_{\RR^d} \sigma^0 (\partial_\mu \calU)(t,x,m_t)(v) m_t(v) \, dv,
	\end{aligned}
\end{equation}	
	where the last equality is integration by parts.
	Observe this is exactly the formula from Corollary 2.12 of Cardaliaguet-Delarue-Lasry-Lions \cite{cdll1} for the process $v_t(x)$.\footnote{The factor $\sigma^0$ is due to scaling differently than the corresponding system (31) of Cardaliaguet-Delarue-Lasry-Lions \cite{cdll1}.}
The punchline of the above is we have the formula
\begin{equation} \label{keyformula}
 \partial_y\hat{u}(t,x,\omega,y) = \int_{\RR^d} \sigma^0 (\partial_\mu \calU)(t,x,\hat{\mu}(t,\omega,y))(v) \, \hat{m}(t,v,\omega,y) \, dv.
\end{equation}

A similar calculation applies to the compensator 
of $\hat{u}(t,x,\omega, y)$.
First note, 
by comparing the forward 
equation in \eqref{liftedsystemintro} 
with the equation \eqref{forwardPDEintro}, we can identify the compensator of $\hat{m}(t,x,\omega,y)$ as
\begin{equation}
	\label{forwardcomp}
	\calD_\omega^y \hat{m}(t,x,\omega,y) = -\text{div}_x \Bigg [ \hat{m}(t,x,\omega,y)  \Big ( F(t,x,\omega,y) - F(t,x,\omega,\omega_t)  \Big)   \Bigg ]
\end{equation}
where
\[
	F(t,x,\omega,y) := \partial_x \hat{u}(t,x+\sigma^0(\omega_t - y),\omega,y).
\]
Combining \eqref{forwardcomp} with a calculation 
similar to \eqref{ibpcalc}
implies the compensator in the backward equation of \eqref{liftedsystemintro} can 
be expressed as
\begin{equation}
\calD_\omega^y \hat{u}(t,x,\omega,y) = \int_{\RR^d} (\partial_\mu \calU)(t,x,\hat{\mu}(t,\omega,y))(v) \, \hat{m}(t,v,\omega,y) \Big ( F(t,v,\omega,y) - F(t,v,\omega,\omega_t)  \Big) \, dv.
\end{equation}
The main issue with this formula is that one in general may not have access to $\partial_\mu \calU(t,x,\mu)(v)$.
\end{rem}

\begin{rem}	 \label{roughmfgremark}
	Adopting a combination of the perspectives of rough path theory and path dependent PDEs, one could introduce an alternative 
	notion of ``pathwise solution'' that consists of a pair of merely non-anticipative
	functionals $(u(t,x,\omega), m(t,x,\omega))$ on $[0,T] \times \RR^d \times \Omega$ 
	such that, for almost every (with respect to Wiener measure) $\alpha$-H\"older geometric rough path 
	$\bomega = (\omega, \bbomega)$ (i.e., $t \mapsto \omega_t$ is a fixed realization of $\bW^0$ and $(s,t) \mapsto \bbomega_{s,t}$
	is a fixed realization of the iterated Stratonovich integral $\int_s^t W^0_{s,r} \otimes \circ d W^0_r$),
	the pair of functions $(t,x) \mapsto u(t,x,\omega), m(t,x,\omega)$ satisfies
	the \emph{rough MFG system}\footnote{See, e.g., Cosso-Russo \cite{cosso2019crandall} for the definition of the vertical derivative $\partial_\omega = \partial_\omega^V$, which is simply the spatial path dependent derivative found in most any reference from the path dependent PDE literature.}
	\begin{align}
		\label{roughMFGsystem}
				d u(t,x,\omega) =&\ \bigg (- \frac{1}{2} \tr[ \big(\sigma\sigma^\top+\sigma^0(\sigma^0)^\top\big) \ \partial_{xx}^2 u(t,x,\omega)] - \tr[\sigma^0 \partial_x \partial_\omega u(t,x,\omega) ] - \frac{1}{2} \tr[\partial_{\omega \omega}^2 u(t,x,\omega)]   \bigg) dt \nonumber\\
				&   \ \ \ \ \ \ \ \ \ \ \ \ \ \ \ \ \ \ \ \ \ \ \ \ \ \ \ \ \ \  + H(t,x , \partial_x u(t,x,\omega),\mu(t,\omega) )  dt + \partial_\omega u(t,x,\omega) \cdot d\bomega_t, \nonumber\\
				 d m(t,x,\omega)  =&\ \bigg ( \frac{1}{2}\tr[\sigma \sigma^\intercal \partial_{xx}^2 m(t,x,\omega)] + \text{div}_x  [ m(t,x,\omega) \partial_p H(t,x , \partial_x u(t,x,\omega),\mu(t,\omega) ) ] \bigg ) dt \\ 
				 & \ \ \ \ \ \ \ \ \ \ \ \ \ \ \ \ \ \ \ \ \ \ \ \ \ \ \ \ \ \  - \text{div}_x  [ m(t,x,\omega) \sigma^0 d\bomega_t  ], \nonumber\\
				 u(T,x,\omega) =&\ g(x, \mu(T,\omega)), \ m(0,x,\omega) = \ell(x), \ \mu(t,\omega)(dx):=m(t,x,\omega)dx,
	\end{align}
	where, as indicated, the bold differential ``$d\bomega_t$'' can be understood in the sense of rough path theory 
	(see, e.g., Friz-Victoir \cite{friz2010multidimensional}).
	In particular, 
	the stochastic term ``$v_t(x) \cdot dW^0_t$''
	from the backward equation in \eqref{sMFG}
	corresponds to the two terms ``$- \frac{1}{2} \tr[\partial_{\omega \omega}^2 u(t,x,\omega)]  dt + \partial_\omega u(t,x,\omega) \cdot d\bomega_t$'' 
	in \eqref{roughMFGsystem}.

	Fortunately, our compensated solutions 
	of \eqref{liftedsystemintro} will furnish such an intermediate notion of pathwise solution to 
	\eqref{roughMFGsystem}
	by calculations parallel to
	\eqref{keycalc2}, \eqref{keycalc3} above, but
	based instead on a \emph{pathwise (lifted) functional It\^o formula} of the form
	\[
		d \varphi(t,\omega) = (\partial_t \hat{\varphi})(t,\omega,\omega_t) \, dt + (\partial_y \hat{\varphi})(t,\omega,\omega_t) \cdot d\bomega_t,
	\]
	given a suitable lifted functional $\hat{\varphi}(t,\omega,y)$ of $\varphi(t,\omega)$.
	This formula follows as a consequence of Keller-Zhang \cite{keller2016pathwise}, recited as (2.5)
	and (2.11)
	of Buckdahn-Keller-Ma-Zhang \cite{buckdahn2015fully}.
	
	However, getting to the point of this remark, we otherwise omit this intermediate path-by-path notion since (besides being less straightforward for calculations and estimates for a fixed point in our opinion)
	it does not exhibit that there is an underlying duality between the two equations at a pathwise level, 
	as our compensated system \eqref{liftedsystemintro} does.
	Indeed, 
	eliminating the ``$d\bomega_t$''
	term would seem to require averaging the paths over Wiener measure, thus leaving the pathwise formulation.
	\end{rem}

\section{Problem formulations} \label{problemformulations}

Now that we have reviewed the lifted functional approach in the setting of a typical mean field game with common noise,
we step back to review the various settings where we will apply the lifted functional method.
For the sake of clarity, we state these formulations somewhat informally and with straightforward data (in particular, these problems will be solved with more general data below).
More precisely, we illustrate the lifted functional approach for four problems, each of which admits an exact solution when the data fits into the framework of linear-quadratic-Gaussian control theory:
\begin{enumerate}
	\item a stochastic control problem with a path dependent terminal cost
	\item a mean field game with common noise
	\item a mean field game with controlled common noise 
	\item a mean field game with common noise and partial information
\end{enumerate}

{\bf Problem 1: }
As a warm-up, we start by considering a
stochastic control problem with a path dependent terminal cost as follows:
given an initial condition $x_0 \in \RR^d$,
\[
\text{minimize} \ \ \ \ \ \ \mathbb{E}\Big[\int_0^T \frac{1}{2}|\alpha_t|^2\, dt +  X_T\cdot \int_0^T W_s^0\, ds\Big]
\]
over $\FF^{\bW,\bW^0}$-adapted bounded continuous processes $(\alpha_t)_\T$ (here the control set is $\mathbb{A}=\mathbb{R}^d$), subject to the dynamical constraint
\[
	dX_t = \alpha_t\, dt + dW_t, \ \ \ X_0 = x_0.
\]
Intuitively, the controller will drive the process away from the anticipated random cost, $\int_0^t W_s^0\, ds$.  Indeed, we find the controller to be given as a linear feedback of $\int_0^t W_s^0\, ds$.
The explicit solution to this problem is covered in Section \ref{pathdepcostprob}. We note this is a simple example of the type of problem studied in \cite{zhou2023viscosity}, whereas our next examples go in a different direction from that work.

{\bf Problem 2: }
We consider a linear-quadratic mean field game in the spirit of Section 3.5 of Carmona-Delarue \cite{bible1}: let $s>0$, $b,\bar{b},\bar{q},\sigma,\sigma^0\in \mathbb{R}^{d\times d}$, $\mathbb{A}=\mathbb{R}^d$,
given an initial law $\lambda \in \calP_2(\RR^d)$ and a $\FF^{\bW^0}$-adapted flow of probability measures $\bmu = (\mu_t)_\T$, 
we first solve, writing $\bar{\mu}_t = \int_{\RR^d} x\, \mu_t(dx)$ for the mean position of players,
\[
\text{minimize} \ \ \ \ \ \ \mathbb{E}\Big[\int_0^T \frac{1}{2}|\alpha_t|^2\, dt + \frac{1}{2}(X_T- s\, \bar{\mu}_T)\cdot \bar{q}\, (X_T- s\, \bar{\mu}_T) \Big]
\]
over $\FF^{\bW,\bW^0}$-adapted bounded continuous processes $(\alpha_t)_\T$, subject to the dynamical constraint
\[
	dX_t = (b_t X_t+ \bar{b}_t \bar{\mu}_t + \alpha_t) \, dt + \sigma dW_t + \sigma^0\, dW_t^0, \ \ \ X_0 \sim \lambda.
\]
We denote by $(X_t^{\alpha^*})_\T$ the solution of the dynamical constraint with optimal control $(\alpha^*_t)_\T$ and second solve the fixed point problem $\mu_t = \calL(X_t^{\alpha^*}|\calF^{\bW^0}_t)$, $t \in [0,T]$, i.e., $\mu_t$ will be the conditional law of an optimally controlled
process $X_t^{\alpha^*}$ given the common noise $\bW^0 = (W^0_t)_\T$.

In this problem, the mean position of players $\bar{\mu}_t$ is translated by a Brownian common noise.  The solution we find is a linear function of the player's position and the mean position of players.
The explicit solution of this problem for a class of linear-quadratic data is covered in Section \ref{MFGtypical}.
\begin{rem}
Note this problem implicitly involves a term of the form ``$X_T \bar{\mu}_T$,''
and in turn we will see $\bar{\mu}_T$ will involve ``$\int_0^T W^0_s ds$''. Thus,
this problem features the basic structure of the path dependent cost problem
of Section \ref{pathdepcostprob}, which motivated its inclusion in this paper.
\end{rem}

{\bf Problem 3: }
We consider a similar setting as the previous problem but with a controlled volatility coefficient of the common noise, $\mathbb{A}=\mathbb{R}$:
first, given an initial law $\lambda \in \calP_2(\RR^d)$, a parameter $(\bar{a}_t)_\T$, and a flow of probability measures $\bmu = (\mu_t)_\T$, 
we first solve
\[
\text{minimize} \ \ \ \ \ \ \mathbb{E}\Big[\int_0^T \frac{1}{2}|\alpha_t - \bar{a}_t|^2\, dt + \frac{1}{2}(X_T- s\, \bar{\mu}_T)\cdot \bar{q}\, (X_T- s\, \bar{\mu}_T) \Big]
\]
over $\FF^{\bW,\bW^0}$-adapted bounded contiuous processes $(\alpha_t)_\T$, subject to the dynamical constraint
\[
dX_t = (b_t X_t + \bar{b}_t \bar{\mu}_t) \, dt + \sigma dW_t + \alpha_t \, dW_t^0, \ \ \ X_0 \sim \lambda.
\]
Second, we solve the fixed point problem $\mu_t = \calL(X_t^{\alpha^*}|\calF^{\bW^0}_t)$, $t \in [0,T]$.

The solution we find is a deterministic time dependent multiple of the parameter $\bar{a}_t$, similar to examples in the literature (see, e.g., Proposition 5.1 of Ankirchner-Fromm \cite{ankirchner2018optimal}). 
However, the factor we get reflects parameters not only from the diffusion coefficient, 
but also from
the so-called It\^{o}-Wentzell correction term, which involves the control against the unknown process ``$v_t(x)$'' that enforces the $\FF^{\bW^0}$-adaptivity constraint in the stochastic backward HJB equation in \eqref{sMFG}.
The explicit solution of this problem for a class of linear-quadratic data is covered in Section \ref{controlledvolatilitysection}.

{\bf Problem 4: }
Our final problem considers a mean field game with common noise and partial
information: let $f,g,b$ be given functions as in Section \ref{liftedfuntionalsection} with $\mathbb{A}\subset \mathbb{R}^e$, and $h:[0,T]\times \mathbb{R}^d\times \mathcal{P}_2(\mathbb{R}^d)\rightarrow \mathbb{R}^{d'}$, first, given an initial law $\lambda \in \calP_2(\RR^d)$ and a $\FF^{\bW^0}$-adapted flow of probability measures $\bmu = (\mu_t)_\T$, 
we solve
\begin{align}\label{eqn:random_cost_cost}
	\text{minimize} \ \ \ \ \ \ \mathbb{E}\Big[\int_0^T f(t,X_t,\mu_t, \alpha_t) dt + g(X_T, \mu_T)\Big],
\end{align}
subject to a dynamical constraint
\begin{align}\label{eqn:partial_information_state}
	dX_t = b(t,X_t,\mu_t, \alpha_t) dt + \sigma dW_t + \sigma^0 dW^0_t, \ \ X_0 = x_0;
\end{align}
however, there is an additional constraint that one must optimize over bounded continuous controls $\balpha = (\alpha_t)_\T$ that are adapted to $\FF^{\bW^0, \bZ}$, where $\bZ = (Z_t)_\T$ is the $d'$-dimensional \emph{observation process}
\begin{align}\label{eqn:partial_information_observations}
	dZ_t = h(t,X_t, \mu_t) dt + d\tilde{\theta}_t
\end{align}
with $\tilde{\btheta} = (\tilde{\theta}_t)_\T$ a $d'$-dimensional diffusion process with positive definite covariance $\tilde{\Theta}$ that is independent of $\bW = (W_t)_\T$ and $\bW^0 = (W_t^0)_\T$.
Second, one solves the fixed point problem $\mu_t = \calL(X_t^{\alpha^*}|\calF^{\bW^0}_t)$, $t \in [0,T]$.

Finally, we recall the mean field problem with common noise and partial information above can be interpreted as the limit of an $N$-player dynamical game: given a strategy profile $\vec{\balpha}^N = (\balpha^{N,i})_{i=1}^N$,
the $i$th player, $1 \leq i \leq N$, in the search for Nash optimality, solves the optimal control problem
\begin{align}\label{eqn:cost_N}
	\text{minimize} \ \ \ \ \ \ \mathbb{E}\Big[\int_0^T f(t,X_t^{N,i},\mu_{\bX^N_t},\beta_t)\, dt + g(X_T^{N,i},\mu_{\bX^N_T})\Big].
\end{align}
over $\FF^{\bW^0, \bZ}$-adapted bounded continuous controls $\bbeta = (\beta_t)_\T$,
subject to the dynamical constraint
\begin{align}\label{eqn:partial_information_state_N}
	dX_t^{N,k} = 
	\begin{cases}
	b(t,X_t^{N,i},\mu_{\bX^N_t} ,\beta_t)\, dt + \sigma  dW^i_t + \sigma^0 dW_t^0, & k=i, \\
	b(t,X_t^{N,k},\mu_{\bX^N_t},\alpha^{N,k}_t) \, dt + \sigma dW^k_t + \sigma^0 dW_t^0, & k\neq i,
	\end{cases}
\end{align}
and subject to the observation process
\begin{align}\label{eqn:partial_information_observations_N}
	dZ_t^i = h(t,X_t^{N,i}, \mu_{\bX^N_t})\, dt + d\tilde{\theta}_t^i, 
\end{align}
where $\mu_{\bX^N_t} := \frac{1}{N}\sum_{j=1}^N \delta_{X_t^{N,j}}$ is the empirical measure of players. 

We emphasize that players have knowledge of the common noise and their individual observation process.
Also, one can reason from this $N$-player setting that we
expect the $N \to \infty$ limit of
the empirical measures $\mu_{\bX^N_t}$ should 
converge to the conditional law of the state given the common noise $\bW^0 = (W^0_t)_\T$ with respect to $\PP$,
thus justifying the formulation made above. 

As just reviewed, the partially observed control problem is made difficult by the necessity to consider non-Markovian controls that incorporate the entirety of the history of the observation process.  
As such, the problem does not satisfy an ordinary dynamic programming principle.  With the compensated HJB equation,
a dynamic programming principle is recovered in some sense.  
Despite the mean field coupling, 
we illustrate how the solution for a linear-quadratic-Gaussian problem is still solved by \emph{the Kalman filter} and the \emph{separation principle}, as classically expected.
See Section \ref{MFGpartialinfosection}, especially equation \eqref{kalmancontrol} and nearby discussion, for more on these concepts and the explicit solution of this problem for a class of linear-quadratic data.

\section{A path dependent cost problem} \label{pathdepcostprob}

As a warm-up, we first consider a simple scenario
where there is no coupling between the
forward and backward equations of \eqref{liftedsystemintro},
which thus reduces to a classical optimal control problem.
The interest in this example
is that we can observe, in a simple setting, how our method
is consistent with the classical optimal control theory literature.
Accordingly, we first consider the solution to the path dependent cost Problem 1
reviewed in the previous Section \ref{problemformulations}.  

In the compensated HJB approach, we will solve for the lifted functional determining 
the random value function.  
The lifted value function is expected to satisfy a dynamic programming principle, i.e.,
$$
	\hat{u}(t,x,\omega,y) = \inf_{(\alpha_s)_{t\leq s\leq T}}\mathbb{E}\Big[\int_t^T \frac{1}{2}|\alpha_s|^2ds + X_T^\alpha \cdot \Big(\int_0^t \omega_s\, ds + \int_t^T[y+W^0_s - W_t^0]\, ds\Big)\Big],
$$
where
\[
dX_s^\alpha = \alpha_s\, ds+dW_s, \ \ \ X^\alpha_t = x, \ \ \ t \leq s \leq T.
\]
The compensated system of equations will have the form
\begin{align}\label{eqn:problem_1_system}
		\mathcal{D}_{\omega}^y\hat{u}(t,x,\omega,y)=&\ -\partial_t\hat{u}(t,x,\omega,y) -\frac{1}{2}\partial^2_{xx} \hat{u}(t,x,\omega,y) -\frac{1}{2}\partial^2_{yy} \hat{u}(t,x,\omega,y) + \frac{1}{2}|\partial_x \hat{u}(t,x,\omega,y)|^2, \nonumber \\
		-\mathcal{D}_{\omega}^y  \hat{m}(t,x,\omega,y)=&\ \partial_t \hat{m}(t,x,\omega,y) -\frac{1}{2}\partial^2_{xx} \hat{m}(t,x,\omega,y) -\frac{1}{2}\partial^2_{yy} \hat{m}(t,x,\omega,y) -\partial_x\big(\hat{m}(t,x,\omega,y)\, \partial_x \hat{u}(t,x,\omega,y)\big), \nonumber\\
	\hat{u}(T,x,\omega, y)=&\ x\cdot \int_0^T \omega_s\, ds, \ \ \ \hat{m}(t,x,\omega,y) = \ell(x).
\end{align}

\begin{thm}\label{thm:prob1}
A solution to (\ref{eqn:problem_1_system}), i.e., the HJB/Fokker-Planck system for Problem 1, is given by
\begin{align*}
	\hat{u}(t,x,\omega,y) =&\ -\frac{1}{2} (T-t)^3\, y^2 +  (T-t)\, x\, y - \frac{1}{8}(T-t)^4  \\ 
	& \quad -\frac{1}{2} (T-t)\Big(\int_0^t\omega_s\, ds\Big)^2 + x\, \int_0^t\omega_s\, ds -  (T-t)^2\, y\, \int_0^t\omega_s\, ds,\\
\hat{m}(t,x,\omega,y) =&\ \int_{\mathbb{R}}\ell\Big(x -\int_0^t\big[-(T-s)\omega_s -\int_0^s\omega_rdr\big]ds-z\Big)g_t(z)dz,
\end{align*}
where $g_t$ is Gaussian with mean $0$ and variance $t$.
\end{thm}

\begin{proof}
We make the ansatz
$$
	\hat{u}(t,x,\omega,y) = a_t\, x^2 + b_t\, y^2 + 2\, c_t\, x\, y + d_t + e_t\, \Big(\int_0^t \omega_s\, ds\Big)^2 + 2\, f_t\, x\, \int_0^t\omega_s\, ds + 2\, g_t\, y\, \int_0^t\omega_s\, ds.
$$
Note the terminal condition $\hat{u}(T,x,\omega, y)= \ x\cdot \int_0^T \omega_s\, ds$ is satisfied with the parameter terminal conditions 
\[
a_T = b_T = c_T = d_T = e_T = g_T = 0, \ \ f_T = \frac{1}{2}.
\]
We then compute as in Remark \ref{rem:first}, recalling the compensated time derivative of \eqref{compensatedtimederivative},
$$
	 D_{t,\omega}^y  \Big ( \int_0^t \omega_s\, ds \Big ) =  y,
$$
and plugging the ansatz into the compensated HJB equation we get
\begin{align*}
0=&\ -a_t'\, x^2 - b_t'\, y^2 - 2\, c_t'\, x\, y - d_t' - e_t'\, \big(\int_0^t \omega_s\, ds\big)^2 - 2\, f_t'\, x\, \int_0^t\omega_s\, ds - 2\, g_t'\, y\, \int_0^t\omega_s\, ds\\
&\ - 2\, e_t\, y\, \int_0^t\omega_s\, ds -2\, f_t\, x\, y - 2\, g_t\, y^2 - a_t -  b_t\\
&\ +2\, a_t^2\, x^2 + 2\, c_t^2\, y^2 + 2\, f_t^2\Big(\int_0^t\omega_s\ ds\Big)^2 +4\, a_t\, c_t\, x\, y + 4\, a_t\, f_t\, x\, \int_0^t\omega_s\, ds + 4\, c_t\, f_t\, y\, \int_0^t\omega_s\, ds.
\end{align*}

By collecting terms corresponding to $x^2, y^2, xy,\Big(\int_0^t\omega_s\, ds\Big)^2, x\, \int_0^t\omega_s\, ds, y\, \int_0^t\omega_s\, ds$, we arrive at the following system of ordinary differential equations:
\begin{itemize}
	\item $|x|^2 \ : \ a'_t = 2\, a_t^2$,
	\item $|y|^2 \ : \ b'_t = -2\, g_t + 2\, c_t^2$,
	\item $|x\, y| \ : \ c'_t = -f_t + 2\, a_t\, c_t$,
	\item $1 \ : \ d_t' = -a_t - b_t$,
	\item $\Big(\int_0^t\omega_s\, ds\Big)^2 \ : \ e_t' = 2\, f_t^2$,
	\item $x\int_0^t\omega_s\, ds \ : \ f_t' = 2\, a_t\, f_t$,
	\item $y\int_0^t\omega_s\, ds \ : \ g_t' = -e_t + 2\, c_t\, f_t$.
\end{itemize} 

We first solve $a_t=0$, thus $f_t=\frac{1}{2}$ is constant.  Now we can see that $c_t' = -\frac{1}{2}$ so $c_t=\frac{1}{2}(T-t)$, and $e_t' = \frac{1}{2}$ so $e_t = -\frac{1}{2}(T-t)$.  
We can solve for $g_t'=(T-t)$ as $g_t = -\frac{1}{2}(T-t)^2$.  
Now $b_t' = \frac{3}{2}(T-t)^2$ and $b_t = -\frac{1}{2}(T-t)^3$. We finally have that 
$d_t' = \frac{1}{2}(T-t)^3$ so $d_t=-\frac{1}{8}(T-t)^4$.
Putting everything together, we have
$$
\begin{aligned}
	\hat{u}(t,x,\omega,y) & = -\frac{1}{2} (T-t)^3\, y^2 +  (T-t)\, x\, y - \frac{1}{8}(T-t)^4  \\ 
	& \quad -\frac{1}{2} (T-t)\Big(\int_0^t\omega_s\, ds\Big)^2 + x\, \int_0^t\omega_s\, ds -  (T-t)^2\, y\, \int_0^t\omega_s\, ds.
\end{aligned}
$$

To prove the forward equation part of the theorem, we can compute that
\begin{align*}
D_{t,\omega}^y\hat{m}(t,x,\omega,y) =&\ \int_{\mathbb{R}}l'\Big(x -\int_0^t\big[-(T-s)\omega_s -\int_0^s\omega_rdr\big]ds-z\Big)g_t(z)dz\big[(T-t)y+\int_0^t\omega_sds\big]\\
&\  + \int_{\mathbb{R}}l\Big(x -\int_0^t\big[-(T-s)\omega_s -\int_0^s\omega_rdr\big]ds-z\Big)\frac{d}{dt}g_t(z)dz.
\end{align*}
We have
$$
\partial_x \hat{u}(t,x,\omega, y) =(T-t)\, y + \, \int_0^t\omega_s\, ds,
$$
so
$$
\partial_x \big(\hat{m}(t,x,\omega,y)\, \partial_x\hat{u}(t,x,\omega,y)\big) = \int_{\mathbb{R}}l'\big(x -\int_0^t\big[-(T-s)\omega_s -\int_0^s\omega_rdr\big]ds-z\big)g_t(z)dz\big[(T-t)y+\int_0^t\omega_sds\big],
$$
which coincides with the first term of $D_{t,\omega}^y\hat{m}$. We compute
$$
\frac{1}{2}\partial^2_{xx}\hat{m}(t,x,\omega,y) = \int_{\mathbb{R}}l\big(x -\int_0^t\big[-(T-s)\omega_s -\int_0^s\omega_rdr\big]ds-z\big)g''_t(z)dz,
$$
which coincides with the second term of $D_{t,\omega}^y\hat{m}$ as $g_t$ solves the heat equation. Finally, $\partial^2_{yy}\hat{m} = 0$, completing the proof.
\end{proof}

\begin{rem}
The optimal $\FF^{\bW^0}$-adapted feedback control $\hat{\balpha} = (\alpha^*_t)_\T$ for Problem 1 is given by
$$
	\alpha^*_t = -\partial_x \hat{u}(t,x,\mathbf{W}^0, W_t^0) =-(T-t)W_t^0 - \, \int_0^tW_s^0\, ds.
$$
Further, the optimal expected value at time zero is $u(0,x_0,\omega,0) = d_0 =-\frac{1}{8}T^4$, 
which is notably independent of the initial position $X_0=x_0$.
\end{rem}

\subsection*{Comparison with the literature}
	A more classical approach to the path dependent cost problem might be to make the problem Markovian by introducing the new state
variables $Y_t:=W_t^0$ and $\Xi_t := \int_0^t Y_s\, ds$.  
In these variables, the problem turns into a stochastic control problem with value function $v(t,x,y,\xi)$ solving the degenerate HJB equation
\[
\begin{cases}
	-\partial_t v(t,x,y,\xi) -\frac{1}{2}\partial^2_{xx} v(t,x,y,\xi) -\frac{1}{2}\partial^2_{yy} v(t,x,y,\xi) - y\cdot \partial_\xi v(t,x,y,\xi) + \frac{1}{2}|\partial_x v(t,x,y,\xi)|^2= 0\\
	v(T,x,y,\xi)=\  x\cdot \xi.
\end{cases}
\]
Observe the correspondence between this approach
with the lifted functional approach is
$$
	\hat{u}(t,x,\omega,y) = v\Big(t,x,y, \int_0^t\omega_s\, ds\Big).
$$
Then we can note that the compensated time derivative satisfies
$$
	D_{t,\omega}^y \hat{u}(t,x,\omega,y) = \partial_t v\Big(t,x,y, \int_0^t\omega_s\, ds\Big)+ y\cdot \partial_\xi v\Big(t,x,y, \int_0^t\omega_s\, ds\Big),
$$
establishing consistency between the two approaches. 

We remark, however, that this more classical reasoning does not seem to work in general for 
	the other more complicated problems we study. Indeed, the desired structure to make the problem Markovian as above cannot be easily determined in advance. 
	
Finally, given the lifted functional approach was motivated
by concepts from the literatures on the functional It\^o formula
and path dependent PDE theory,
we mention that there is a path dependent PDE that 
the functional $u(t,x,\omega) := \hat{u}(t,x,\omega, \omega_t)$
will satisfy that one may work with instead to arrive at the same solution.
Again, we refer the reader to Chapter 11 of Zhang \cite{zhang2017backward}.

\section{Mean field game with common noise} \label{MFGtypical}
\subsection*{The linear-quadratic data for the MFG problem}

Let us recall the linear-quadratic data 
from Problem 2 in Section \ref{problemformulations}:
writing $\bar{\mu}:= \int_{\RR^d} \xi\, \mu(d\xi)$, we set
\begin{equation} \label{LQdataMFG}
\begin{aligned}
b(t,x,\mu,\alpha) & : = b_t\, x+\bar{b}_t\, \bar{\mu} + \alpha, \\
f(t,x,\mu,\alpha) & := \frac{1}{2}  |\alpha|^2, \\
g(x,\mu) & := \frac{1}{2} (x - s\, \bar{\mu})\cdot \bar{q}\, (x - s\, \bar{\mu}),
\end{aligned}
\end{equation}
where we refer to Section \ref{notations} for the description of these given parameters. The Hamiltonian is then given by
$$
H(t,x,p,\mu) = -p\cdot(b_t\, x+\bar{b}_t\, \bar{\mu}) + \frac{1}{2}|p|^2 .
$$

The compensated HJB equation of (\ref{liftedsystemintro}) takes the form
\begin{align*}
	 \mathcal{D}_{\omega}^y  \hat{u}(t,x,\omega,y)=&\ -  \partial_t \hat{u}(t,x,\omega,y) + \frac{1}{2} | \partial_x \hat{u}(t,x,\omega,y) |^2
	 - \partial_x \hat{u}(t,x,\omega,y) \cdot \big ( b_t x + \bar{b}_t \bar{\mu}(t,\omega, y) \big)\\
	 &\ -\frac{1}{2} \tr\big[(\sigma\sigma^\top + \sigma^0(\sigma^0)^\top) \ \partial_{xx}^2 \hat{u}(t,x,\omega,y)  + \partial^2_{yy} \hat{u}(t,x,\omega,y)
	+ 2 \sigma^0 \partial_x \partial_y \hat{u}(t,x,\omega,y) \big] 
\end{align*}
with terminal condition
\[
\hat{u}(T,x,\omega,y) =  \frac{1}{2} \big(x - s\, \bar{\mu}(T,\omega, y)\big)\cdot \bar{q} \big(x - s\, \bar{\mu}(T,\omega, y)\big).
\]

\begin{thm}\label{thm:prob2}
Suppose that the following system of ODE's for $(\Gamma_t,\Gamma_t^0,\Lambda_t^0,\Delta_t)$ has a solution on $[0,T]$:
\begin{align}\label{eqn:prob2_riccati}
	\dot{\Gamma}_t =&\ \Gamma_t^\top \Gamma_t - \Gamma_t^\top b_t - b_t^\top \Gamma_t  , & \Gamma_T = \bar{q},\nn\\
	\dot{\Lambda}_t^0 =&\ (\Lambda^0_t)^\top \Lambda^0_t -(\Lambda^0_t)^\top\, (b_t + \bar{b}_t - \Gamma_t) + \Gamma_t^\top  \Lambda^0_t  - \Gamma_t^\top\, \bar{b}_t  - b_t^\top\, \Lambda_t^0,  & \Lambda^0_T = - \bar{q}\, s,\nn\\
	 \dot{\Gamma}^0_t =&\  -(b_t + \bar{b}_t - \Gamma_t - \Lambda^0_t)^\top \Gamma^0_t -(\Gamma^0_t)^\top(b_t + \bar{b}_t - \Gamma_t - \Lambda^0_t) &\nn\\
	& \quad \quad \quad \ \ +  (\Lambda^0_t)^\top \Lambda^0_t - (\Lambda^0_t)^\top \bar{b}_t - \bar{b}_t^\top \Lambda^0_t,& \Gamma^0_T = s^\top \bar{q}\, s,\nn\\
	\dot{\Delta}_t =&\  - \frac{1}{2}\tr[(\sigma \sigma^\top + \sigma^0 (\sigma^0)^\top) \ \Gamma_t ] - \frac{1}{2} \tr[  \sigma^0 (\sigma^0)^\top \Gamma^0_t] 
	- \tr[\sigma^0 (\sigma^0)^\top \Lambda^0_t], & \Delta_T = 0.
\end{align} 
Then let
\[
\bar{\mu}(t,\omega, y)
= \Phi_t  \bigg (\bar{\mu}_0 
+ \int_0^t \Phi_s^{-1} (b_s + \bar{b}_s - \Gamma_s - \Lambda^0_s)\, \sigma^0\, \omega_s\, ds \bigg )
+ \sigma^0\, y, 
\]
and a solution to (\ref{liftedsystemintro}) for Problem 2, is given by
\begin{align*}
\hat{u}(t,x,\omega,y) =&\  
\frac{1}{2} \left ( x \cdot \Gamma_t\, x + \bar{\mu}(t,\omega, y) \cdot \Gamma^0_t\, \bar{\mu}(t,\omega, y) +
2x \cdot \Lambda^0_t\, \bar{\mu}(t,\omega, y) \right )
+ \Delta_t,
\end{align*}
and $\hat{m}(t,x,\omega,y) = g\big(t,x-\bar{\mu}(t,\omega,y)\big)$, where 
$g$ solves $g(0,\cdot) = \ell$ and 
\begin{align*}
\partial_tg(t,z) =  \frac{1}{2}{\rm tr}\big[\sigma \sigma^\top \partial^2_{zz} g(t,z)\big] + {\rm div}_z \big[ g(t,z)(\Gamma_t-b_t) z\big].
\end{align*}
\end{thm}

\begin{proof}
We make the ansatz that the solution of (\ref{sMFG}) has the form
\[
u_t(x) = 
\frac{1}{2} \left ( x \cdot \Gamma_t\, x + \bar{\mu}_t \cdot \Gamma^0_t\, \bar{\mu}_t +
2x \cdot \Lambda^0_t\, \bar{\mu}_t \right )
+ \Delta_t
\]
so that the optimal feedback function is given by
\begin{equation}
\label{feedbackfullinfo}
\begin{aligned}
- \partial_x u_t(x)
= - \Gamma_t\, x - \Lambda^0_t\, \bar{\mu}_t.
\end{aligned}
\end{equation}
Hence, we have
\[
dX_t = \left ( (b_t - \Gamma_t)\,  X_t +(\bar{b}_t - \Lambda^0_t )\, \bar{\mu}_t \right ) dt + \sigma\, dW_t + \sigma^0\, dW^0_t,
\]
and taking expectations of this equation conditional on $\calF^{\bW^0}_t$
yields
\[
d\bar{\mu}_t = \left ( b_t + \bar{b}_t - \Gamma_t - \Lambda^0_t \right ) \bar{\mu}_t\, dt + \sigma^0\, dW^0_t,
\]
which has an explicit solution
of the form $\bar{\mu}_t = \bar{\mu}(t,\bW^0, W^0_t)$
where
\[
\bar{\mu}(t,\omega, y)
= \Phi_t  \bigg (\bar{\mu}_0 
+ \int_0^t \Phi_s^{-1} (b_s + \bar{b}_s - \Gamma_s - \Lambda^0_s)\, \sigma^0\, \omega_s\, ds \bigg )
+ \sigma^0\, y, 
\]
where $(\Phi_t)_\T$ is the solution of the matrix-valued ODE
\[
\dot{\Phi}_t = (b_t + \bar{b}_t -\Gamma_t - \Lambda^0_t)\, \Phi_t, \ \ \Phi_0 = 1.
\]
Thus, the ansatz for the lifted value function becomes
\[
\hat{u}(t,x,\omega,y) = 
\frac{1}{2} \left ( x \cdot \Gamma_t\, x + \bar{\mu}(t,\omega, y) \cdot \Gamma^0_t\, \bar{\mu}(t,\omega, y) +
2x \cdot \Lambda^0_t\, \bar{\mu}(t,\omega, y) \right )
+ \Delta_t.
\]

Now we may begin computing the terms appearing in the lifted functional backward equation \eqref{liftedsystemintro}. We first compute
\[
\begin{aligned}
 D_{t,\omega}^y  \bar{\mu}(t,\omega, y)
& = 
(b_t + \bar{b}_t - \Gamma_t - \Lambda^0_t)\, \bar{\mu}(t,\omega, y).
\end{aligned}
\]
\[
\partial_y \bar{\mu}(t,\omega, y)
 =  \sigma^0 
\]
Then we have
\[
\begin{aligned}
 D_{t,\omega}^y  \hat{u}(t,x,\omega,y) 
& = \frac{1}{2} \left ( x \cdot \dot{\Gamma}_t\, x + 2x \cdot \dot{\Lambda}^0_t\, \bar{\mu}(t,\omega, y) + \bar{\mu}(t,\omega, y) \cdot \dot{\Gamma}^0_t\, \bar{\mu}(t,\omega, y) \right )
+ \dot{\Delta}_t \\
& \quad \ +  (b_t + \bar{b}_t - \Gamma_t - \Lambda^0_t)\, \bar{\mu}(t,\omega, y) \cdot \left ( \Gamma^0_t\, \bar{\mu}(t,\omega, y) 
 + \Lambda^0_t\, x  \right ),
\end{aligned}
\]
and can further compute
\[
\partial_x \hat{u}(t,x,\omega,y)  = 
\Gamma_t\, x 
+ \Lambda^0_t\, \bar{\mu}(t,\omega,y), \ \ \ 
\partial_{xx}^2 \hat{u}(t,x,\omega,y)  = 
\Gamma_t,
\]
\[
\partial_y \hat{u}(t,x,\omega,y) 
= (\sigma^0)^{\top}  \left ( \Gamma^0_t\,	 \bar{\mu}(t,\omega,y)
+ \Lambda^0_t\, x  \right ), \ \ \ 
\partial_{yy}^2 \hat{u}(t,x,\omega,y) 
= (\sigma^0)^\top \Gamma^0_t\, \sigma^0,
\]
\[
\partial_y \partial_x\hat{u}(t,x,\omega,y)
=
\Lambda^0_t\sigma^0.
\]

Inputting the above calculations in the compensated equation gives
\[
\begin{aligned}
	 0=&\ - \frac{1}{2} \left ( x \cdot \dot{\Gamma}_t\, x + \bar{\mu}(t,\omega, y) \cdot \dot{\Gamma}^0_t\, \bar{\mu}(t,\omega, y) +
2x \cdot \dot{\Lambda}^0_t\, \bar{\mu}(t,\omega, y) \right )
- \dot{\Delta}_t \\
&\  -(b_t + \bar{b}_t - \Gamma_t - \Lambda^0_t)\, \bar{\mu}(t,\omega, y) \cdot \left ( \Gamma^0_t\, \bar{\mu}(t,\omega, y) 
 + \Lambda^0_t\, x  \right )    \\
 &\  - \frac{1}{2} \tr\big[ (\sigma\sigma^\top + \sigma^0(\sigma^0)^\top) \Gamma_t  \sigma^0\, (\sigma^0)^\top \Gamma^0_t
	+ 2 \sigma^0\, (\sigma^0)^\top \Lambda^0_t\big] 	\\
	&\ 
	 + \frac{1}{2} | \Gamma_t\, x + \Lambda^0_t\, \bar{\mu}(t,\omega,y) |^2  - \left ( \Gamma_t\, x + \Lambda^0_t\, \bar{\mu}(t,\omega,y) \right )\cdot \left ( b_t\, x + \bar{b}_t\, \bar{\mu}(t,\omega,y) \right ).
\end{aligned}
\]
We now collect terms (symmetrizing for the squared terms) to arrive at the following closed system of Riccati equations:
\begin{itemize}
	\item $|x|^2 \ : \ \dot{\Gamma}_t = \Gamma_t^\top \Gamma_t - \Gamma_t^\top b_t - b_t^\top \Gamma_t   , \ \ \ \Gamma_T = \bar{q} $ ,
	\item $x \bar{\mu}_t \ : \ \dot{\Lambda}_t^0 = (\Lambda^0_t)^\top \Lambda^0_t -(\Lambda^0_t)^\top\, (b_t + \bar{b}_t - \Gamma_t) + \Gamma_t^\top  \Lambda^0_t  - \Gamma_t^\top\, \bar{b}_t  - b_t^\top\, \Lambda_t^0 \, s,  \ \ \ \Lambda^0_T = - \bar{q}\, s  $,
	\item $\bar{\mu}_t^2 \ \  :  \ \dot{\Gamma}^0_t =  -(b_t + \bar{b}_t - \Gamma_t - \Lambda^0_t)^\top \Gamma^0_t -(\Gamma^0_t)^\top(b_t + \bar{b}_t - \Gamma_t - \Lambda^0_t)$ 
	\subitem $\quad \quad \quad \ \ +  (\Lambda^0_t)^\top \Lambda^0_t - (\Lambda^0_t)^\top \bar{b}_t - \bar{b}_t^\top \Lambda^0_t   ,\quad \quad \Gamma^0_T = s^\top \bar{q}\, s $,
	\item $1 \ \ \ \ : \   \dot{\Delta}_t =  - \frac{1}{2}\tr[(\sigma \sigma^\top + \sigma^0 (\sigma^0)^\top) \ \Gamma_t ] - \frac{1}{2} \tr[  \sigma^0 (\sigma^0)^\top \Gamma^0_t] 
	- \tr[\sigma^0 (\sigma^0)^\top \Lambda^0_t], \  \ \Delta_T = 0.
   $
\end{itemize} 
Notice that the equations for $\Gamma_t, \Lambda^0_t$ are quadratic
Riccati equations, while the equation for $\Gamma^0_t$
is linear. 

For the forward equation, we can write
$$
\hat{m}(t,x,\omega,y) = g\big(t,x- \bar{\mu}(t,\omega,y)\big).
$$
We may now compute, letting $Z = x-\bar{\mu}(t,\omega,y)$,
\begin{align*}
D_{t,\omega}^y \hat{m}(t,x,\omega,y) =&\ \partial_t g(t,Z) - \partial_z g(t,Z)\cdot D_{t,\omega}^y\bar{\mu}(t,\omega)\\
=&\  \partial_t g(t,Z) - \partial_z g(t,Z)\cdot  (b_t + \bar{b}_t - \Gamma_t - \Lambda^0_t)\, \bar{\mu}(t,\omega).
\end{align*}

Also,
\begin{align*}
\partial_y \hat{m}(t,x,\omega,y) =&\  - (\partial_y\, \bar{\mu}(t,\omega,y))^\top \partial_z g(t,Z)\\
=&\  -(\sigma^0)^{\top} \partial_z g(t,Z)\\
\partial^2_{yy} \hat{m}(t,x,\omega,y) =&\ {\rm tr}\big[\sigma^0 (\sigma^0)^\top \partial^2_{zz} g(t,Z)\big].
\end{align*}
We see that this makes the second order terms become
$$
\frac{1}{2} \tr[(\sigma\sigma^\top + \sigma^0(\sigma^0)^{\top} \ \partial_{xx}^2 \hat{m}(t,x,\omega,y) + \partial^2_{yy} \hat{m}(t,x,\omega,y) 
	+ 2 \sigma^0 \partial_x \partial_y \hat{m}(t,x,\omega,y) \big]  = \frac{1}{2}{\rm tr}\big[\sigma \sigma^\top \partial^2_{zz} g(t,Z)\big].
$$

We also have that
\begin{align*}
&\ \text{div}_x \left [ \hat{m}(t,x,\omega,y) \partial_p H\big(t,x,\partial_x \hat{u}(t,x,\omega, y),\hat{\mu}(t,\omega,y)\big)   \right ]\\
 =&\ \text{div}_x \left [ g(t,Z)\big] \big(\Gamma_t\, x + \Lambda^0_t\, \bar{\mu}(t,\omega,y)-b_t\, x-\bar{b}_t\, \bar{\mu}(t,\omega,y) \right ] \\
=&\ \text{div}_z \big[ g(t,Z)(\Gamma_t-b_t) Z\big] - \partial_z g(t,Z)\cdot (b_t + \bar{b}_t - \Gamma_t - \Lambda^0_t)\, \bar{\mu}(t,\omega).
\end{align*}

Plugging them together, we see that $g$ solves
\begin{align*}
\partial_tg(t,z) =&\  \frac{1}{2}{\rm tr}\big[\sigma \sigma^\top \partial^2_{zz} g(t,z)\big] + \text{div}_z \big[ g(t,z)(\Gamma_t-b_t) z\big].
\end{align*}
\end{proof}

\subsection*{Discussion of the Solvability of the Riccati Equations}

Standard ODE theory applies to guarantee there exists a unique solution to the system of equations for at least a short time.  The only barrier to global existence is if the matrices $\Gamma_t$ or $\Lambda^0_t$ diverge (since the $\mu^2$ equation is linear in $\Gamma_t^0$ it does not pose a barrier to global existence). An upper bound, in the sense of positive semidefinite matrices, for $\Gamma_t$ will always hold by a Gronwall argument: that $\Gamma_t\leq M_t$ where $M_t$ solves the linear ODE
$$
	\dot{M}_t = -M_t\, b_t - b_t^\top M_t , \ \ M_T=\bar{q}.
$$

A lower bound of $\Gamma_t\geq 0$ holds so long as $\bar{q}$ is positive semidefinite.

For $\Lambda_t^0$, we consider $\tilde{\Lambda}_t= \Gamma_t+\Lambda_t^0$, which solves:
$$
	\dot{\tilde{\Lambda}}_t = \tilde{\Lambda}_t^\top\, \tilde{\Lambda}_t  -b_t^\top\, \tilde{\Lambda}_t - \tilde{\Lambda}_t^\top\, (b_t+\bar{b}_t),
$$
with $\tilde{\Lambda}_T = \bar{q}- \bar{q}\, s$.  We assume that $\bar{b}_t$ is a scalar times the identity matrix, so  that $\tilde{\Lambda}_t$ remains symmetric.

Similar to the argument for $\Gamma_t$, there is a global solution so long as $\bar{q}- \bar{q}\, s$ is positive semidefinite ($s\leq 1$).  This same result appears in \cite{bible1}, where an example is also given that shows how solutions exist only for a finite time period if the positive semidefinite condition fails for the problem data (that is, $\bar{q}-\bar{q}\, s\not\geq 0$).

\subsection*{Comparison with the literature}

The mean field game system with common noise
can be interpreted as the system of characteristics
for the master equation set on the Wasserstein space $\calP_2(\RR^d)$
of probability measures with finite second moment.
For the linear-quadratic data of \eqref{LQdataMFG},
the master equation has the form
(see display (4.41) of Carmona-Delarue \cite{bible2}): 
\begin{equation}
	\label{mastergeneral}
	\begin{cases}
		\begin{aligned}
			& - \partial_t U(t,x,\mu) - \frac{1}{2} \tr \left [(\sigma \sigma^\top + \sigma^0 (\sigma^0)^\top) \partial_{xx}^2 U(t,x,\mu) \right ] - (b_t x + \bar{b}_t \bar{\mu}) \cdot \partial_x U(t,x,\mu)  + \frac{1}{2}| \partial_x U(t,x,\mu) |^2 \\
			&	- \int_{\RR^d} \tr \left [ \sigma^0 (\sigma^0)^\top \partial_{x} \partial_\mu U(t,x,\mu)(v) \right ] \mu(dv)   	- \frac{1}{2}  \int_{\RR^d} \tr \left [(\sigma \sigma^\top + \sigma^0 (\sigma^0)^\top) \partial_{v} \partial_\mu U(t,x,\mu)(v) \right ] \mu(dv) 
			\\
			& - \frac{1}{2}  \int_{\RR^d} \int_{\RR^d}  \tr \left [\sigma^0 (\sigma^0)^\top  \partial_{\mu \mu}^2 U(t,x,\mu)(v,v') \right ] \mu(dv) \mu(dv') \\
			 & -  \int_{\R^d}\partial_\mu U(t,x,\mu)(v) \cdot  \Big ( b_t v + \bar{b}_t \bar{\mu}  -(\partial_x U)(t,v,\mu)   \Big )  \mu(dv)= 0, \ (t,x,\mu) \in [0,T) \times \RR^d \times \calP_2(\RR^d) \\
			 & U(T,x,\mu) = \frac{1}{2} (x - s\bar{\mu})\cdot \bar{q}\, (x - s\bar{\mu}), \ \ \ \  (x,\mu) \in \RR^d \times \calP_2(\RR^d).
		\end{aligned}
	\end{cases} 
\end{equation}
Here, ``$\partial_\mu$'' is the gradient on the Wasserstein space $\mc{P}_2(\R)$,
which can \emph{formally} be interpreted as ``$\partial_v \frac{\delta}{\delta \mu} U(t,x,\mu)(v)$,'' with $\frac{\delta}{\delta \mu}$ denoting the linear functional (i.e., Fr\'echet) derivative
in the vector space of all finite signed measures.

As mentioned above and as in display (22) of Cardaliaguet-Delarue-Lasry-Lions \cite{cdll1}, 
the relationship between a solution $(u_t(x), v_t(x), m_t(x))_\T$ of the characteristic equations \eqref{sMFG}
and a solution $U(t,x,\mu)$ of
\eqref{mastergeneral} should be given by $u_t(x) = U(t,x,\mu_t)$ with $\mu_t(dz) =m_t(z)dz$. 
Hence, we expect to have the same ansatz
\[
U(t,x,\mu)=\frac{1}{2} \left ( x \cdot \Gamma_t x + \bar{\mu} \cdot \Gamma^0_t \bar{\mu} +
2x \cdot \Lambda^0_t \bar{\mu} \right )
+ \Delta_t.
\]
We then compute
\[
\partial_t U(t,x,\mu)
= 
\frac{1}{2} \left ( x \cdot \dot{\Gamma}_t x + \bar{\mu} \cdot \dot{\Gamma}^0_t \bar{\mu} +
2x \cdot \dot{\Lambda}^0_t \bar{\mu} \right )
+ \dot{\Delta}_t
\]
\[
\partial_x U(t,x,\mu)
= 
\Gamma_t x + \Lambda^0_t \bar{\mu}, \ \ \
\partial_{xx}^2 U(t,x,\mu)
= 
\Gamma_t,
\]
\[
\partial_\mu U(t,x,\mu)(v) =
\Gamma_t^0 \bar{\mu} + (\Lambda^0_t)^\top x, \ \ \ 
\partial_{\mu \mu}^2 U(t,x,\mu)(v) =
\Gamma_t^0,
\]
\[
\partial_x \partial_\mu U(t,x,\mu)(v) =
 (\Lambda^0_t)^\top, \ \ \ 
 \partial_v \partial_\mu U(t,x,\mu)(v) = 0.
\]
Plugging these calculations
in the equation gives
\[
\begin{aligned}
	& -  \frac{1}{2} \left ( x \cdot \dot{\Gamma}_t x + \bar{\mu} \cdot \dot{\Gamma}^0_t \bar{\mu} +
2x \cdot \dot{\Lambda}^0_t \bar{\mu}  \right )
- \dot{\Delta}_t  \\ 
	&- \frac{1}{2} \tr \left [(\sigma \sigma^\top + \sigma^0 (\sigma^0)^\top) \Gamma_t \right ] - (b_t x + \bar{b}_t \bar{\mu}_t) \cdot \Big ( \Gamma_t x + \Lambda^0_t \bar{\mu} \Big )  + \frac{1}{2}| \Gamma_t x + \Lambda^0_t \bar{\mu} |^2 \\
			&	-  \tr \left [ \sigma^0 (\sigma^0)^\top \Lambda^0_t \right ]  
		- \frac{1}{2}   \tr \left [\sigma^0 (\sigma^0)^\top  \Gamma_t^0 \right ] -  \Big ( \Gamma_t^0 \bar{\mu} +  (\Lambda^0_t)^\top x \Big ) \cdot  \Big ( b_t  + \bar{b}_t - \Gamma_t - \Lambda^0_t  \Big )\bar{\mu} \\
		& \quad \quad  \quad  \quad  \quad  \quad  \quad   = 0.
\end{aligned}
\]
We then arrive at the same set of equations as in Section \ref{MFGtypical}.

\section{Mean field game with controlled common noise}
\label{controlledvolatilitysection}

Suppose we have a more general state process $(X^{\balpha}_t)_\T$ with dynamics of the form
\[
dX_t = b(t,x,\mu_t, \alpha_t) dt
+ \sigma(t,x,\mu_t,\alpha_t) dW_t
+ \sigma^0(t,x,\mu_t,\alpha_t) dW^0_t.
\]
Write  $a(t,x,\mu,\alpha):= \left ( \sigma \sigma^\top  +\sigma^0(\sigma^0)^\top \right ) (t,x,\mu,\alpha)$
and define 
\[
\alpha^*(t,x,\mu, p,X,Q):=
\argmin_\alpha \bigg \{ \frac{1}{2} \tr \Big [a(t,x,\mu,\alpha) X   \Big ]  + \tr \Big [\sigma^0(t,x,\mu,\alpha) Q  \Big ] + p \cdot b(t,x,\mu,\alpha)  + f(t,x,\mu,\alpha)   \bigg \}.
\]
Then, given an $\FF^{\bW^0}$-adapted measure flow $\bmu = (\mu_t)_\T$, we need to 
find a pair $(u_t(x),v_t(x))_\T$ of $\FF^{\bW^0}$-adapted random fields satisfying the stochastic HJB
\begin{equation}
	\label{generalsHJB}
	\begin{aligned}
	d u_t(x) & = -  \frac{1}{2} \tr \Big [a\Big(t,x,\mu_t,\alpha^*(t,x,\mu_t, \partial_x u_t(x), \partial_{xx}^2 u_t(x) ,\partial_x v_t(x)) \Big )  \ \partial_{xx}^2 u_t(x) \Big ]        dt \\
	& \quad -  \tr \Big [\sigma^0\Big(t,x,\mu_t,\alpha^*(t,x,\mu_t, \partial_x u_t(x), \partial_{xx}^2 u_t(x) ,\partial_x v_t(x)) \Big)  \ \partial_{x} v_t(x) \Big ]      dt \\
	& \quad - \partial_x u_t(x) \cdot b\Big( t,x,\mu_t,\alpha^*(t,x,\mu_t, \partial_x u_t(x), \partial_{xx}^2 u_t(x) ,\partial_x v_t(x)) \Big) dt \\
	& \quad  - f\Big( t,x,\mu_t,\alpha^*(t,x,\mu_t, \partial_x u_t(x), \partial_{xx}^2 u_t(x) ,\partial_x v_t(x)) \Big) dt + v_t(x) \cdot dW^0_t.
	\end{aligned} 
\end{equation}
Besides being fully nonlinear, this stochastic HJB poses a new difficulty
of the optimizer $\alpha^*$ potentially introducing additional nonlinearities based
on the unknown random field $v_t(x)$. 

Fortunately, the lifted functional approach shows how to reduce consideration to a more 
classical-looking scenario.
Indeed, the (fully nonlinear) compensated HJB equation involves finding a lifted functional $\hat{u}(t,x,\omega,y)$ satisfying
\begin{align}
	\label{generalcompensatedHJB}
	 \mathcal{D}_{\omega}^y \hat{u}(t,x,\omega,y)=&\ -\partial_t \hat{u}(t,x,\omega,y)  -  \frac{1}{2} \tr \Big [A\Big(t,x,\hat{m},\alpha^*(t,x,\hat{m}, \partial_x \hat{u}, \partial_{xx}^2 \hat{u} ,\partial_{xy} \hat{u}) \Big )  \partial^2_{(x,y)} \hat{u} \Big ] 
	 \\
	 &     - \partial_x \hat{u} \cdot b\Big( t,x,\hat{m},\alpha^*(t,x,\hat{m} , \partial_x \hat{u}, \partial_{xx}^2 \hat{u} ,\partial_{xy} \hat{u}) \Big)  - f\Big( t,x,\hat{m},\alpha^*(t,x,\hat{m}, \partial_x \hat{u}, \partial_{xx}^2 \hat{u} ,\partial_{xy} \hat{u}) \Big),\nonumber
\end{align}
where $\partial^2_{(x,y)}$ is the Hessian in $(x,y)$ and where
\[
A(t,x,\mu,\alpha) := 
\begin{pmatrix}
\big( \sigma \sigma^\top +\sigma^0 (\sigma^0)^\top \big)(t,x,\mu,\alpha) & \sigma^0(t,x,\mu,\alpha) \\
\sigma^0(t,x,\mu,\alpha) & I_d
\end{pmatrix}.
\]

\subsection*{Linear Quadratic data for controlled volatility}

For simplicity, we work in dimension $d=1$, though the manipulations below may be generalized to higher dimensions.
Set the linear-quadratic cost data similarly as before to
\[	
f(t,x,\mu,\alpha) := \frac{1}{2} (\alpha -\bar{a}_t)^2,
\]
\[ 
g(x,\mu) := \frac{\bar{q}}{2} (x - s\bar{\mu})^2,
\]
(so the only difference is that we add the given parameter $\bar{a}_t$).
For the dynamics, we take
\[
b(t,x,\mu,\alpha) : = b_t x+\bar{b}_t \bar{\mu}, \ \ \ 
\sigma(t,x,\mu,\alpha):= \sigma, \ \  \ 
\sigma^0(t,x,\mu,\alpha):= \alpha.
\]
The optimality condition then becomes
\[
\begin{aligned}
\alpha^*(t,x,\mu, p,X,Q):=
\argmin_\alpha \bigg \{ \frac{1}{2} \alpha^2 X  + \alpha Q  + \frac{1}{2} (\alpha-\bar{a}_t)^2  \bigg \}
= \frac{\bar{a}_t-Q}{1+X},
\end{aligned}
\]
in the case that $X>-1$ and a minimizer exists.  The compensated HJB becomes 
\begin{equation}
	\label{LQgeneralcompensatedHJB}
	\begin{aligned}
	 -& D_{t,\omega}^y  \hat{u} -  
	 \frac{1}{2} \left (\sigma^2 + \left ( \frac{\bar{a}_t-\partial_{xy} \hat{u}}{1+\partial_{xx}^2 \hat{u}}    \right )^2   \right ) \partial_{xx}^2 \hat{u} 
	 - \frac{\bar{a}_t-\partial_{xy} \hat{u}}{1+\partial_{xx}^2 \hat{u}}   \partial_{xy}\hat{u} - \frac{1}{2} \partial^2_{yy} \hat{u}
	 \\
	 &     - \partial_x \hat{u} \left ( b_t x+\bar{b}_t \bar{\mu} \right ) = \frac{1}{2} \Big |\frac{\bar{a}_t-\partial_{xy} \hat{u}}{1+\partial_{xx}^2 \hat{u}}    -\bar{a}_t \Big |^2.
	\end{aligned} 
\end{equation}

\begin{thm}\label{thm:prob3}
Suppose that 
$$
\dot{\Gamma}_t =  -2\, b_t\, \Gamma_t  , \ \ \ \Gamma_T = \bar{q} 
$$
has a solution on $[0,T]$ with $\Gamma_t>-1$. We let $\Gamma^0_t,\Lambda^0_t,\Delta_t$ denote the solutions to
\begin{align*}
	\dot{\Lambda}_t^0 =&\ -\Lambda^0_t b_t - \Gamma_t \, \bar{b}_t - \Lambda^0_t(b_t + \bar{b}_t) ,  & \Lambda^0_T = - s\, \bar{q}  ,\\
	 \dot{\Gamma}^0_t =&\  -2\, \Gamma_t^0(b_t + \bar{b}_t) - 2\, \Lambda^0_t \bar{b}_t  , & \Gamma^0_T = s\, \bar{q} s,\\
	\dot{\Delta}_t =&\  - \frac{1}{2} \left ( \frac{\bar{a}_t}{1+\Gamma_t + \Lambda^0_t} \right )^2 \, \left ( \Gamma_t + 2\Lambda^0_t + \Gamma_t^0 \right )  -  \frac{1}{2} ( \sigma^2 + \bar{a}_t^2)  , & \Delta_T = 0.
\end{align*}
A solution to (\ref{LQgeneralcompensatedHJB})  is given by
\begin{align*}
\hat{u}(t,x,\omega,y) =&\  
\frac{1}{2} \left ( x \cdot \Gamma_t\, x + \bar{\mu}(t,\omega, y) \cdot \Gamma^0_t\, \bar{\mu}(t,\omega, y) +
2x \cdot \Lambda^0_t\, \bar{\mu}(t,\omega, y) \right )
+ \Delta_t,
\end{align*}
where
\[
\bar{\mu}(t,\omega, y)
=\Phi_t \bigg (\bar{\mu}_0 
+ \int_0^t \Phi^{-1}_s \left [  (b_s + \bar{b}_s) \alpha_s - \dot{\alpha}_t  \right ] \omega_s ds \bigg )
+ \alpha_t y 
\]
and the optimal control is
$$
\alpha_t = (1+\Gamma_t + \Lambda^0_t)^{-1} \bar{a}_t.
$$
The forward equation is solved by
$$
\hat{m}(t,x,\omega,y) = g\big(t,x-\bar{\mu}(t,\omega,y)\big),
$$
where $g$ solves $g(0,\cdot) = \ell$ and
\begin{align}\label{eqn:forward_controlled_common_noise}
\partial_t g(t,z) = \frac{\sigma^2}{2}\partial^2_{zz} g(t,z) - \partial_z \big(g(t,z)\, b_t\, z\big).
\end{align}
\end{thm}

\begin{proof}
We adopt a similar ansatz as before, namely,
\[
\hat{u}(t,x,\omega,y) = 
\frac{1}{2} \left ( x\, \Gamma_t\, x + \bar{\mu}(t,\omega, y) \, \Gamma^0_t\, \bar{\mu}(t,\omega, y) +
2x \, \Lambda^0_t\, \bar{\mu}(t,\omega, y) \right )
+ \Delta_t,
\]
so that the optimal feedback function has the lifted form
\[
\hat{\alpha}^*(t,\omega,y)=
\frac{\bar{a}_t-\Lambda^0\, \partial_y \bar{\mu}(t,\omega,y)}{1+\Gamma_t}  
\]
But this expression is a bit problematic because the term ``$\partial_y \bar{\mu}(t,\omega,y)$''
will likely involve the control itself.  

To resolve this issue, let us search for the optimal control among 
deterministic $C^1$ functions of time $\bbeta = (\beta_t)_\T$.
Indeed, given such a function, the associated state dynamics will have the form
\[
dX^\beta_t = ( b_t X^\beta_t +\bar{b}_t \bar{\mu}^\beta_t) dt
+ \sigma dW_t
+ \beta_t dW^0_t.
\]
As before, we can take expectations of 
this equation conditional on 
 $\calF^{\bW^0}_t$ to get
\[
d\bar{\mu}^\beta_t = (b_t + \bar{b}_t) \bar{\mu}^\beta_t \ dt + \beta_t dW^0_t.
\]
And again, as before, the lifted functional of $\bar{\mu}^\beta_t = \bar{\mu}^\beta(t,\bW^0,W^0_t)$ can be solved explicitly as
\[
\bar{\mu}^\beta(t,\omega, y)
= \Phi_t \bigg (\bar{\mu}_0 
+ \int_0^t \Phi^{-1}_s \left [  (b_s + \bar{b}_s) \beta_s - \dot{\beta}_t  \right ] \omega_s ds \bigg )
+ \beta_t y, 
\]
where $(\Phi_t)_\T$ is the solution of 
\[
\dot{\Phi}_t = (b_t + \bar{b}_t)\Phi_t, \ \ \Phi_0 = 1.
\]
From this last expression, we can then compute directly
\[
 D_{t,\omega}^y  \bar{\mu}^\beta(t,\omega,y)= (b_t + \bar{b}_t) \bar{\mu}^\beta(t,\omega,y), \ \ \ 
\partial_y \bar{\mu}^\beta(t,\omega,y) = \beta_t.
\]
Hence, given a flow of measures $\bmu = (\mu_t)_\T$ determined by a deterministic $C^1$ control $\bbeta = (\beta_t)_\T$,
the optimal control will satisfy (now removing the dependence on $\omega,y$)
\[
\hat{\alpha}^*_t =
\frac{\bar{a}_t-\Lambda^0_t \beta_t }{1+\Gamma_t}.  
\]
But the mean field game consistency condition implies we have $\beta_t = \hat{\alpha}^*_t$, resulting in a readily solved equation for $\hat{\alpha}^*_t$, namely,

\[
\hat{\alpha}^*_t =
\frac{\bar{a}_t-\Lambda^0_t \hat{\alpha}^*_t }{1+\Gamma_t}, \ \ \text{so that} \ \ 
\hat{\alpha}^*_t = (1+\Gamma_t + \Lambda^0_t)^{-1} \bar{a}_t.
\]
In particular, the optimal control $\hat{\alpha}^*_t$ is a deterministic function of time
and the lifted function ``$\bar{\mu}(t,\omega,y)$'' appearing in the ansatz for $\hat{u}(t,x,\omega,y)$ may be taken to satisfy:
\[
 D_{t,\omega}^y  \bar{\mu}(t,\omega,y)= (b_t + \bar{b}_t) \bar{\mu}(t,\omega,y),  \ \ \ 
\partial_y \bar{\mu}(t,\omega,y) = \hat{\alpha}^*_t = (1+\Gamma_t + \Lambda^0_t)^{-1} \bar{a}_t.
\]

At last, we can plug all these considerations into the compensated HJB to get
\[
	\begin{aligned}
	 -& \frac{1}{2} \left ( x\, \dot{\Gamma}_t x + \bar{\mu}(t,\omega, y) \, \dot{\Gamma}^0_t \bar{\mu}(t,\omega, y) +
2x\, \dot{\Lambda}^0_t \bar{\mu}(t,\omega, y) \right )
- \dot{\Delta}_t \\
& - \Gamma^0_t (b_t + \bar{b}_t) \, \bar{\mu}^2(t,\omega, y)  -
 \Lambda^0_t (b_t + \bar{b}_t)\,  x\bar{\mu}(t,\omega, y)   \\ 
	 & -  \frac{1}{2} \left (\sigma^2 + \left ( \frac{\bar{a}_t}{1+\Gamma_t + \Lambda^0_t}    \right )^2   \right ) \Gamma_t 
	 - \left ( \frac{\bar{a}_t}{1+\Gamma_t + \Lambda^0_t} \right )^2  \Lambda^0_t - \frac{1}{2} \Gamma^0_t \left ( \frac{\bar{a}_t}{1+\Gamma_t + \Lambda^0_t} \right )^2
	 \\
	 &     - \Big ( \Gamma_t x + \Lambda^0_t \bar{\mu}(t,\omega, y) \Big ) \, \left ( b_t x+\bar{b}_t \bar{\mu} \right ) = \frac{1}{2}\Big |\frac{\bar{a}_t}{1+\Gamma_t + \Lambda^0_t}    -\bar{a}_t \Big |^2  ,
	\end{aligned} 
\]
with terminal condition
\[
\hat{u}(T,x,\omega,y) = \frac{\bar{q}}{2} (x - s\bar{\mu}(t,\omega,y))^2.
\]
This leads to the following system of ODEs (that can be solved in the order presented):
\begin{itemize}
	\item $|x|^2 \ : \ \dot{\Gamma}_t =  -2\, b_t\, \Gamma_t  , \ \ \ \Gamma_T = q+\bar{q} $, 
	\item $x \bar{\mu}_t \ : \ \dot{\Lambda}_t^0 = -\Lambda^0_t b_t - \Gamma_t \, \bar{b}_t - \Lambda^0_t(b_t + \bar{b}_t),  \ \ \ \Lambda^0_T = - s\, \bar{q}  $,
	\item $\bar{\mu}_t^2 \ \  :  \ \dot{\Gamma}^0_t =  -2\, \Gamma_t^0(b_t + \bar{b}_t) - 2\, \Lambda^0_t \bar{b}_t   ,\quad \quad \Gamma^0_T = s\, \bar{q} s $,
	\item $1 \ \ \ \ : \   \dot{\Delta}_t =  - \frac{1}{2} \left ( \frac{\bar{a}_t}{1+\Gamma_t + \Lambda^0_t} \right )^2 \, \left ( \Gamma_t + 2\Lambda^0_t + \Gamma_t^0 \right )  -  \frac{1}{2} ( \sigma^2 + \bar{a}_t^2)  , \  \ \Delta_T = 0
   $.
\end{itemize} 

For the forward equation, we can write
$$
\hat{m}(t,x,\omega,y) = g(t,Z),
$$
letting $Z = x-\bar{\mu}(t,\omega,y)$, and we compute
\begin{align*}
D_{t,\omega}^y \hat{m}(t,x,\omega,y) =&\ \partial_t g(t,Z) - \partial_z g(t,Z) D_{t,\omega}^y\bar{\mu}(t,\omega)\\
=&\  \partial_t g(t,Z) - \partial_z g(t,Z)\cdot  (b_t + \bar{b}_t)\, \bar{\mu}(t,\omega).
\end{align*}

Also, we have
\begin{align*}
\partial_y \hat{m}(t,x,\omega,y) =&\  - \big(\partial_y\, \bar{\mu}(t,\omega,y)\big) \partial_z g(t,Z)\\
=&\ - \alpha_t\, \partial_z g(t,Z),\\
\partial^2_{yy}\hat{m}(t,x,\omega,y) =&\  \alpha_t^2\, \partial^2_{zz} g(t,Z).
\end{align*}
We see that this makes the second order terms become
$$
\frac{1}{2} \left ( (\sigma^2+\alpha_t^2) \ \partial_{xx}^2 \hat{m}(t,x,\omega,y) + \partial^2_{yy} \hat{m}(t,x,\omega,y) 
	+ 2 \alpha_t\, \partial_x \partial_y \hat{m}(t,x,\omega,y) 
	  \right )  = \frac{1}{2}\sigma^2 \partial^2_{zz} g(t,Z).
$$

We also have that
\begin{align*}
&\ \partial_x \left [ \hat{m}(t,x,\omega,y) \partial_p H\big(t,x,\partial_x \hat{u}(t,x,\omega, y),\hat{\mu}(t,\omega,y)\big)   \right ]\\
 =&\ \partial_x \left [ g(t,Z)\big] \big(-b_t\, x-\bar{b}_t\, \bar{\mu}(t,\omega,y) \right ] \\
=&\ \partial_z \big[ g(t,Z)(-b_t) Z\big] - \partial_z g(t,Z)\cdot (b_t + \bar{b}_t)\, \bar{\mu}(t,\omega).
\end{align*}

Putting the equations together, we see that $g$ solves (\ref{eqn:forward_controlled_common_noise}).
\end{proof}

\subsection*{Discussion of Solvability of Riccati Equations}

As the system of ODEs for $\Gamma$, $\Lambda^0$, and $\Gamma^0$ is linear, there always exists a unique solution. We require $\Gamma_t>-1$ in order for $\alpha^*$ to correspond to the minimum in the Hamiltonian.  We then require $\Gamma_t+\Lambda_t^0\not=-1$ so that there exists a fixed point.  Both of these conditions hold in the case considered in Section \ref{MFGtypical}, where we assume that $\bar{q}\geq 0$ and $s\leq 1$, which implies that $\Gamma\geq 0$ and $\Lambda^0\geq 0$.

\subsection*{Comparison with the literature}

As mentioned in the introduction, we do not know many references on mean field games
with control in the volatility 
coefficient of the common noise except for the recent theoretical paper of Barrasso-Touzi 
\cite{barrasso2020controlled}
and sporadic statements throughout Carmona-Delarue \cite{bible1, bible2}.
However, we can still compare with an existing
explicitly solvable 
model of controlled
volatility in a more
classical stochastic control setting.
For example, 
Proposition 5.1 of Ankirchner-Fromm \cite{ankirchner2018optimal}
arrives at an optimal control that in our notation
would correspond to ``$ \frac{\bar{a}_t}{1+\Gamma_t}$''.
It is interesting
that we instead arrive at a slightly
modified form ``$ \frac{\bar{a}_t}{1+\Gamma_t + \Lambda^0_t}$,''
since the control is entangled with the additional unknown process ``$v_t(x)$'', as is clear from the stochastic HJB equation \eqref{generalsHJB}.

\section{Mean field game with common noise and partial information} \label{MFGpartialinfosection}

Recall we formulated the mean field game with common noise and partial information
as Problem 4 of Section \ref{problemformulations} with the observation process $\bZ = (Z_t)_\T$ driven by a drift $h(t, X_t, \mu_t)$ and an $\RR^{d'}$-valued Brownian motion $\tilde{\btheta}$ with covariance $\tilde{\Theta}$ that is independent of $\bW$ and $\bW^0$ under $\PP$. 
In this section, we start fresh: we introduce a new probability measure $\tilde{\PP}$ under which $\bW^0$ is still an $\RR^d$-valued Brownian motion, but now $\bZ = (Z_t)_\T$ is an $\RR^{d'}$-Brownian motion with covariance $\tilde{\Theta}$.

For given measures $\mu, \eta \in \calP(\RR^d)$ and an $\RR^d$-valued function $p\in C(\RR^d)^d$, 
we define the optimal feedback function under partial information as
\begin{equation}
	\label{controlfeedback}
\alpha^*(t,\mu,\eta,p) := \argmin_{\alpha \in \RR^d} \left \{  \int_{\RR^d} \eta(d\xi) \big [f(t,\xi,\mu,\alpha)+ p(\xi) \cdot b(t,\xi,\mu,\alpha)  \big]   \right \}.
\end{equation}
Next, for a given  $\FF^{\bW^0}$-adapted flow $\bmu = (\mu_t)_\T$ of probability measures in $\calP(\RR^d)$
and a $\FF^{\bW^0,\bZ}$-adapted flow $\bp = (p_t)_\T$ of functions in $C(\RR^d)^d$,
let $\boldsymbol{\eta}^{\bmu, \bp} = ({\eta}^{\bmu,\bp}_t)_\T$ denote the solution to the \emph{Kushner equation}
\begin{align*}
d \eta_t(dx) =&\ \bigg ( \frac{1}{2} \tr \left [(\sigma\sigma^\top + \sigma^0(\sigma^0)^\top)\partial_{xx}^2   \eta_t(dx) \right ] - \text{div}_x \left [\eta_t(dx) b(t,x, \alpha(t,\mu_t,\eta_t,p_t),\mu_t)  \right ]  \bigg ) \, dt \\
        &\ \quad - \text{div}_x \big [ \eta_t(dx) \sigma^0   \, dW^0_t  \big ] + \eta_t(dx)\left ( h(t,x,\mu_t) -  \bar{h}(t,\mu_t)  \right )^\top\tilde{\Theta}^{-1} \left ( dZ_t -   \bar{h}(t,\mu_t) \,  dt \right ),
\end{align*}
where
$$
\bar{h}(t,\mu_t) := \int_{\RR^d} \eta_t(d\xi)\,  h(t,\xi,\mu_t).
$$
The flow $\boldsymbol{\eta}^{\bmu, \bp}$ is the conditional law, with respect to the filtration $\FF^{\bW^0,\bZ}$ and measure $\tilde{\mathbb{P}}$, of the dynamics $\bX^{\bmu, \bp} = (X_t^{\bmu, \bp})_\T$ that solves
\[
dX_t = b(t,X_t,\mu_t,\alpha(t,\mu_t, \eta_t^{\bmu, \bp},p_t)) \, dt + \sigma \, dW_t + \sigma^0 \, dW^0_t.
\]
Lastly, given $\bmu$ and $\bX^{\bmu, \bp} = (X_t^{\bmu, \bp})_\T$ of sufficient regularity we define 
\[
M_t^{\bmu,\bp} := \exp \left \{ \int_0^t h(s,X_s^{\bmu,\bp}, \mu_s) \cdot \tilde{\Theta}^{-1}\, dZ_s -\frac{1}{2} \int_0^t h(s,X_s^{\bmu,\bp}, \mu_s)\cdot \tilde{\Theta}^{-1} h(s,X_s^{\bmu,\bp}, \mu_s) \, ds  \right \},
\]
i.e., the solution of
$$
dM_t = M_t\, h(t,X_t^{\bmu,\bp},\mu_t)\cdot \tilde{\Theta}^{-1} \, dZ_t, \quad M_0 = 1,
$$
which is a martingale under $\tilde{\PP}$ (where we recall that $\bZ $ is an $\RR^{d'}$-Brownian motion).
Then define a probability measure $\tilde{\PP}^{\bmu, \bp}$ by 
\begin{equation}
	\label{equivalentmeasures}
d \tilde{\PP}^{\bmu,\bp} = M_T^{\bmu,\bp} d\tilde{\PP}.
\end{equation}
By Girsanov's theorem, under the measure $\tilde{\PP}^{\bmu,\bp}$, the observation noise process $\tilde{\theta}^{\bmu,\bp}$ that solves
$$
d\tilde{\theta}_t = dZ_t-h(t, X_t^{\bmu,\bp}, \mu_t) \, dt,
$$
is a $d'$-dimensional Brownian motion with covariance $\tilde{\Theta}$, independent of $\bW$ and $\bW^0$. In other words, the measure $\tilde{\PP}^{\bmu,\bp}$ corresponds to the initial measure $\mathbb{P}$ under the change of variables from paths of $\tilde{\boldsymbol{\theta}}$ to the paths of $\mathbf{Z}$ using the solution of (\ref{eqn:partial_information_observations}).

We may now articulate the \emph{mean field game system with common noise and partial information}. \\
Given any probability measure $\lambda \in \calP_2(\RR^d)$,
find an $\FF^{\bW^0,\bZ}$-adapted quintuple 
\[
(u_t,v_t,k_t, \mu_t, \eta_t)\in C^2(\RR^d)\times C^1(\RR^d)^d \times C(\RR^d)^{d'}\times \mathcal{P}(\RR^d) \times \mathcal{P}(\RR^d),
\]
for $t\in [0,T]$, satisfying the following system, consisting of a stochastic HJB equation coupled with a forward Kushner equation: 
\begin{equation} \label{general}
	\begin{cases}
		d u_t(x) = \Big(- \frac{1}{2} \tr \left[ (\sigma\sigma^\top + \sigma^0(\sigma^0)^\top) \partial_{xx}^2 u_t(x) \right] - \partial_x u_t(x) \cdot b(t,x, \alpha^*(t,\mu_t,\eta_t,\partial_x u_t),\mu_t)\Big)dt \\
		 \quad  \quad  \quad \quad \quad + \left( v_t(x) \cdot dW^0_t - \tr\big[ \sigma^0 \partial_x v_t(x) \big] \, dt \right)   + k_t(x) \cdot  \left( dZ_t - h(t,x,\mu_t) \, dt \right)   - f(t,x,\alpha^*(t,\mu_t, \eta_t,\partial_x u_t),\mu_t)  dt\\
		d \eta_t(dx) = \bigg ( \frac{1}{2} \tr \left[(\sigma\sigma^\top + \sigma^0(\sigma^0)^\top)\partial_{xx}^2  \eta_t(dx) \right] - \text{div}_x \left[\eta_t(dx) b(t,x, \alpha^*(t,\mu_t, \eta_t,\partial_x u_t),\mu_t)  \right]  \bigg ) \, dt \\
		\quad  \quad  \quad \quad \quad - \text{div}_x \big [ \eta_t(dx) \sigma^0   dW^0_t  \big] + \eta_t(dx)\left( h(t,x,\mu_t) -  \bar{h}(t,\mu_t)  \right)^\top\tilde{\Theta}^{-1} \left( dZ_t -   \bar{h}(t,\mu_t) \, dt \right)  \\
		u_T(x) = g(x,\mu_T), \ \  \eta_0 = \lambda, \ \ \mu_t(dx) = \tilde{\EE}^{\bmu,\bp} [ \eta_t(dx)| \calF_t^{\bW^0}],	\ \  \bar{h}(t,\mu_t):= \int_{\RR^d} h(t,\xi,\mu_t) \eta_t(d\xi),
	\end{cases}
\end{equation}
where for the fixed point condition ``$\mu_t(dx) = \tilde{\EE}^{\bmu,\bp} [ \eta_t(dx)| \calF_t^{\bW^0}]$'', we let $\bp = (\partial_x u_t)_{\T}$. This equation can be justified by the so-called \emph{Kallianpur-Striebel formula}, which realizes $\eta_t$ as the conditional law of the state given $\mathcal{F}_t^{\bW^0,\bZ}$. Alternatively, 
conditioning on $\bW^0$ 
reduces the Kushner equation to the McKean–Vlasov equation for the law of the state process $\bX_t^{\bmu, \bp}$.
Hence, 
this fixed point condition can be seen as a consequence of the implicit consistency required of the forward-backward solution loop in (\ref{general}). To summarize, compared with the concrete control formulation of Problem 4 in Section \ref{problemformulations}, we are trading
the implicit condition required of the partial information constraint on the controls for the fixed point
condition required of the solution loop of the system (\ref{general}), just as in Bensoussan-Yam \cite{benyam}.

\subsection*{Lifted functional formulation}

The lifted functional form of the solution quintuple of the system \eqref{general} consists of $(\hat{u},\hat{n},\hat{\mu})$ such that
\begin{align}\label{eq:partialinfounlift}
&
u_t(x) = \hat{u}(t,x,\bW^0,\bZ,W^0_t,Z_t), \ \ \ v_t(x) = (\partial_y \hat{u})(t,x,\bW^0,\bZ,W^0_t,Z_t), \ \ \ k_t(x) = (\partial_z \hat{u})(t,x,\bW^0,\bZ,W^0_t,Z_t),\nonumber
\\
& \eta_t(dx) = \hat{n}(t,x,\bW^0,\bZ,W^0_t,Z_t)dx, \ \ \ \mu_t = \hat{\mu}(t, \bW^0, W^0_t).
\end{align}
where $(\hat{u}(t,x,\omega,\gamma,y,z), \hat{n}(t,x,\omega,\gamma,y,z), \hat{\mu}(t,\omega,y))$ solves the following system of equations:
\begin{align} \label{liftedsystemgeneral}
        \ -\partial_t\hat{u}&(t,x,\omega,\gamma,y,z) - \partial_x \hat{u}(t,x,\omega,\gamma,y,z) \cdot b\big(t,x, \hat{\alpha}(t,\omega,\gamma,y,z), \hat{\mu}(t,\omega,y)\big)\nonumber\\
		&\ -\frac{1}{2} \tr\Big[(\sigma\sigma^\top+\sigma^0(\sigma^0)^\top) \ \partial_{xx}^2 \hat{u}(t,x,\omega,\gamma,y,z)  + \partial^2_{yy} \hat{u}(t,x,\omega,\gamma,y,z)  +2\sigma^0 \partial_x \partial_y \hat{u}(t,x,\omega,\gamma,y,z)\Big]\nonumber\\
		&\ - \frac{1}{2} \tr\Big[\tilde{\Theta} \ \partial_{zz}^2 \hat{u}(t,x,\omega,\gamma,y,z)\Big]  - (\partial_z\hat{u})(t,x,\omega,\gamma,y,z) \cdot h\big(t,x,\hat{\mu}(t,\omega,y)\big)  \nonumber \\
		&\ - f\big(t,x, \hat{\alpha}(t,\omega,\gamma,y,z), \hat{\mu}(t,\omega,y)\big) =\mathcal{D}_{\omega,\gamma}^{y,z}\hat{u}(t,x,\omega,\gamma,y,z) , \nonumber
		\\
\ \partial_t \hat{n}&(t,x,\omega,\gamma,y,z) +\text{div}_x \big[\hat{n}(t,x,\omega,\gamma,y,z)\, b\big(t,x, \hat{\alpha}(t,\omega,\gamma,y,z), \hat{\mu}(t,\omega,y)\big)\big]\nonumber\\
&\ -\frac{1}{2} \tr\Big[(\sigma\sigma^\top+\sigma^0(\sigma^0)^\top) \ \partial_{xx}^2 \hat{n}(t,x,\omega,\gamma,y,z)  + \partial^2_{yy} \hat{n}(t,x,\omega,\gamma,y,z)  +2\sigma^0 \partial_x \partial_y \hat{n}(t,x,\omega,\gamma,y,z)\Big]\nonumber\\
		&\ - \frac{1}{2} \tr\Big[\tilde{\Theta} \ \partial_{zz}^2 \hat{n}(t,x,\omega,\gamma,y,z)\Big]  + (\partial_z\hat{n})(t,x,\omega,\gamma,y,z) \cdot h\big(t,x,\hat{\mu}(t,\omega,y)\big) = -\mathcal{D}_{\omega,\gamma}^{y,z} \hat{n}(t,x,\omega,\gamma,y,z),\nonumber\\
		\ \hat{u}&(T,x,\omega,\gamma,y,z) = \hat{g}\big(x,\hat{\mu}(T,\omega,y)\big),\ \hat{n}(0,x,\omega,\gamma,y,z) = \ell(x-\sigma^0y),\nonumber\\
		\ \hat{\mu}&(t,\omega,y)(dx) = \mathbb{E}[\hat{n}(t,x,\bW^0,\bZ,y,Z_t)dx | \mathcal{F}_t^{\bW^0}](\omega).
\end{align}
Here, the $\hat{\alpha}$ is given by
\begin{align}\label{eq:partialinfocontrol}
\hat{\alpha}&(t,\omega,\gamma,y,z):=  \\
&\ \argmin_{\alpha \in \RR^d} \left \{  \int_{\RR^d} \hat{n}(t,\xi,\omega,\gamma,y,z) \left [f(t,\xi,\hat{\mu}(t,\omega,y),\alpha)+ \partial_x \hat{u}(t,\xi,\omega,\gamma,y,z) \cdot b(t,\xi,\hat{\mu}(t,\omega,y),\alpha)  \right ] d\xi  \right \}.\nonumber
\end{align}
The equation for $\hat{\mu}$ in \eqref{liftedsystemgeneral} indicates that, in contrast to \eqref{liftedsystemintro},
the system \eqref{liftedsystemgeneral} is not quite path-by-path
in the sense that it requires an average over $\gamma$ to determine $\hat{\mu}$. To take this average, we must first solve
\begin{align}\label{eqn:info_dynamic}
dX_t =&\  b(t,X_t,\mu_t,\alpha^*_t) \, dt + \sigma \, dW_t + \sigma^0 \, dW^0_t,\nonumber\\
dZ_t =&\ h(t,X_t,\mu_t)dt + d\tilde{\theta}_t,
\end{align}
with $\mu_t = \hat{\mu}(t,\bW^0,W^0_t)$ and $\alpha^*_t = \hat{\alpha}(t,\bW^0,\bZ,W^0_t,Z_t)$.
Note the compensated time derivative $ D_{t,\omega,\gamma}^{y,z}$ will now involve both path variables $\omega$, $\gamma$, and present locations $y,z$:
\begin{align*}
	 D_{t,\omega,\gamma}^{y,z}  \hat{\varphi}(t,\omega,\gamma,y,z)  := \lim_{\epsilon \downarrow 0}  \epsilon^{-1} \left [\hat{\varphi}\Big(t+\epsilon,\omega\, 1_{[0,t)} + y\, 1_{[t,T]},\gamma\, 1_{[0,t)} + z\, 1_{[t,T]}, y,z\Big) - \hat{\varphi}(t,\omega,\gamma, y,z) \right ].
\end{align*}

The equations (\ref{liftedsystemgeneral}) should be understood with the flow transformation condition and an additional observation condition \eqref{eqn:consistency} to prescribe the relations between the state and common noise as well as the observation process.

We now restate one of our main results, Theorem \ref{thm:strong_solutions} in the context of partial information, which is that solutions of the lifted system provide strong solutions to the mean-field game system. We consider now the lifted functionals in an extend space $\hat{\mathbb{C}}^{1,2,2,2}$ that is twice differentiable in the spatial variables and lifted functionals of paths of common noise as well as the observation process.
\begin{thm}
Suppose that $\hat{u},\hat{n},\hat{\mu}\in \hat{\mathbb{C}}^{1,2,2,2}$ solve (\ref{liftedsystemgeneral}) and let the quintuple $(u,v,k,\mu,\eta)$ be defined as in (\ref{eq:partialinfounlift}). Suppose that $(X,Z)$ solve the stochastic differential equations of (\ref{eqn:info_dynamic}). 

Furthermore, suppose the flow transformation and observation conditions hold that
\begin{align}\label{eqn:consistency}
&\ (\sigma^0)^\top\, \partial_x \hat{n}(t,x,\bW^0,\bZ,W^0_t,Z_t) + \partial_y \hat{n}(t,x,\bW^0,\bZ,W^0_t,Z_t) = 0,\nonumber\\
\ -\hat{n}&(t,x,\bW^0,\bZ,W^0_t,Z_t)\tilde{\Theta}^{-1}\Big(h(t,x,\hat{\mu}(t,\bW^0,W^0_t)) - \int_{\mathbb{R}^d}\hat{n}(t,\xi,\bW^0,\bZ,W^0_t,Z_t)h(t,\xi,\hat{\mu}(t,\bW^0,W^0_t))d\xi\Big)\nonumber\\
&\ +\partial_z \hat{n}(t,x,\bW^0,\bZ,W^0_t,Z_t)= 0.
\end{align}
Then the quintuple $(u,v,k,\mu,\eta)$ solves (\ref{general}).
\end{thm}
\begin{proof}
We first note that from the definitions we have
$$
\alpha^*(t,\mu_t, \eta_t, \partial_x u_t) = \hat{\alpha}(t,\bW^0,\bZ,W^0_t,Z_t).
$$
We can then calculate using the lifted It\^{o} formula of Lemma \ref{lem:ito}
\begin{align*}
d u_t(x)=&\    \partial_t \hat{u}(t,x,\bW^0,\bZ, W^0_t, Z_t)dt + \partial_y\hat{u}(t,x,\bW^0,\bZ, W^0_t, Z_t)\cdot dW^0_t + \partial_z\hat{u}(t,x,\bW^0,\bZ, W^0_t, Z_t)\cdot dZ_t \\
&\ +\frac{1}{2} {\rm tr}\Big[\tilde{\Theta}\partial^2_{zz} \hat{u}(t,x,\bW^0,\bZ, W^0_t, Z_t) + \partial^2_{yy}\hat{u}(t,x,\bW^0,\bZ, W^0_t, Z_t) \Big]dt \\
=&\  - \frac{1}{2} \tr \Big[ (\sigma\sigma^\top+\sigma^0(\sigma^0)^\top) \  \partial_{xx}^2 u_t(x) + 2\sigma^0 \partial_x v_t(x)\Big]dt - \partial_x u_t(x) \cdot b(t,x, \alpha^*(t,\mu_t,\eta_t,\partial_x u_t),\mu_t)dt  \\
		 &\ + v_t(x) \cdot  dW^0_t    + k_t(x) \cdot (dZ_t-h(t,x,\mu_t)dt\big)- f(t,x,\alpha^*(t,\mu_t, \eta_t,\partial_x u_t),\mu_t)  \bigg )  dt.
\end{align*}
For the forward equation, we also make use of the consistency equations
\begin{align*}
d \eta_t(dx)=&\    \partial_t \hat{n}(t,x,\bW^0,\bZ, W^0_t, Z_t)dx\, dt \\
&\ + \partial_y \hat{n}(t,x,\bW^0,\bZ, W^0_t,Z_t)dx\cdot dW^0_t +  \partial_z \hat{n}(t,x,\bW^0,\bZ, W^0_t, Z_t)dx\cdot dZ_t\\
&\ + \frac{1}{2} \tr\Big[\tilde{\Theta}\partial^2_{zz} \hat{n}(t,x,\bW^0,\bZ, W^0_t, Z_t) +\partial^2_{yy}\hat{n}(t,x,\bW^0,\bZ, W^0_t, Z_t) \Big] dx\, dt\\
=&\ \bigg ( \frac{1}{2} \tr \left [(\sigma\sigma^\top+\sigma^0(\sigma^0)^\top)\partial_{xx}^2 \eta_t(dx) \right ] - \text{div}_x \left [\eta_t(dx) b(t,x, \alpha^*(t,\mu_t, \eta_t,\partial_x u_t),\mu_t)  \right ]  \bigg ) \, dt \\
		&\ - \text{div}_x \big [ \eta_t(dx) \sigma^0   dW^0_t  \big ] + \eta_t(dx)\left ( h(t,x,\mu_t) -  \bar{h}(t,\mu_t)  \right )^\top\tilde{\Theta}^{-1} \left ( dZ_t -   \bar{h}(t,\mu_t) \, dt \right ).
\end{align*}
Recall, the compensators are zero along the common noise and observation process, so plugging in equation \eqref{liftedsystemgeneral} implies the quintuple $(u,v,k,\mu,\eta)$ solves \eqref{general}.

The initial and terminal conditions of \eqref{general} follow immediately. The consistency equation of \eqref{general} for $\mu_t$ follows from the definition of (\ref{equivalentmeasures}) to relate the measure $\tilde{\mathbb{P}}^{\bmu,\mathbf{p}}$ to the original measure $\mathbb{P}$ where $\tilde{\boldsymbol\theta}$ is a Brownian motion.
\end{proof}

Another approach based on the Zakai-Stratonovich equations and flow transformation method of Lions-Souganidis is given in Appendix \ref{apx:Zakai}.

\subsection*{Linear-quadratic MFG with common noise and partial information}

We consider now the linear-quadratic assumptions that
\begin{align*}
b(t,x,\mu,\alpha) =&\  b\, x+\bar{b}\, \bar{\mu} + \alpha\\
f(t,x,\mu,\alpha) =&\ \frac{1}{2}|\alpha|^2\\
h(t,x,\mu) =&\ H\, x\\
g(t,x,\mu) =&\ \frac{1}{2}(x-s\, \bar{\mu})\cdot \bar{q}(x-s\, \bar{\mu}).
\end{align*}

\begin{thm}
Suppose that data is linear-quadratic and $\ell_0$ is Gaussian with mean $\bar{\mu}_0$ and covariance $\Pi_0$, and the following system of Riccati equations has a solution on $[0,T]$:
\begin{align*}
	\dot{\Sigma}_t =&\  - \Sigma_t^\top b_t - b_t^\top \Sigma_t   + \Lambda^1_t\, \Pi_t\, H^\top\, \tilde{\Theta}^{-1}\, H+ H^\top\, \tilde{\Theta}^{-1}\, H\, \Pi_t\, {\Lambda_t^1}^\top  , & \Sigma_T = \bar{q} \\
	\dot{\Lambda}_t^0 =&\ -{\Lambda^0_t}^\top\, (b_t + \bar{b}_t -\Sigma_t-\Lambda_t^1-\Lambda^0_t ) -{\Lambda^1_t}\, (\bar{b}_t - \Lambda^0_t)-\Sigma_t^\top\, (\bar{b}_t - \Lambda^0_t)- b_t^\top \Lambda^0_t   , & \Lambda^0_T = -  \bar{q}\, s  \\
	\dot{\Sigma}^0_t =&\  -( b_t + \bar{b}_t-\Sigma_t -\Lambda_t^1-\Lambda^0_t )^\top \Sigma^0_t-{\Sigma_t^0}^\top\, ( b_t + \bar{b}_t-\Sigma_t-\Lambda_t^1 -\Lambda^0_t ) - {\Lambda^0_t}^\top\, ( \bar{b}_t-\Lambda^0_t) & \\
	&\quad \quad \quad \quad  -( \bar{b}_t-\Lambda^0_t)^\top \Lambda^0_t-{\Lambda_t^0}^\top\, \Lambda_t^0, & \Sigma^0_t = s^\top \bar{q} s \\
	\dot{\Delta}_t =&  - \frac{1}{2}\tr[(\sigma \sigma^\top + \sigma^0 (\sigma^0)^\top) \ \Sigma_t ] - \frac{1}{2} \tr[  \sigma^0 (\sigma^0)^\top (\Sigma^0_t + \Sigma^1_t)  ] 
	- \tr[\sigma^0 (\sigma^0)^\top (\Lambda^0_t + \Lambda^1_t)] &\\
	& \quad \quad \quad \quad  - \frac{1}{2}  \tr[ \ \Pi_t H^\top H \Pi_t^\top \tilde{\Theta}^{-1} \Sigma^1_t  ] , &   \Delta_T = 0   \\
	\dot{\Sigma}^1_t =&\ - \left ( L_t^\top \Sigma^1_t +  (\Sigma^1_t)^\top L_t \right )  + \left ( (\Sigma_t + \Lambda^1_t)^\top \Lambda^1_t + ( \Lambda^1_t)^\top (\Sigma_t + \Lambda^1_t) \right ) - (\Sigma_t + \Lambda^1_t)^\top (\Sigma_t + \Lambda^1_t), & \Sigma^1_T = 0 \\
	\dot{\Lambda}^1_t =&\ -  \Lambda^1_t\, L_t -  H^\top\, \tilde{\Theta}^{-1}\, H\, \Pi_t^\top\, {\Sigma^1_t}^\top + \Sigma_t^\top (\Sigma_t + \Lambda^1_t) - b_t^\top \Lambda^1_t, &  \Lambda^1_T = 0,\\
	0=&\ -\Sigma^1_t (\bar{b}_t - \Lambda^0_t)  - \Lambda^1_t (\bar{b}_t - \Lambda^0_t) + \Sigma_t \Lambda^0_t - \Sigma_t \Lambda^0_t + (\Sigma^1_t + \Lambda^1_t) (\bar{b}_t - \Lambda^0_t).
\end{align*}

 Then a solution to (\ref{liftedsystemgeneral}) and (\ref{eqn:consistency}) is given by
\begin{align*}
\hat{u}(t,x,\omega,\gamma,y,z)=&\ \frac{1}{2}\, {x} \cdot \Sigma_t\, {x} +  {x} \cdot \Lambda^0_t\, \bar{\mu}(t,\omega,y)+\frac{1}{2}\,  \bar{\mu}(t,\omega,y) \cdot \Sigma^0_t\, \bar{\mu}(t,\omega,y)\\
& \quad  +  \frac{1}{2}\, \bar{\eta}(t,\omega,\gamma,y,z) \cdot \Sigma^1_t\, \bar{\eta}(t,\omega,\gamma,y,z) + {x}\cdot  \Lambda^1_t\, \bar{\eta}(t,\omega,\gamma,y,z)  + \Delta_t,
\end{align*}
with 
\[
\begin{aligned}
\bar{\eta}(t,\omega, \gamma, y,z)
& = \Psi_t\,  \bigg( \bar{\mu}_0 + \int_0^t  \Psi_s^{-1} (\bar{b}_s - \Lambda^0_s) \bar{\mu}(s,\omega, \omega_s) ds - \int_0^t L_s\, \Psi_s^{-1}\, \sigma^0\, \omega_s \, ds  \bigg )
 \\
& \quad + \Psi_t \, \bigg( \int_0^t \Big ( L_s \Psi_s^{-1}\Pi_s + \Psi_s^{-1} \dot{\Pi}_s \big )  H^\top \tilde{\Theta}^{-1} \gamma_s ds    \bigg) + \sigma^0 y + \Pi_t H^\top \tilde{\Theta}^{-1} z
\end{aligned}
\]
and
\[
\bar{\mu}(t,\omega, y)
= \Phi_t  \bigg (\bar{\mu}_0 
+ \int_0^t \Phi_s^{-1}\, \hat{L}_s\, \sigma^0\, \omega_s\, ds \bigg )
+ \sigma^0 y
\]
and 
\begin{equation*}
\dot{\Pi}_t = \sigma \sigma^\top + \sigma^0 (\sigma^0)^\top
+ b_t^\top \Pi_t + \Pi_t b_t - \Pi_t\, H^\top\, \tilde{\Theta}^{-1}\, H\, \Pi_t.
\end{equation*}
with
\begin{align*}
L_t:=&\ b_t  - (\Sigma_t  + \Lambda^1_t) -\Pi_t\, H^\top\, \tilde{\Theta}^{-1}\, H,\\
\hat{L}_t :=&\ (b_t + \bar{b}_t -\Sigma_t - \Lambda^0_t  -\Lambda^1_t),
\end{align*}
and $(\Psi_t)_\T, (\Phi_t)_\T$ the solutions of
\[
\dot{\Psi}_t = L_t \Psi_t, \ \ \Psi_0 = 1 \quad \hbox{and}\ \dot{\Phi}_t = \hat{L}_t\, \Phi_t, \ \ \Phi_0 = 1.
\]
The forward equation is solved by $\hat{n}(t,x,\omega,\gamma,y,z)$ that is Gaussian in $x$ with mean $\bar{\eta}(t,\omega,\gamma,y,z)$ and covariance $\Pi_t$.

\end{thm}

\begin{proof}

In the case of partial information, we can proceed
very similarly as in the case of full information in Section \ref{MFGtypical},
but with a few important modifications.

We can begin by checking (\ref{eqn:consistency}). We have the Gaussian form
$$
\hat{n}(t,x,\omega,\gamma,y,z) = c_t\, e^{-\frac{1}{2}\big(x-\bar{\eta}_t(\omega,\gamma,y,z)\big)\cdot \Pi_t^{-1}\big(x-\bar{\eta}_t(\omega,\gamma,y,z)\big)}.
$$
We see that $x-\sigma^0\, y$ always appears together so that we have the first equation of (\ref{eqn:consistency}). For the second equation of (\ref{eqn:consistency}) we have
\begin{align*}
\partial_z \hat{n}(t,x,\omega,\gamma,y,z) =&\ \big(\partial_z\bar{\eta}(\omega,\gamma,y,z)\big)^\top \Pi_t^{-1} \big(x-\bar{\eta}(\omega,\gamma,y,z)\big)\, \hat{n}(t,x,\omega,\gamma,y,z)\\
=&\ \tilde{\Theta}^{-1} H\big(x-\bar{\eta}(\omega,\gamma,y,z)\big)\hat{n}(t,x,\omega,\gamma,y,z).
\end{align*}

We now make the ansatz (suppressing dependence of $\hat{u}, \hat{\alpha}$, and $\bar{\eta}$ on $(\omega,\gamma,y,z)$ and $\bar{\mu}$ on $(\omega,y)$)
\[
\hat{u}(t,x) = 
\frac{1}{2} \left ( x \cdot \Sigma_t\, x + \bar{\mu}_t \cdot \Sigma^0_t\, \bar{\mu}_t + \bar{\eta}_t \cdot \Sigma^1_t\, \bar{\eta}_t +
2x \cdot \Lambda^0_t\, \bar{\mu}_t  + 2 x \cdot \Lambda^1_t\, \bar{\eta}_t  \right ) 
+ \Delta_t
\]
so that 
\[
\partial_x \hat{u}(t,x)
= \Sigma_t\, x + \Lambda^0_t\, \bar{\mu}_t + \Lambda^1_t\, \bar{\eta}_t
\]
and thus the control of \eqref{eq:partialinfocontrol} is given by
\[
\begin{aligned}
\hat{\alpha}(t) & := - \int_{\RR^d} \partial_x \hat{u}(t,\xi) \eta_t(d \xi)
= - (\Sigma_t  + \Lambda^1_t) \bar{\eta}_t - \Lambda^0_t \bar{\mu}_t.
\end{aligned}
\]
Thus, we have 
\[
dX_t = \left (b_tX_t - (\Sigma_t+\Lambda^1_t)\bar{\eta}_t + (\bar{b}_t - \Lambda^0_t) \bar{\mu}_t \right ) dt + \sigma dW_t + \sigma^0 dW^0_t.
\]
We then take the conditional expectation given 
 $\calF^{\bW^0}_t$ to get
\[
d\bar{\mu}_t = (b_t + \bar{b}_t -\Sigma_t - \Lambda^0_t  -\Lambda^1_t) \bar{\mu}_t \ dt + \sigma^0 dW^0_t.
\]
Letting $\hat{L}_t := (b_t + \bar{b}_t -\Sigma_t - \Lambda^0_t  -\Lambda^1_t)$,
the lifted functional of $\bar{\mu}_t$ has the form
\[
\bar{\mu}(t,\omega, y)
= \Phi_t  \bigg (\bar{\mu}_0 
+ \int_0^t \Phi_s^{-1}\, \hat{L}_s\, \sigma^0\, \omega_s\, ds \bigg )
+ \sigma^0 y, 
\]
where $(\Phi_t)_\T$ is the solution of 
\[
\dot{\Phi}_t = \hat{L}_t\, \Phi_t, \ \ \Phi_0 = 1.
\]

 The equation for $\bar{\eta}_t$ needs to be
derived by computing the first moment
directly from the forward Kushner equation in \eqref{general}, which gives
\begin{align}\label{eqn:kalman_stochastic}
d \bar{\eta}_t
= \left ( (b_t  - (\Sigma_t + \Lambda^1_t) - \Pi_t H^\top \tilde{\Theta}^{-1}H ) \bar{\eta}_t + (\bar{b}_t - \Lambda^0_t) \bar{\mu}_t \right ) dt
+ \sigma^0 dW^0_t
+\Pi_t H^\top \tilde{\Theta}^{-1} dZ_t
\end{align}
where we define the variance $\Pi_t := \int_{\RR^d} \xi^2 \eta_t(d\xi) - \bar{\eta}_t^2$.
If the initial condition is Gaussian, then this quantity is deterministic and classically satisfies
\begin{equation} \label{variance}
\dot{\Pi}_t = \sigma \sigma^\top + b_t^\top \Pi_t + \Pi_t b_t - \Pi_t\, H^\top\, \tilde{\Theta}^{-1}\, H\, \Pi_t.
\end{equation}
This procedure of estimating the state with $\bar{\eta}_t = \mathbb{E}\big[X_t|\mathcal{F}_t^{\mathbf{W}^0,\mathbf{Z}}\big]$ (classically without the presence of the common noise) 
is commonly known as the \emph{Kalman filter} in a discrete time context or as the \emph{Kalman-Bucy filter} in a continuous time context (see the seminal work
Kalman-Bucy \cite{kalman1961new}).  

So relying on this strong consequence of the Gaussian assumption,
we can solve the resulting linear equation for $\bar{\eta}_t$ explicitly. 
More precisely, let 
\[
L_t:= b_t  - (\Sigma_t  + \Lambda^1_t) -\Pi_t\, H^\top\, \tilde{\Theta}^{-1}\, H
\]
and consider the solution $(\Psi_t)_\T$ of
\[
\dot{\Psi}_t = L_t \Psi_t, \ \ \Psi_0 = 1.
\]
Then
\[
\begin{aligned}
\bar{\eta}(t,\omega, \gamma, y,z)
& = \Psi_t\,  \bigg( \bar{\mu}_0 + \int_0^t  \Psi_s^{-1} (\bar{b}_s - \Lambda^0_s) \bar{\mu}(s,\omega, \omega_s) ds - \int_0^t L_s\, \Psi_s^{-1}\, \sigma^0\, \omega_s \, ds  \bigg )
 \\
& \quad + \Psi_t \, \bigg( \int_0^t \Big ( L_s \Psi_s^{-1}\Pi_s + \Psi_s^{-1} \dot{\Pi}_s \big )  H^\top \tilde{\Theta}^{-1} \gamma_s ds    \bigg) + \sigma^0 y + \Pi_t H^\top \tilde{\Theta}^{-1} z
\end{aligned}
\]
where $\dot{\Pi}_t$ is given by \eqref{variance}. The circuitous derivation has not yet justified that $\hat{n}$ satisfies (\ref{liftedsystemgeneral}). From the Gaussian structure we compute
\begin{align*}
D_{t,\omega,\gamma}^{y,z}\hat{n} =&\ -{\rm div}_x \big[(D_{t,\omega}^y\bar{\eta}_t) \hat{n}\big]+\frac{1}{2}{\rm tr}\big[ \dot{\Pi}_t \partial^2_{xx} \hat{n}\big]\\
=&\ -{\rm div}_x \big[ (b_t\bar{\eta}_t  +\bar{b}_t\bar{\mu}_t + \alpha_t)\hat{n}\big]+ \Pi_t H^\top \tilde{\Theta}^{-1}H \bar{\eta}_t \cdot \partial_x\hat{n}\\
&\ +\frac{1}{2}{\rm tr}\Big[\big(\sigma \sigma^\top + b_t^\top \Pi_t + \Pi_t b_t - \Pi_t\, H^\top\, \tilde{\Theta}^{-1}\, H\, \Pi_t\big)\partial_{xx}^2\hat{n}\Big].
\end{align*}
We also have
\begin{align*}
\partial_z\hat{n} =&\ - (\partial_z \bar{\eta}_t)^\top \partial_x \hat{n}\\
=&\ -\tilde{\Theta}^{-1}\, H\, \Pi_t \partial_x \hat{n}.
\end{align*}
Combining this and $\partial_x\hat{n} = -\Pi_t^{-1} (x-\bar{\eta}_t)\hat{n}$ we get
\begin{align*}
D_{t,\omega,\gamma}^{y,z}\hat{n} =&\ -{\rm div}_x \big[ (b_t\, x  +\bar{b}_t\bar{\mu}_t + \alpha_t)\hat{n}\big]-H\, \bar{\eta_t}\cdot \partial_z\hat{n}\\
&\ +\frac{1}{2}{\rm tr}\Big[\sigma \sigma^\top \partial^2_{xx}\hat{n} + \tilde{\theta}\partial^2_{zz} \hat{n}\Big] - H(x-\bar{\eta}_t) \cdot \partial_z \hat{n}.
\end{align*}
Finally, the first equation of (\ref{eqn:consistency}) implies that
\begin{align*}
\frac{1}{2}{\rm tr}\Big[\sigma \sigma^\top \partial^2_{xx}\hat{n}\Big]= \frac{1}{2} \tr\Big[(\sigma\sigma^\top+\sigma^0(\sigma^0)^\top) \ \partial_{xx}^2 \hat{n} + \partial^2_{yy} \hat{n}  +2\sigma^0 \partial_x \partial_y \hat{n}\Big].
\end{align*}
Together these show that $\hat{n}$ solves (\ref{liftedsystemgeneral}).

The lifted value function is given by
\[
\begin{aligned}
\hat{u}(t,x,\omega,\gamma,y,z) & = 
\frac{1}{2}\, {x} \cdot \Sigma_t\, {x} +  {x} \cdot \Lambda^0_t\, \bar{\mu}(t,\omega,y)+\frac{1}{2}\,  \bar{\mu}(t,\omega,y) \cdot \Sigma^0_t\, \bar{\mu}(t,\omega,y)\\
& \quad  +  \frac{1}{2}\, \bar{\eta}(t,\omega,\gamma,y,z) \cdot \Sigma^1_t\, \bar{\eta}(t,\omega,\gamma,y,z) + {x}\cdot  \Lambda^1_t\, \bar{\eta}(t,\omega,\gamma,y,z)  + \Delta_t.
\end{aligned}
\]

Then the compensated HJB equation will take the form
\[
\begin{aligned}
	-&  D_{t,\omega,\gamma}^{y,z} \hat{u}(t,x,\omega,\gamma,y,z)  -\frac{1}{2} \tr\Big[(\sigma\sigma^\top + \sigma^0(\sigma^0)^\top) \ \partial_{xx}^2 \hat{u}(t,x,\omega,\gamma,y,z)   + \partial^2_{yy} \hat{u}(t,x,\omega,\gamma,y,z)  \Big]\\
	& -  \tr[\sigma^0 \partial_x \partial_y \hat{u}(t,x,\omega,\gamma,y,z)  ]  -  \frac{1}{2}  \tr[ \tilde{\Theta} \ \partial_{zz}^2 \hat{u}(t,x,\omega,\gamma,y,z)]  - \partial_z\hat{u}(t,x,\omega,\gamma,y,z) \cdot Hx \\
	& - \partial_x \hat{u}(t,x,\omega,\gamma,y,z) \cdot \Big ( b_t x + \bar{b}_t \bar{\mu}(t,\omega, y) + \hat{K}(t,\omega,\gamma,y,z)\Big )   =  \frac{1}{2} \big | \hat{K}(t,\omega,\gamma,y,z) \big |^2,
\end{aligned}
\]
where 
\[ 
\begin{aligned}
\hat{K}(t,\omega,\gamma,y,z) & := -(\Sigma_t  + \Lambda^1_t) \bar{\eta}(t,\omega,\gamma,y,z) - \Lambda^0_t \bar{\mu}(t,\omega,y)
\end{aligned}
\]
with terminal condition
\[
\hat{u}(T,x,\omega,\gamma,y,z) =  \frac{1}{2} (x - s\bar{\mu}(T,\omega, y))\cdot \bar{q}\, (x - s\bar{\mu}(T,\omega, y)).
\]

Now we may begin computing 
the terms appearing in the lifted functional
backward equation \eqref{liftedsystemgeneral}. 
We first compute a system of equations that provide a lifted form of the Kalman filter equations (\ref{eqn:kalman_stochastic})
\begin{equation}
	\label{liftedkalman}
\begin{aligned}
 D_{t,\omega,\gamma}^{y,z}  \bar{\mu}(t,\omega, y)
&  = 
(b_t + \bar{b}_t - (\Sigma_t  + \Lambda^1_t) - \Lambda^0_t) \bar{\mu}(t,\omega, y) \\
D_{t,\omega,\gamma}^{y,z}  \bar{\eta}(t,\omega,\gamma,y,z) & = L_t \bar{\eta}(t,\omega,\gamma,y,z) 
+ (\bar{b}_t - \Lambda^0_t) \bar{\mu}(t,\omega, y)
\end{aligned}
\end{equation}
\[	
\partial_y \bar{\mu}(t,\omega, y) = \partial_y \bar{\eta}(t,\omega,\gamma,y,z)
 =  \sigma^0, \ \ \ \
 \partial_z \bar{\eta}(t,\omega,\gamma,y,z)
 =   \Pi_t H^\top \tilde{\Theta}^{-1}.
\]
We then compute
\[
\begin{aligned}
 D_{t,\omega,\gamma}^{y,z}  \hat{u}(t,x,\omega,\gamma,y,z) 
& = \frac{1}{2} \left ( x \cdot \dot{\Sigma}_t x + 2x \cdot \dot{\Lambda}^0_t \bar{\mu}(t,\omega, y) + \bar{\mu}(t,\omega, y) \cdot \dot{\Sigma}^0_t \bar{\mu}(t,\omega, y) \right ) \\
& \quad + \frac{1}{2} \bar{\eta}(t,\omega,\gamma,y,z)^\top \dot{\Sigma}^1_t \bar{\eta}(t,\omega,\gamma,y,z)) + x\cdot  \dot{\Lambda}^1_t \bar{\eta}(t,\omega,\gamma,y,z) + \dot{d}_t 
\\
&\quad + (b_t + \bar{b}_t -(\Sigma_t + \Lambda^1_t) - \Lambda^0_t)\bar{\mu}(t,\omega, y) \cdot \left (
 \Lambda^0_t x + \Sigma^0_t \bar{\mu}(t,\omega, y)   \right ) \\
 & \quad + ({\Sigma^1_t}^\top\, \bar{\eta}(t,\omega,\gamma,y,z) + {\Lambda^1_t}^\top\, x)\cdot \left ( L_t \bar{\eta}(t,\omega,\gamma,y,z) + (\bar{b}_t - \Lambda^0_t) \bar{\mu}(t,\omega, y) \right ) 
\end{aligned}
\]
and can further compute
\[
\partial_x \hat{u}(t,x,\omega,\gamma,y,z)  = 
\Sigma_t x 
+ \Lambda^0_t \bar{\mu}(t,\omega,y)+ \Lambda^1_t \bar{\eta}(t,\omega,\gamma,y,z), \ \ \ 
\partial_{xx}^2 \hat{u}(t,x,\omega,y)  = 
\Sigma_t
\]
\[
\partial_y \hat{u}(t,x,\omega,\gamma,y,z) 
= (\sigma^0)^\top  \left ( \Sigma^0_t	 \bar{\mu}(t,\omega,y) + \Sigma^1_t \bar{\eta}(t,\omega,\gamma,y,z)
+ (\Lambda^0_t + \Lambda^1_t)x \right ),
\]
\[
\partial_{yy}^2 \hat{u}(t,x,\omega,\gamma,y,z) 
= (\sigma^0)^\top (\Sigma^0_t+ \Sigma^1_t) \sigma^0,
\]
\[
\partial_z \hat{u}(t,x,\omega,\gamma,y,z) = (\Pi_t H^\top \tilde{\Theta}^{-1})^\top\, ({\Sigma^1_t}^\top\, \bar{\eta}(t,\omega,\gamma,y,z) + {\Lambda^1_t }^\top\, x) ,
\]
\[
\partial_{zz}^2 \hat{u}(t,x,\omega,\gamma,y,z) = (\Pi_t H^\top \tilde{\Theta}^{-1})^\top \Sigma^1_t \Pi_t H^\top \tilde{\Theta}^{-1}
\]
\[
\partial_x \partial_y \hat{u}(t,x,\omega,\gamma,y,z)
=
(\sigma^0)^\top (\Lambda^0_t+\Lambda^1_t).
\]
Inputting these calculations in the compensated equation gives
\[
\begin{aligned}
	&  \frac{1}{2} \left ( x \cdot \dot{\Sigma}_t x + 2x \cdot \dot{\Lambda}^0_t \bar{\mu}(t,\omega, y) + \bar{\mu}(t,\omega, y) \cdot \dot{\Sigma}^0_t \bar{\mu}(t,\omega, y) \right ) \\
& \quad + \frac{1}{2} \bar{\eta}(t,\omega,\gamma,y,z)\cdot \dot{\Sigma}^1_t \bar{\eta}(t,\omega,\gamma,y,z)) + x\cdot  \dot{\Lambda}^1_t \bar{\eta}(t,\omega,\gamma,y,z) + \dot{d}_t 
\\
&\quad + (b_t + \bar{b}_t -\Sigma_t -\Lambda_t^1- \Lambda^0_t)\bar{\mu}(t,\omega, y) \cdot \left (
 \Lambda^0_t x + \Sigma^0_t \bar{\mu}(t,\omega, y)   \right ) \\
 & \quad + ({\Sigma^1_t}^\top\, \bar{\eta}(t,\omega,\gamma,y,z)  + {\Lambda^1_t}^\top\, x) \cdot \left ( L_t \bar{\eta}(t,\omega,\gamma,y,z) + (\bar{b}_t - \Lambda^0_t) \bar{\mu}(t,\omega, y) \right )  \\
 & \quad +\frac{1}{2}  \tr[(\sigma \sigma^\top + \sigma^0 (\sigma^0)^\top)  \ \Sigma_t  ] + \frac{1}{2}\tr[\sigma^0 (\sigma^0)^\top (\Sigma^0_t+ \Sigma^1_t)]   +  \tr[\sigma^0 (\sigma^0)^\top (\Lambda^0_t+\Lambda^1_t)  ]  +  \frac{1}{2}  \tr[ \ \Pi_t H^\top H \Pi_t^\top \tilde{\Theta}^{-1} \Sigma^1_t  ]  \\
 &\quad + ({\Sigma^1_t}^\top\, \bar{\eta}(t,\omega,\gamma,y,z) + {\Lambda^1_t}^\top\, x )\cdot \Pi_t\, H^\top \tilde{\Theta}^{-1} H\, x \\
	&\quad + \Big ( \Sigma_t x 
+ \Lambda^0_t \bar{\mu}(t,\omega,y)+ \Lambda^1_t \bar{\eta}(t,\omega,\gamma,y,z) \Big ) \cdot \Big ( b_t x + (\bar{b}_t- \Lambda^0_t) \bar{\mu}(t,\omega, y) - (\Sigma_t + \Lambda^1_t) \bar{\eta}(t,\omega,\gamma,y,z) \Big ) \\
	 & \quad    + \frac{1}{2} \left(  (x - s \bar{\mu}(t,\omega, y))\cdot \bar{q}\, (x - s \bar{\mu}(t,\omega, y)) \right) + \frac{1}{2} \big | (\Sigma_t  + \Lambda^1_t)\bar{\eta}(t,\omega,\gamma,y,z) + \Lambda^0_t \bar{\mu}(t,\omega,y) \big |^2 = 0.
\end{aligned}
\]
We now collect terms (symmetrizing for the squared terms) to arrive at the following closed system of Riccati equations (note we anticipate the coefficient of ``$\bar{\mu}\bar{\eta}$'' is $0$):
\begin{itemize}
	\item $|x|^2 \ : \ \dot{\Sigma}_t =  - \Sigma_t^\top b_t - b_t^\top \Sigma_t    + \Lambda^1_t\, \Pi_t\, H^\top\, \tilde{\Theta}^{-1}\, H+ H^\top\, \tilde{\Theta}^{-1}\, H\, \Pi_t\, {\Lambda_t^1}^\top  , \ \ \ \Sigma_T = \bar{q} $ 
	\item $x \bar{\mu} \ : \ \dot{\Lambda}_t^0 = -{\Lambda^0_t}^\top\, (b_t + \bar{b}_t -\Sigma_t-\Lambda_t^1-\Lambda^0_t ) -{\Lambda^1_t}\, (\bar{b}_t - \Lambda^0_t)-\Sigma_t^\top\, (\bar{b}_t - \Lambda^0_t)- b_t^\top \Lambda^0_t  , $
	\subitem $\quad \quad \quad \quad \Lambda^0_T = -  \bar{q}\, s  $
	\item $\bar{\mu}^2 \ \  :  \ \dot{\Sigma}^0_t =  -( b_t + \bar{b}_t-\Sigma_t -\Lambda_t^1-\Lambda^0_t )^\top \Sigma^0_t-{\Sigma_t^0}^\top\, ( b_t + \bar{b}_t-\Sigma_t-\Lambda_t^1 -\Lambda^0_t ) - {\Lambda^0_t}^\top\, ( \bar{b}_t-\Lambda^0_t)$\\
	\subitem$\quad \quad \quad \quad  -( \bar{b}_t-\Lambda^0_t)^\top \Lambda^0_t-{\Lambda_t^0}^\top\, \Lambda_t^0$, $\quad \quad \quad \quad \Sigma^0_T = s^\top \bar{q} s $
	\item $1 \ \ \ \ : \   \dot{\Delta}_t =  - \frac{1}{2}\tr[(\sigma \sigma^\top + \sigma^0 (\sigma^0)^\top) \ \Sigma_t ] - \frac{1}{2} \tr[  \sigma^0 (\sigma^0)^\top (\Sigma^0_t + \Sigma^1_t)  ] 
	- \tr[\sigma^0 (\sigma^0)^\top (\Lambda^0_t + \Lambda^1_t)]$
	\subitem $ \quad \quad \quad \quad  - \frac{1}{2}  \tr[ \ \Pi_t H^\top H \Pi_t^\top \tilde{\Theta}^{-1} \Sigma^1_t  ] , \ \ \ \ \   \Delta_T = 0   \  \ $
	\item $\bar{\eta}^2 \ : \ \dot{\Sigma}^1_t = - \left ( L_t^\top \Sigma^1_t +  (\Sigma^1_t)^\top L_t \right )  + \left ( (\Sigma_t + \Lambda^1_t)^\top \Lambda^1_t + ( \Lambda^1_t)^\top (\Sigma_t + \Lambda^1_t) \right ) - (\Sigma_t + \Lambda^1_t)^\top (\Sigma_t + \Lambda^1_t),$
	\subitem $\quad \quad \quad \quad \quad \Sigma^1_T = 0   $
	\item $x \bar{\eta} \ : \ \dot{\Lambda}^1_t = -  \Lambda^1_t\, L_t -  H^\top\, \tilde{\Theta}^{-1}\, H\, \Pi_t^\top\, {\Sigma^1_t}^\top + \Sigma_t^\top (\Sigma_t + \Lambda^1_t) - b_t^\top \Lambda^1_t $, $\quad \quad  \Lambda^1_T = 0 $,
	\item $  \bar{\mu} \bar{\eta} \ : \  0=-\Sigma^1_t (\bar{b}_t - \Lambda^0_t)  - \Lambda^1_t (\bar{b}_t - \Lambda^0_t) + \Sigma_t \Lambda^0_t - \Sigma_t \Lambda^0_t + (\Sigma^1_t + \Lambda^1_t) (\bar{b}_t - \Lambda^0_t) $.
\end{itemize}

\end{proof}

\subsection*{Discussion of Solvability of the Riccati Equations}

Notice that the equations for $\Sigma_t, \Lambda^0_t$ are quadratic
Riccati equations, while the equation for $\Sigma^0_t$
is linear.
Now observe that 
\[
L_t+ \Pi_t H^\top \tilde{\Theta}^{-1} H = b_t  - (\Sigma_t + \Lambda^1_t).
\]
Hence, if we add the $\Sigma_t^1$ and ${\Lambda^1_t}^\top$
equations together, we see that $\Upsilon_t := \Sigma_t^1 + {\Lambda^1_t}^\top$ solves
\[
\begin{aligned}
\dot{\Upsilon}_t & = -  L_t^\top \Upsilon_t - {\Sigma^1_t}^\top(b_t  - (\Sigma_t + \Lambda^1_t)) + (\Sigma_t + \Lambda^1_t)^\top\, \Sigma_t - {\Lambda^1_t}^\top\, b_t  \\
& \quad + \left ( (\Sigma_t + \Lambda^1_t)^\top \Lambda^1_t + ( \Lambda^1_t)^\top (\Sigma_t + \Lambda^1_t) \right ) - (\Sigma_t + \Lambda^1_t)^\top (\Sigma_t + \Lambda^1_t) \\
&= -  L_t^\top \Upsilon_t + \Upsilon_t\, (\Sigma_t + \Lambda^1_t) -\Upsilon_t\, b_t.
\end{aligned}
\]
Since the terminal condition is $\Upsilon_T = \Sigma^1_T + {\Lambda^1_T}^\top = 0$,  we have $\Upsilon_t \equiv 0$
and thus $\Sigma^1_t = - (\Lambda^1_t)^\top = -\Lambda_t^1$.
Note this implies that the coefficient
for the $\bar{\mu} \bar{\eta}$ term disappears,
thus reducing the problem to the case of the first six
equations.

To solve these equations, we can look at $\Gamma_t := \Sigma_t - \Sigma^1_t$
and sum together the first and sixth equation 
along with $\Lambda^1_t = -\Sigma^1_t$ and the definition of $L_t$ to get
\[
\dot{\Gamma}_t = \Gamma_t^\top \Gamma_t - b_t^\top \Gamma_t - \Gamma_t^\top b_t, \ \ \ \Gamma_T = \bar{q}.
\]
which is the same quadratic Riccati equation as 
for ``$\Gamma_t$'' in the case of full information 
in Section \ref{MFGtypical}.
Similarly, the equation for $\Lambda_t^0$ can be rewritten as:
\begin{align*}
\dot{\Lambda}_t^0 =&\ -{\Lambda^0_t}^\top\, (b_t + \bar{b}_t -\Sigma_t-\Lambda_t^1-\Lambda^0_t ) -{\Lambda^1_t}\, (\bar{b}_t - \Lambda^0_t)-\Sigma_t^\top\, (\bar{b}_t - \Lambda^0_t)- b_t^\top \Lambda^0_t  \\
=& (\Lambda^0_t)^\top \Lambda^0_t -(\Lambda^0_t)^\top\, (b_t + \bar{b}_t - \Gamma_t) + \Gamma_t^\top  \Lambda^0_t  - \Gamma_t^\top\, \bar{b}_t  - b_t^\top\, \Lambda_t^0 ,
\end{align*}
which is the same quadratic Riccati equation as 
for ``$\Lambda^0_t$'' in the case of full information 
in Section \ref{MFGtypical}. 

Now, following the approach in Section \ref{MFGtypical}, we can consider $\tilde{\Lambda}_t = \Gamma_t + \Lambda_t^0$, which satisfies:
$$
	\dot{\tilde{\Lambda}}_t = \tilde{\Lambda}_t^\top\tilde{\Lambda_t} - b_t^\top\, \tilde{\Lambda}_t - \tilde{\Lambda}_t^\top\, (b_t+\bar{b}_t).
$$
Once again, under the conditions that $\bar{q}$ is symmetric and positive semidefinite, $s\leq 1$, and that $\bar{b}_t$ is a scalar times the identity, we see that $\tilde{\Lambda}_t$ is also symmetric and positive semidefinite and thus a unique global solution exists.

Finally and most importantly,
the above manipulations embody the \emph{separation principle}:
to go from the optimal feedback function \eqref{feedbackfullinfo} in the case of full information
in Section \ref{MFGtypical} to the case of partial information here,
one just needs to replace the state with the best guess of the state given 
the common noise $\bW^0$ and the
partial observation $\bZ$.
This is exactly what we have just established. 
\begin{thm*}[Separation Principle]
	The optimal feedback control for a mean field game with common noise and a partial information constraint in the linear-quadratic framework with Gaussian initial condition has the linear feedback form
	\begin{equation}
		\label{kalmancontrol}
	\hat{\alpha}(t,\eta_t,\partial_x u_t) = - \int_{\RR^d} \partial_x u_t(\xi) \eta_t(d \xi) =  -\Gamma_t\, \bar{\eta}_t - \Lambda^0_t \, \bar{\mu}_t,
	\end{equation}
	where the coefficients $\Gamma_t$ and $\Lambda^0_t$
	satisfy the same equations as for the optimal feedback function
	\[
		\alpha^*(t,x) := -\Gamma_t\, x - \Lambda^0_t \, \bar{\mu}_t
	\] 
	in the case of full information.\footnote{We remind the reader that we also checked the consistency of these equations with the literature at the end of Section \ref{MFGtypical}.}
\end{thm*}
Thus, the optimal control is determined ``separately'' from the partial observation in that
the latter only enters in the former through the conditional expectation $\bar{\eta}_t = \mathbb{E}\big[X_t|\mathcal{F}_t^{\mathbf{W}^0,\mathbf{Z}}\big]$, which solves the so-called \emph{Kalman filtering problem} (again, see the seminal work \cite{kalman1961new} of Kalman-Bucy).

\subsubsection*{Comparison with the literature}

The recent article \cite{benyam} of Bensoussan-Yam has demonstrated a close connection between the partial observation control problem and mean field theory; more precisely, they use a master equation approach for the linear quadratic partial information problem without mean field interactions.  This approach allows them to prove a \emph{separation principle}, i.e., that the optimal control is a linear feedback of the expected state given the observation process, but without requiring the standard simplifying assumption that the initial distribution is Gaussian (a significant assumption that our calculations of Section \ref{MFGpartialinfosection} notably rely on).  
It is noted that the complications from non-Gaussian initial conditions only arise in the Kalman filter equations to determine the distribution conditioned on the observations, whereas the fact one arrives at a linear feedback control \eqref{kalmancontrol} should not change.  

In the context of our paper, the Kalman filter corresponds to the mean flow $(\bar{\eta}_t)_\T$
of the solution $\boldsymbol{\eta} = (\eta_t)_\T$ 
to the forward Kushner equation in \eqref{general}.  
In the Gaussian case with linear-quadratic data, 
when computing the covariance $\Pi_t$ from the Kushner equation, 
a term involving the third moment naturally arises, but this can
be expressed in terms of the second moment, 
thus leading to the deterministic Riccati equation \eqref{variance} for the covariance $\Pi_t$.  
 It is not yet clear to the authors whether the approach of \cite{benyam} could be adapted
 to the mean field game problem with partial information to generalize the solution outside of the case of Gaussian initial conditions.

 \pagebreak

 \appendix

\section{Derivation of the lifted functional approach}\label{apx:derivation}

For the sake of simplicity, we take $\sigma = \sigma^0 = I_d$
and a quadratic separable Hamiltonian $H(t,x,p,\mu) := \frac{1}{2}|p|^2 - f(t,x,m)$ when $m$ is the density of $\mu$.
We then perform the change of variables $x \mapsto x-y$
so that we can reduce the form of the lifted functional system \eqref{liftedsystemintro} to
finding a pair $(\widetilde{u}(t,x,\omega,y), \widetilde{r}(t,x,\omega))$
satisfying 
\begin{equation} \label{liftedsystemsimplified}
\begin{cases}
 -\partial_t \widetilde{u}(t,x,\omega,y) - \frac{1}{2}  ( \Delta_x + \Delta_y  ) \widetilde{u}(t,x,\omega,y) + \frac{1}{2} | \partial_x \widetilde{u}(t,x,\omega,y)|^2  - \widetilde{f}(t,x,y,\widetilde{r}(t,\cdot, \omega))  =  \calD_\omega^y \widetilde{u}(t,x,\omega,y),  \\
		  \partial_t  \widetilde{r}(t,x,\omega)
 = \frac{1}{2} \Delta_x \widetilde{r}(t,x,\omega)   + \text{div}_x \left [\widetilde{r}(t,x,\omega)\partial_x \widetilde{u}(t,x,\omega,\omega_t) \right ],\\
\widetilde{u}(T,x,\omega,y) = \widetilde{g}(x,y,\widetilde{r}(T,\cdot,\omega)), \ \ \ 
 \widetilde{r}(0,x,\omega) = \ell(x), 
\end{cases}
\end{equation}
where 
\[
\widetilde{f}(t,x,y,\rho):= f(t,x+y,\rho(\cdot - y)), \ \ 
\widetilde{g}(x,y,\rho):= g(x+y,\rho(\cdot - y)).
\]
To make explicit the connection, once we find a solution to \eqref{liftedsystemsimplified}, 
we immediately recover a solution to the original system \eqref{liftedsystemintro}
by setting
\[
\hat{u}(t,x,\omega,y) := \widetilde{u}(t,x-y,\omega,y), \ \ 
\hat{m}(t,x,\omega,y):= \widetilde{r}(t,x-y,\omega).
\]

Now let $\bB = (B_t)_\T$ be a $d$-dimensional Brownian motion independent
of $\bW$.
Also, given
a path $\omega \in \Omega$,
we write
\begin{equation}
\label{rightcont_apx}
W^{t,\omega,y}_s := \omega_s 1_{[0,t)}(s) + [y + W_s - W_t ] 1_{[t,T]}(s).
\end{equation}
We consider the candidate solution of the compensated backward HJB of \eqref{liftedsystemsimplified} given by the BSDE representation
\begin{equation}
	\label{liftedsolnofPDHJB}
	\widetilde{u}(t,x,\omega,y) := Y^{t,x,\omega,y}_t, \ \ (t,x,\omega,y) \in [0,T] \times \RR^d 
	\times \Omega \times \RR^d,
\end{equation}
where for each $(t,x,\omega,y) \in [0,T) \times \RR^d \times \Omega \times \RR^d$, the 
triple $(Y^{t,x,\omega,y}_s, (Z^{t,x,\omega,y}_s, \Gamma_s^{t,x,\omega,y}))_{t \leq s\leq T}$
satisfies the ``lifted'' BSDE (see Peng \cite{peng2011note})
\begin{equation}
	\label{liftedBSDE}
	\begin{aligned}
	Y^{t,x,\omega,y}_s  & = \widetilde{g}(B^{t,x}_T, W^{t,\omega,y}_T, \widetilde{r}\big(T,\cdot, \bW^{t,\omega,   y }) \big) 
	+ \int_s^T \widetilde{f}(\theta,B^{t,x}_\theta, W^{t,\omega,y}_\theta, \widetilde{r}(\theta, \cdot, \bW^{t,\omega,   y } ) ) d\theta \\
	& \quad - \frac{1}{2} \int_s^T |Z^{t,x,\omega,y}_\theta|^2 d\theta  - \int_s^T \left [Z^{t,x,\omega,y}_\theta \cdot dB_\theta + \Gamma_\theta^{t,x,\omega,y} \cdot dW_\theta   \right ].
	\end{aligned}
\end{equation}

We now \emph{sketch} how to go 
from the candidate solution $\widetilde{u}(t,x,\omega,y) := Y^{t,x,\omega,y}_t$
as in \eqref{liftedsolnofPDHJB} to the form of compensated backward HJB
of the system \eqref{liftedsystemsimplified}.
First, it is readily seen that the concatenated path $W^{t,\omega,y}_\cdot$ still satisfies the flow property: 
\[
W^{s,W^{t,\omega,y}_\cdot,W^{t,\omega,y}_s}_r = W^{t,\omega,y}_r, \ \ \text{for} \ \  t \leq s \leq r \leq T.
\]
This in turn will ensure we have
a corresponding flow property at the level of the BSDE:
\begin{equation}
	\label{BSDEflowproperty}
\widetilde{u}(t+\epsilon,B^{t,x}_{t+\epsilon},\bW^{t,\omega,y},W^{t,y}_{t+\epsilon}) = Y^{t,x,\omega,y}_{t+\epsilon}.
\end{equation}
Arguing as in Theorem 3.2 of Pardoux-Peng \cite{pardoux1992backward}, this leads us to consider 
the decomposition
\begin{equation}
	\label{liftedflowproperty}
	\begin{aligned}
		\widetilde{u}(t+\epsilon,x,\omega,y) - \widetilde{u}(t,x,\omega,y) & = [\widetilde{u}(t+\epsilon,B^{t,x}_{t+\epsilon},\bW^{t,\omega,y},W^{t,y}_{t+\epsilon}) - \widetilde{u}(t,x,\omega,y) ] \quad \text{(1st difference)} \\
		& + [\widetilde{u}(t+\epsilon,x,\omega,y) - \widetilde{u}(t+\epsilon,B^{t,x}_{t+\epsilon},\omega, W^{t,y}_{t+\epsilon}) ] \quad \text{(2nd difference)} \\
		& +  [\widetilde{u}(t+\epsilon,B^{t,x}_{t+\epsilon},\omega, W^{t,y}_{t+\epsilon}) - \widetilde{u}(t+\epsilon,B^{t,x}_{t+\epsilon}, \bW^{t,\omega,y}, W^{t,y}_{t+\epsilon}) ]   \quad \text{(3rd difference)}
	\end{aligned}
\end{equation}
The first difference in \eqref{liftedflowproperty} can be expressed in terms of the BSDE by the flow property \eqref{BSDEflowproperty}, and thus upon dividing by $\epsilon>0$, taking expectations, and letting $\epsilon \to 0$, 
it will contribute the term ``$\frac{1}{2}| \partial_x \widetilde{u}(t,x,\omega)|^2 - \widetilde{f}(t,x,y,\widetilde{r}(t,\cdot, \omega))$''.
Next, by an application of the completely classical It\^o formula, the second
difference in \eqref{liftedflowproperty} will contribute ``$-\frac{1}{2}(\Delta_x + \Delta_y)\widetilde{u}(t,x,\omega,y)$''. 
Finally, for the third difference of \eqref{liftedflowproperty}, write $X(\epsilon):= B^{t,x}_{t+\epsilon}$, $Y(\epsilon):= W^{t,y}_{t+\epsilon}$, and
\[
H^{t,\omega,y}_s(\epsilon) := [y-\omega_s+W_s - W_t]1_{[t,t+\epsilon)}(s).
\]
Then the third difference of \eqref{liftedflowproperty} can be rewritten as
\[
\widetilde{u}(t+\epsilon,X(\epsilon),\omega, Y(\epsilon)) - \widetilde{u}(t+\epsilon,X(\epsilon),\omega+H^{t,\omega,y}(\epsilon), Y(\epsilon))
\]
where $X(\epsilon),Y(\epsilon) \to x,y$ as $\epsilon \to 0$. Hence, up to stochastic arguments that 
will not contribute given suitable joint regularity, we identify the compensator \eqref{compensator} as the limit of the final difference in \eqref{liftedflowproperty}:
\begin{equation}
	\label{compensator2} 
	\lim_{\epsilon \to 0} \epsilon^{-1} \EE \left [ \widetilde{u}(t+\epsilon,X(\epsilon),\omega+H^{t,\omega,y}(\epsilon), Y(\epsilon)) - \widetilde{u}(t+\epsilon,X(\epsilon),\omega, Y(\epsilon))   \right ] = \calD_\omega^y \widetilde{u}(t,x,\omega,y).
\end{equation}

\section{Unique extensions to paths with a single jump}\label{apx:unique_extensions_to_paths_with_a_single_jump}

There are various ways to deal with the various
forthcoming derivatives 
on path space. 
A common 
approach in the path dependent PDE literature is to work on 
Skorokhod space as we have done in the main body of this paper. An alternative is to adapt the more recent seminorm topology approach of Section 2.2 from Cosso-Russo \cite{cosso2016functional} to our setting.
It formalizes the notion of a path dependent functional being ``not too sensitive to a possible jump nearby any given fixed time $t \in [0,T]$.''

Fix $t \in [0,T]$. Then for each fixed $M>0$,
consider the space $\calC_{t,M}([0,T];\RR^d)$
of paths bounded by $M$ and continuous on $[0,T]$ except for
possibly a jump at time $t$, where the path has a limit from the left and is continuous from the right (c\`{a}dl\`{a}g).
Endow $\calC_{t,M}([0,T];\RR^d)$ with the topology
associated to the metric\footnote{Note here we use a more standard looking metric since the restriction to bounded paths allows us to avoid the arguably more abstract Frechet-type metric construction ``$ \sum_{k=1}^\infty 2^{-k} \frac{[\omega-\eta]_{t,k}}{1+[\omega-\eta]_{t,k}}$,'' which does not appear as good for checking estimates. } 
\begin{equation}
\label{metricseminorm}
\bd_{t}(\omega,\eta):= \sum_{k=1}^\infty 2^{-k} [\omega-\eta]_{t,k},
\end{equation}
induced by an increasing countable family of seminorms of the form
\begin{equation}
\label{seminorms}
[\omega]_{t,k} := \sup_{s\in [0, \max\{t-2^{-k}, 0\}]\cup [t,T] } |\omega_s|.
\end{equation}

Then 
consider the space $\calC_{t}([0,T];\RR^d):= \cup_{M>0} \, \calC_{t,M}([0,T];\RR^d)$ endowed 
with the smallest topology such that all the inclusions $\calC_{t,M}([0,T];\RR^d) \hookrightarrow \calC_{t}([0,T];\RR^d)$ are continuous.\footnote{We remark that this ``inductive topology'' on $\calC_{t}([0,T];\RR^d)$ is not metrizable.}
More concretely, $\eta^N$ converges to $\eta$ in $\calC_{t}([0,T];\RR^d)$
if there is an $M>0$ such that $\sup_{s\in [0,T]} |\eta^N(s)| \leq M$ for all $N$, and for all $k \geq 1$, $[\eta^N - \eta]_{t,k} \to 0$ as $N \to \infty$. 

Intuitively, for each $t \in [0,T]$,
paths in $\calC_{t}([0,T];\RR^d)$
are continuous other than 
a possible jump
discontinuity at time $t$,
and convergent sequences of such paths cannot form arbitrarily large jumps.

Finally, let $\CC_T$
denote the set 
of all
functionals $\psi: \Omega \to \RR$ on $\Omega = C([0,T]; \R^d)$ uniformly 
continuous with respect to $\bd_t$ of \eqref{metricseminorm}
for all $t \in [0,T]$. Despite $\psi \in \CC_T$ only being defined on continuous paths $\Omega$,
the stronger 
continuity property of $\psi$ being 
a member of $\CC_T$ ensures that 
$\psi$ admits a unique extension to each $\calC_{t,M}([0,T];\RR^d)$ for any $M>0$, and thus to $\calC_{t}([0,T];\RR^d)$, for any $t \in [0,T]$.\footnote{
	This stronger continuity assumption for functionals $\phi \in \CC_T$ of continuous paths can also be shown to be compatible with the general Arzela-Ascoli criterion (Theorem 47.1 of Munkres \cite{munkres2000topology}), which would be 
	convenient for fixed point arguments
	for the main lifted functional mean field game system \eqref{liftedsystemintro} below.
}

\section{Zakai-Stratonovich equation}\label{apx:Zakai}
We now assume that $h(t,x,\mu_t) = h(x)$, i.e., independent of $t$ and $\mu_t$. 
To begin, we first trade the
nonlinear forward Kushner equation of \eqref{general}
for its unnormalized counterpart, the so-called Zakai-Stratonovich equation:
\begin{equation} \label{zakai}
	\begin{cases}
		d q_t(x) = \bigg ( \frac{1}{2} \tr \left [ \partial_{xx}^2 \left ( \sigma \sigma^\top  q_t(x) \right ) \right ] - \text{div}_x \left [q_t(x) b(t,x, \alpha(t,\mu_t,\eta_t,\partial_x u_t),\mu_t)\right ]  - q_t(x) \frac{1}{2} h(x)^\top \tilde{\Theta}^{-1} h(x)   \bigg )dt \\
		\quad  \quad  \quad \quad \quad - \text{div}_x \big [ q_t(x) \sigma^0 \circ  dW^0_t  \big ] +  q_t(x)h(x)^\top\tilde{\Theta}^{-1} \circ dZ_t   \\
		q_0(x) \in L^1_+(\RR^d), \ \ \eta_t(dx) := \frac{q_t(x)}{\int_{\RR^d} q_t(\xi) d\xi } dx,
	\end{cases}
\end{equation}
where $\circ$ denotes Stratonovich integration.

The point of the flow transformation method is to remove the noises in the above equation
via a suitable change of variables, thus reducing its solution
to a more classical, albeit random, PDE. 
Following Section 3.4.2 of Souganidis \cite{souganidis2019pathwise}, 
we will look for a solution of the form $q_t(x) = S(t)w_t(x)$,
where $S(t)$ is the solution map of the linear equation
\[
d \mathfrak{p}_t(x) =  - \text{div}_x \big [ \mathfrak{p}_t(x) \sigma^0 \circ  dW^0_t  \big ] +  \mathfrak{p}_t(x) h(x)^\top\tilde{\Theta}^{-1} \circ dZ_t.
\] 
More explicitly, this solution map is given by
\begin{equation} \label{solutionmap}
S(t)f(x) = f(x-\sigma^0W^0_t) \exp \left ( \int_0^t h(x+\sigma^0(W^0_s - W^0_t))^\top\tilde{\Theta}^{-1} \circ dZ_s \right )
\end{equation}
and thus
\[
S^{-1}(t)f(x) = f(x+\sigma^0W^0_t) \exp \left ( - \int_0^t h(x+\sigma^0W^0_s)^\top\tilde{\Theta}^{-1} \circ dZ_s \right ).
\]
Now define
\[
\begin{aligned}
F(X,p,u,x,t) & :=  \frac{1}{2} \tr \left [\sigma \sigma^\top  X  \right ] - p \cdot K_t(x) + u \cdot \left ( \text{div}_x  K_t(x)  - \frac{1}{2} h(x)^\top \tilde{\Theta}^{-1} h(x) \right )
\end{aligned}
\]
where here we employ the generic notation $K_t(x) := b(t,x,\mu_t, \alpha(t,\mu_t,\eta_t,\partial_x u_t))$. 
Then Section 3.4.2 of Souganidis \cite{souganidis2019pathwise} shows that $w_t(x)$ is a solution of the random PDE
\begin{equation}
	\label{flowtransformation}
\begin{aligned}
\partial_t & w_t(x) = S^{-1}(t) F\big(\partial_{xx}^2 [S(t)w_t(x)] , \partial_x [S(t)w_t(x)] , S(t)w_t(x), x,t \big).
\end{aligned}
\end{equation}
We thus have a functional dependence of the form
\[
w_t(x) = \widetilde{w}(t,x,\bW^0,\bZ) :=  \widetilde{w}(t,x,(W^0_s)_{0\leq s < t},(Z_s)_{0\leq s < t});
\]
in particular, this line of reasoning shows how $w_t(x)$ can be expressed as a functional of the paths  of the strict prior history of the noises.
We will write out this dependence more explicitly in the system \eqref{liftedsystemgeneral} below,
where we will expand out the equation \eqref{flowtransformation}.

Now further assume $h(x) = Hx$. Then the normalized measure takes the form
\[
n_t(x) = 
\frac{q_t(x)}{\int_{\RR^d} q_t(\xi) d\xi}
= 
\frac{w_t(x-\sigma^0W^0_t) \exp \left(Hx \cdot \tilde{\Theta}^{-1}Z_t \right) }{\int_{\RR^d} w_t(\xi-\sigma^0W^0_t) 
\exp \left( H\xi \cdot \tilde{\Theta}^{-1}Z_t  \right) d\xi},
\]
where it is significant that the stochastic integrals
arising from the solution map $S(t)$ of \eqref{solutionmap}
have canceled out.
Indeed, now that these stochastic integrals are gone, 
we can conclude that $n_t(x)$ admits the lifted functional
representation $\hat{n}(t,x,\bW^0,\bZ, W^0_t,Z_t)$,
where
\[
\hat{n}(t,x,\omega,\gamma,y,z) = 
\frac{\widetilde{w}(t,x-\sigma^0y,\omega,\gamma) \exp \left(Hx \cdot \tilde{\Theta}^{-1} z \right) }{\int_{\RR^d} \widetilde{w}(t,\xi-\sigma^0y,\omega,\gamma) 
\exp \left( H\xi \cdot \tilde{\Theta}^{-1} z  \right) d\xi}.
\]

Altogether, we expect the lifted functional form of the solution quintuple of the system \eqref{general} to be
\begin{align*}
&
u_t(x) = \hat{u}(t,x,\bW^0,\bZ,W^0_t,Z_t), \ \ \ v_t(x) = (\partial_y \hat{u})(t,x,\bW^0,\bZ,W^0_t,Z_t), \ \ \ k_t(x) = (\partial_z \hat{u})(t,x,\bW^0,\bZ,W^0_t,Z_t),
\\
&
 m_t(x) = \hat{m}(t,x,\bW^0,W^0_t) \ \ \ 
\text{with} \ \ \ \hat{m}(t,x,\omega,y) = \widetilde{r}(t,x-\sigma^0 y, \omega),
\\
&
 n_t(x) = \hat{n}(t,x,\bW^0,\bZ,W^0_t,Z_t) \ \ \ \text{with} \ \ \ \hat{n}(t,x,\omega,\gamma,y,z) = \frac{\widetilde{w}(t,x-\sigma^0y,\omega,\gamma) \exp \left (Hx \cdot \tilde{\Theta}^{-1} z \right ) }{\int_{\RR^d} \widetilde{w}(t,\xi-\sigma^0y,\omega,\gamma) 
\exp \left ( H\xi \cdot \tilde{\Theta}^{-1} z  \right ) d\xi},
\end{align*}
where the triple $(\hat{u}(t,x,\omega,\gamma,y,z), \widetilde{r}(t,x,\omega), \widetilde{w}(t,x,\omega,\gamma))$ solves the following system of equations:
\begin{equation} \label{liftedsystemgeneral2}
	\begin{cases}
		-D_{t,\omega,\gamma}^{y,z}\hat{u}(t,x,\omega,\gamma,y,z) -\frac{1}{2} \bigg ( \tr[(\sigma \sigma^\top + \sigma^0 (\sigma^0)^\top) \ \partial_{xx}^2 \hat{u}(t,x,\omega,\gamma,y,z)  ] + \Delta_y \hat{u}(t,x,\omega,\gamma,y,z)    \bigg ) \\
	\quad \quad  -  \tr[\sigma^0 \partial_x \partial_y \hat{u}(t,x,\omega,\gamma,y,z)  ]  -  \frac{1}{2}  \tr[ \tilde{\Theta} \ \partial_{zz}^2 \hat{u}(t,x,\omega,\gamma,y,z)]   - (\partial_z\hat{u})(t,x,\omega,\gamma,y,z) \cdot Hx   \\
		\quad \quad - f\Big (t,x, \hat{m}(t,\cdot,\omega,y), \ \hat{\alpha}(t,x,\omega,\gamma,y,z), \partial_x \hat{u}(t,\cdot,\omega,\gamma,y,z)) \Big)   = 0,
		\\
\partial_t \widetilde{w}(t,x,\omega,\gamma) = \frac{1}{2}\tr[\sigma \sigma^\top \partial_{xx}^2 \widetilde{w}(t,x,\omega,\gamma)] + \tr[\sigma \sigma^\top \partial_x \widetilde{w}(t,x,\omega,\gamma) \gamma_t^\top \tilde{\Theta}^{-1}H  ] \\
\quad \quad  \quad  \quad  \quad \quad \quad - \text{div}_x [ \widetilde{w}(t,x,\omega,\gamma)   \hat{b}(t,x,\omega,\gamma)]    - \widetilde{w}(t,x,\omega,\gamma) \,
\hat{b}(t,x+\sigma^0 \omega_t,\omega,\gamma) \cdot H^\top \tilde{\Theta}^{-1} \gamma_t  \\
\quad \quad  \quad  \quad  \quad \quad \quad + \widetilde{w}(t,x,\omega,\gamma) \left ( \frac{1}{2} \tr[\sigma \sigma^\top  H^\top \tilde{\Theta}^{-1} \gamma_t \gamma_t^\top \tilde{\Theta}^{-1} H ]  - \frac{1}{2} (x+\sigma^0\omega_t)^\top H^\top \tilde{\Theta}^{-1} H (x+\sigma^0\omega_t)    \right ), \\
\partial_t  \widetilde{r}(t,x,\omega)
 = \frac{1}{2} \tr  \big [ \sigma \sigma^\top  \partial_{xx}^2  \widetilde{r}(t,x,\omega) \big ]  + \text{div}_x \left [\widetilde{r}(t,x,\omega) \bar{b}(t,x+\sigma^0\omega_t,\omega) \right ], \\
		\hat{u}(t,x,\omega,\gamma,y,z) = g(x, \hat{m}(t,\cdot,\omega,y)), \ \ \ 
		\widetilde{w}(0,x,\omega,\gamma)=\widetilde{r}(0,x,\omega) = \ell(x).		
	\end{cases}
\end{equation}
Here, the ``$\hat{\alpha}$'' is given by
\[ 
\begin{aligned}
& \hat{\alpha}(t,x,\omega,\gamma,y,z) \\
& := \argmin_{\alpha \in \RR^d} \left \{  \int_{\RR^d} \hat{n}(t,\xi,\omega,\gamma,y,z) \left [f(t,\xi,\alpha,\hat{m}(t,\cdot,\omega,y))+ \partial_x \hat{u}(t,\xi,\omega,\gamma,y,z) \cdot b(t,\xi,\alpha,\hat{m}(t,\cdot,\omega,y))  \right ] d\xi  \right \}.
\end{aligned}
\]
We abbreviate the lifted drift by
\begin{align*}
&\hat{b}(t,x,\omega,\gamma) :=b(t,x,\hat{m}(t,\cdot,\omega,\omega_t),\hat{\alpha}(t,x,\omega,\gamma,\omega_t,\gamma_t)),
\end{align*}
and the mean drift is given by, using $\bp = (\partial_x \hat{u}(t,\bW^0,\bZ,W^0_t,Z_t))_{\T}$,
\[
\bar{b}(t,x,\omega):= 
\tilde{\EE}^{\bmu,\bp} \left [\hat{b}(t,x,\cdot,\cdot)\big| \calF_t^{\bW^0}\right ](\omega).
\]

\subsubsection*{Revisiting the calculations for the lifted Kalman filter equations}

To illustrate some advantages of using 
the PDEs we derived from the flow transformation method of Lions-Souganidis,
we now provide details for how one can use (\ref{liftedsystemgeneral}) to derive equations such as for $\bar{\eta}$ in (\ref{liftedkalman}).
Although these calculations may appear more tedious than the approach indicated before, 
they have the advantage of shedding some light on the barriers one needs to overcome for the difficult case of non-Gaussian initial condition.

We first can compute, noting that when taking the compensated time derivative, $\gamma_t=z$,
\begin{align*}
\ D_{t,\omega,\gamma}^{y,z}  \bar{\eta}(t,\omega,\gamma,y,z)
=&\ \int_{\R^d}x\, \partial_t\frac{\widetilde{w}(t,x-\sigma^0y,\omega,\gamma) \exp \left (Hx \cdot \tilde{\Theta}^{-1}z \right ) }{\int_{\RR^d} \widetilde{w}(t,\xi-\sigma^0y,\omega,\gamma) 
\exp \left ( H\xi \cdot \tilde{\Theta}^{-1} z  \right ) d\xi}\, dx\\
=&\ \int_{\R^d}(x-\bar{\eta})\, \frac{ \partial_t\widetilde{w}(t,x-\sigma^0y,\omega,\gamma) \exp \left (Hx \cdot \tilde{\Theta}^{-1}z \right )}{\int_{\RR^d} \widetilde{w}(t,\xi-\sigma^0y,\omega,\gamma) 
\exp \left ( H\xi \cdot \tilde{\Theta}^{-1} z  \right ) d\xi} dx.
\end{align*}
To simplify the last expression, we consider the contribution of the first term in the equation
for $\partial_t\tilde{w}$ corresponding to $\frac{1}{2}\tr[\sigma \sigma^\top \partial_{xx}^2 \widetilde{w}(t,x,\omega,\gamma)]$ (see system \eqref{liftedsystemgeneral}), and integrate by parts to obtain:
\begin{align*}
&\ \int_{\R^d}(x-\bar{\eta})\, \frac{  \frac{1}{2}\tr[\sigma \sigma^\top \partial_{xx}^2 \widetilde{w}(t,x-\sigma^0\, y,\omega,\gamma)]  \exp \left (Hx \cdot \tilde{\Theta}^{-1}z \right ) }{\int_{\RR^d} \widetilde{w}(t,\xi-\sigma^0y,\omega,\gamma) 
\exp \left ( H\xi \cdot \tilde{\Theta}^{-1} z  \right ) d\xi} dx \\
=&\ \int_{\R^d} \frac{ \sigma \sigma^\top\, H^\top\, \tilde{\Theta}^{-1}\, z\, \widetilde{w}(t,x-\sigma^0\, y,\omega,\gamma) \exp \left (Hx \cdot \tilde{\Theta}^{-1}z \right ) }{\int_{\RR^d} \widetilde{w}(t,\xi-\sigma^0y,\omega,\gamma) 
\exp \left ( H\xi \cdot \tilde{\Theta}^{-1} z  \right ) d\xi} dx \\
=&\ \sigma \sigma^\top\, H^\top\, \tilde{\Theta}^{-1}\, z.
\end{align*}

For the second term in the equation for $\partial_t\tilde{w}$, corresponding to $(\sigma \sigma^\top H^\top \tilde{\Theta}^{-1} z)\cdot  \partial_x \widetilde{w}(t,x,\omega,\gamma)$,
\begin{align*}
&\ \int_{\R^d}(x-\bar{\eta})\,\frac{ (\sigma \sigma^\top H^\top \tilde{\Theta}^{-1} z)\cdot  \partial_x \widetilde{w}(t,x-\sigma^0\, y,\omega,\gamma)  \exp \left (Hx \cdot \tilde{\Theta}^{-1}z \right ) }{\int_{\RR^d} \widetilde{w}(t,\xi-\sigma^0y,\omega,\gamma) 
\exp \left ( H\xi \cdot \tilde{\Theta}^{-1} z  \right ) d\xi}\, dx\\
=&\ -\sigma \sigma^\top H^\top \tilde{\Theta}^{-1} z,
\end{align*}
which cancels with the first term we computed above.

Next, we turn to the drift terms in the equation for $\partial_t\tilde{w}$, corresponding to\\
 $-\partial_x \cdot \big(\widetilde{w}(t,x-\sigma^0\, y,\omega,\gamma) \,  b(t,x,\hat{m},\hat{\alpha})\big)-\tilde{\omega}(t,x,\omega,\gamma)\hat{b}(t,x+\sigma^0\omega_t,\omega,\gamma)\cdot H^\top\tilde{\Theta}^{-1}\gamma_t$:
\begin{align*}
&\ \int_{\R^d}(x-\bar{\eta})\,  \frac{-\partial_x \cdot \Big(\widetilde{w}(t,x-\sigma^0\, y,\omega,\gamma) \,  b(t,x,\hat{m},\hat{\alpha})\Big)  \exp \left (Hx \cdot \tilde{\Theta}^{-1}z \right )}{\int_{\RR^d} \widetilde{w}(t,\xi-\sigma^0y,\omega,\gamma) 
\exp \left ( H\xi \cdot \tilde{\Theta}^{-1} z  \right ) d\xi} dx\\
&\ -\int_{\R^d}(x-\bar{\eta})\,  \frac{H^\top\tilde{\Theta}^{-1}z \cdot \Big(\widetilde{w}(t,x-\sigma^0\, y,\omega,\gamma) \,  b(t,x,\hat{m},\hat{\alpha})\Big)  \exp \left (Hx \cdot \tilde{\Theta}^{-1}z \right )}{\int_{\RR^d} \widetilde{w}(t,\xi-\sigma^0y,\omega,\gamma) 
\exp \left ( H\xi \cdot \tilde{\Theta}^{-1} z  \right ) d\xi} dx\\
=&\ \int_{\R^d} \frac{b(t,x,\hat{m},\hat{\alpha})\, \tilde{w}(t,x-\sigma^0\, y,\omega,\gamma)\,  \exp \left (Hx \cdot \tilde{\Theta}^{-1}z \right )}{\int_{\RR^d} \widetilde{w}(t,\xi-\sigma^0y,\omega,\gamma) 
\exp \left ( H\xi \cdot \tilde{\Theta}^{-1} z  \right ) d\xi}dx.
\end{align*}

The term corresponding to $\tilde{\omega}(t,x,\omega,\gamma)\frac{1}{2}(H^\top \tilde{\Theta}^{-1}\gamma_t)\cdot \sigma\sigma^\top  H^\top \tilde{\Theta}^{-1}\gamma_t]$ vanishes because there is no $x$ dependence. 

The final term for $\partial_t\tilde{w}$ contributes
\begin{align*}
&\ -\int_{\R^d}(x-\bar{\eta})\, \frac{ \widetilde{w}(t,x-\sigma^0y,\omega,\gamma)\frac{1}{2} x^\top H^\top \tilde{\Theta}^{-1} H\, x    \, \exp \left (Hx \cdot \tilde{\Theta}^{-1}z \right ) }{\int_{\RR^d} \widetilde{w}(t,\xi-\sigma^0y,\omega,\gamma) 
\exp \left ( H\xi \cdot \tilde{\Theta}^{-1} z  \right ) d\xi} dx\\
=&\ -\Pi_t H^\top \tilde{\Theta}^{-1} H \bar{\eta}_t,
\end{align*}
where we use the Gaussian initial condition to reduce to the expression using the covariance matrix $\Pi_t$ by the Gaussian integration by parts formula, i.e., Isserlis' theorem or Wick's probability theorem.

Finally, we arrive at the expression for $D_{t,\omega,\gamma}^{y,z}  \bar{\eta}(t,\omega,\gamma,y,z)$ in (\ref{liftedkalman}),
\begin{align*}
D_{t,\omega,\gamma}^{y,z}  \bar{\eta}(t,\omega,\gamma,y,z) =&\ \int_{\R^d} \frac{b(t,x,\hat{m},\hat{\alpha})\, \tilde{w}(t,x-\sigma^0\, y,\omega,\gamma)\,  \exp \left (Hx \cdot \tilde{\Theta}^{-1}z \right )}{\int_{\RR^d} \widetilde{w}(t,\xi-\sigma^0y,\omega,\gamma) 
\exp \left ( H\xi \cdot \tilde{\Theta}^{-1} z  \right ) d\xi}dx
 -\Pi_t\, H^\top\, \tilde{\Theta}^{-1}\, H\, \bar{\eta}_t\\
 =&\ L_t \bar{\eta}(t,\omega,\gamma,y,z) 
+ (\bar{b}_t - \Lambda^0_t) \bar{\mu}(t,\omega, y).
\end{align*}

For the expressions $\partial_y \bar{\eta}(t,\omega,\gamma,y,z)$ and  $\partial_z \bar{\eta}(t,\omega,\gamma,y,z)$ of (\ref{liftedkalman}) we have
\begin{align*}
 \partial_y \bar{\eta}(t,\omega,\gamma,y,z) 
 =&\   \int_{\R^d} x\, \partial_y\frac{\widetilde{w}(t,x-\sigma^0y,\omega,\gamma) \exp \left (Hx \cdot \tilde{\Theta}^{-1}z \right ) }{\int_{\RR^d} \widetilde{w}(t,\xi-\sigma^0y,\omega,\gamma) 
\exp \left ( H\xi \cdot \tilde{\Theta}^{-1} z  \right ) d\xi}dx\\
=&\ \int_{\R^d} (x-\bar{\eta})\, \frac{-\sigma^0\, \partial_x\widetilde{w}(t,x-\sigma^0y,\omega,\gamma) \exp \left (Hx \cdot \tilde{\Theta}^{-1}z \right ) }{\int_{\RR^d} \widetilde{w}(t,\xi-\sigma^0y,\omega,\gamma) 
\exp \left ( H\xi \cdot \tilde{\Theta}^{-1} z  \right ) d\xi} dx\\
=&\ \sigma^0,
 \end{align*}
 and
 \begin{align*}
 \partial_z \bar{\eta}(t,\omega,\gamma,y,z)
 =&\   \int_{\R^d} x\, \left(\partial_z\frac{\widetilde{w}(t,x-\sigma^0y,\omega,\gamma) \exp \left (Hx \cdot \tilde{\Theta}^{-1}z \right ) }{\int_{\RR^d} \widetilde{w}(t,\xi-\sigma^0y,\omega,\gamma) 
\exp \left ( H\xi \cdot \tilde{\Theta}^{-1} z  \right ) d\xi}\right)^\top dx\\
=&\ \int_{\R^d} (x-\bar{\eta})\, \frac{ x^\top\, H^\top\, \tilde{\Theta}^{-1} \widetilde{w}(t,x-\sigma^0y,\omega,\gamma) \exp \left (Hx \cdot \tilde{\Theta}^{-1}z \right ) }{\int_{\RR^d} \widetilde{w}(t,\xi-\sigma^0y,\omega,\gamma) 
\exp \left ( H\xi \cdot \tilde{\Theta}^{-1} z  \right ) d\xi} dx\\
=&\ \Pi_t\,  H^\top\, \tilde{\Theta}^{-1}.
\end{align*}
This completes the illustration of how to use the general lifted system (\ref{liftedsystemgeneral}) to derive lifted Kalman filter equations of (\ref{liftedkalman}).

\subsection*{Acknowledgments}

The first author would like to thank many people: Daniel Lacker, for helping identify a crucial error in an early reference, which, in order to fix, led to the discovery of the need for the compensator; 
Andrea Cosso and Francesco Russo, for many helpful correspondences; Nizar Touzi, for pointing out useful references;
Nikiforos Mimikos-Stamatopoulos, for regular discussions of technical concepts; and finally and most importantly,
Takis Souganidis, who helped guide 
the lifted functional perspective from its inception as well as provide many critical suggestions for this paper.

\bibliography{LQcompensator}
\bibliographystyle{plain}

\end{document}